\RequirePackage{fix-cm}
\documentclass[smallcondensed]{svjour3}     
\pdfoutput=1

\usepackage[usenames,dvipsnames]{xcolor}  

\smartqed                  
\usepackage{graphicx}
\usepackage{tikz}
\usetikzlibrary{positioning, calc}

\usepackage{mathptmx}      
\usepackage{mathrsfs}

\usepackage[utf8]{inputenc}

\usepackage[english]{babel}

\usepackage{amsthm}
\usepackage{amssymb,amsmath}
\usepackage{amsmath,subfigure,bm}
\usepackage{multicol}
\usepackage{mathtools}
\usepackage{algorithm,algpseudocode}

\algrenewcommand{\algorithmicrequire}{\textbf{INPUT:}}
\algrenewcommand{\algorithmicensure}{\textbf{OUTPUT:}}

\usepackage{color}

\usepackage{verbatim}

\usepackage[numbers, square, comma, sort&compress]{natbib}

\usepackage[colorlinks=true]{hyperref}
\usepackage[noabbrev,capitalize,nameinlink]{cleveref}
\RequirePackage[capitalize,nameinlink]{cleveref}[0.19]

\crefname{section}{section}{sections}
\crefname{subsection}{subsection}{subsections}

\Crefname{figure}{Figure}{Figures}

\crefformat{equation}{\textup{#2(#1)#3}}
\crefrangeformat{equation}{\textup{#3(#1)#4--#5(#2)#6}}
\crefmultiformat{equation}{\textup{#2(#1)#3}}{ and \textup{#2(#1)#3}}
{, \textup{#2(#1)#3}}{, and \textup{#2(#1)#3}}
\crefrangemultiformat{equation}{\textup{#3(#1)#4--#5(#2)#6}}%
{ and \textup{#3(#1)#4--#5(#2)#6}}{, \textup{#3(#1)#4--#5(#2)#6}}{, and \textup{#3(#1)#4--#5(#2)#6}}

\Crefformat{equation}{#2Equation~\textup{(#1)}#3}
\Crefrangeformat{equation}{Equations~\textup{#3(#1)#4--#5(#2)#6}}
\Crefmultiformat{equation}{Equations~\textup{#2(#1)#3}}{ and \textup{#2(#1)#3}}
{, \textup{#2(#1)#3}}{, and \textup{#2(#1)#3}}
\Crefrangemultiformat{equation}{Equations~\textup{#3(#1)#4--#5(#2)#6}}%
{ and \textup{#3(#1)#4--#5(#2)#6}}{, \textup{#3(#1)#4--#5(#2)#6}}{, and \textup{#3(#1)#4--#5(#2)#6}}

\crefdefaultlabelformat{#2\textup{#1}#3}

\makeatletter
\renewcommand*\env@matrix[1][\arraystretch]{%
  \edef\arraystretch{#1}%
  \hskip -\arraycolsep
  \let\@ifnextchar\new@ifnextchar
  \array{*\c@MaxMatrixCols c}}
\makeatother

\include{markup}        

\def\meq{{M_\text{eqn}}}

\def\mdeg{{M_\text{deg}}}
\def\morder{{M_\text{O}}}
\def\morderth{{M_\text{O}^{\text{th}}}}
\def\melems{{M_\text{elem}}}

\def\mpred{{M_\text{P}}}
\def\mcorr{{M_\text{C}}}

\newcommand{\bunderline}[1]{\underline{#1}}
\renewcommand{\vec}[1]{{\bunderline{#1}}}
\newcommand{\mat}[1]{{\bunderline{\bunderline{#1}}}}

\newcommand{\Tm}{{\mathcal T}}
\newcommand{\reals}{\mathbb{R}}

\newcommand{\half}{1/2}

\def\WS{{\mathcal W}}

\theoremstyle{definition}
\theoremstyle{remark}

\numberwithin{equation}{section}

\begin{document}

\title{A Positivity-Preserving Limiting Strategy for Locally-Implicit Lax-Wendroff Discontinuous Galerkin Methods}

\titlerunning{A Positivity-Preserving Limiting Strategy for Locally-Implicit Lax-Wendroff DG}

\author{Camille Felton \and
        Mariana Harris \and
        Caleb Logemann \and
        Stefan Nelson \and
        Ian Pelakh \and
        James A. Rossmanith\footnote{Corresponding author}}
\institute{Camille Felton \at
           University of Wisconsin--Platteville, Department of Mathematics, 435 Gardner Hall, 1 University Plaza,
           Platteville, WI 53818, USA, \email{\url{camillefelton129@gmail.com}} \\
           \and
           Mariana Harris \at
           Instituto Tecnol\'ogico Aut\'onomo de M\'exico,
           R\'io Hondo 1, Progreso Tizap\'an,
           01080 Alvaro Obreg\'on, CDMX, Mexico,
           \email{\url{harris-94@hotmail.com}} \\
           \and
           Caleb Logemann \at
           Iowa State University,
           Department of Mathematics,
           411 Morrill Road,
           Ames, IA 50011-2104, USA,
           \email{\url{logemann@iastate.edu}} \\
           \and
           Stefan Nelson \at
           Minnesota State University Moorhead,
           Department of Mathematics,
           1104 7th Ave South,
           Moorhead, MN 56563, USA,
           \email{\url{nelsonste@mnstate.edu}} \\
           \and
           Ian Pelakh \at
           University of Florida,
           Department of Mathematics,
           1400 Stadium Rd,
           Gainesville, FL 32611, USA
           \email{\url{ipelakh@ufl.edu}} \\
           \and
           James A. Rossmanith \at
           Iowa State University,
           Department of Mathematics,
           411 Morrill Road,
           Ames, IA 50011-2104, USA,
           \email{\url{rossmani@iastate.edu}}
}


\date{Received: date / Accepted: date}

\maketitle

\begin{abstract}

Nonlinear hyperbolic conservation laws admit singular solutions such as shockwaves (discontinuities in conserved variables), rarefaction waves (discontinuities in derivatives), and vacuum states (loss of strong hyperbolicity). When ostensibly high-order numerical methods are applied in such solution regimes, unphysical oscillations present themselves that can lead to large errors and a breakdown of the numerical simulation. In this work we develop a new Lax-Wendroff discontinuous Galerkin (LxW-DG) method with a limiting strategy that keeps the solution non-oscillatory and positivity-preserving for relevant variables, such as height in the shallow water equations and density and pressure in the compressible Euler equations. The proposed LxW-DG scheme updates the solution over each time-step with a locally-implicit predictor followed by an explicit corrector. The locally-implicit prediction phase is formulated in terms of primitive variables, which greatly simplifies the solver. The resulting system of nonlinear algebraic equations are approximately solved via a Picard iteration, where the number of iterations is equal to the order of accuracy of the method. The correction phase is an explicit evaluation formulated in terms of conservative variables in order to guarantee numerical conservation. In order to achieve full positivity-preservation, limiting is required in both the prediction and correction steps. The resulting scheme is applied to several standard test cases for the shallow water and compressible Euler equations. All of the presented examples are written in a freely available open-source Python code.

          

\end{abstract}




\section{Introduction}
\label{sec:intro}

Hyperbolic conservation laws are systems of  partial differential equations
used to model a variety of phenomena characterized by waves propagating at
finite speeds; examples include the shallow water (gravity waves), compressible Euler (sound waves), Maxwell (light waves), and Einstein (gravitational waves) equations.
An important feature of hyperbolic conservation laws is that initially smooth solutions may become singular in finite time. Examples of such singularities include (1) {\it shockwaves}, which are discontinuities in the solution,
(2) {\it rarefactions}, which contain discontinuities in the derivatives of the solution, and (3) {\it vacuum states}, which are solutions in regions of solution space where the equation fails to be strongly hyperbolic.
For all of these solutions, care must be taken to appropriately define the notion of a {\it weak solution}, and additional criteria must be introduced to select a unique weak solution via
appropriate entropy conditions and vanishing viscosity solutions (e.g., see Lax \cite{article:Lax1973}).

The formation of singularities in finite time cause standard high-order methods -- which are almost always based on some form of polynomial interpolation -- to exhibit unphysical oscillations (i.e., Gibbs phenomena). 
These unphysical oscillations can lead to a loss of numerical stability, which in turn can lead to the complete breakdown of the numerical computation. Even in the event that the instabilities do not lead to a full breakdown of the computation, they are often characterized by large numerical errors.

One potential remedy for the unphysical solutions produced by high-order methods is to introduce a post-processing step known as a {\it limiter}. The idea is that when and where the solution is smooth, the limiter should do nothing, but, when unphysical oscillations or excursion of the solutions outside of the region of hyperbolicity occur,  the limiter should damp the high-order correction terms in order to remove the unwanted behavior. 
The early work of limiters for high-order schemes applied to hyperbolic conservation began
 in the early 1970s with works such as Harten and Zwas \cite{harten1972self}, 
Kolgan \cite{article:Kol72,article:vanLeer2011}, van Leer
\cite{article:vanLeer1,article:VanLeer74}, and Boris and Book \cite{boris1973flux}.
In the more than 40 intervening years, limiters have been developed and generalized for a host of equations and methods, including high-resolution finite volume schemes (e.g., see Chapter 6 of LeVeque \cite{book:Le02}), weighted essentially non-oscillatory (WENO) schemes (e.g., see review article by Shu
\cite{article:ShuWENO2009}), and discontinuous Galerkin (DG) schemes (e.g., 
Krivodonova \cite{article:Kriv07}, Persson and Peraire
\cite{article:PerPer06},  Qiu and Shu \cite{article:QiuShu04},
and Zhang and Shu \cite{article:ZhangShu11}).

The focus of the current paper is on the discontinuous Galerkin (DG) method, which
was first introduced by Reed and Hill \cite{article:ReedHill73} for neutron transport, and then fully developed for time-dependent hyperbolic conservation laws in a series of papers by  Cockburn,  Shu, and collaborators (see \cite{article:CoShu98} and references therein for details). 
DG is a particular flavor of the finite element method that is based on piecewise continuous basis functions (almost always polynomials) that are discontinuous across element faces. These discontinuities have two important consequences when applied to spatial discretizations of \cref{eqn:conslaw}: (1) the associated mass-matrix is block diagonal (the size of these blocks are the number of degrees of freedom on each element), and (2) the discontinuities create a small amount of artificial dissipation that helps stabilize the numerical method (in contrast to continuous Galerkin schemes, which require additional stabilization terms).

In this work we develop a novel variant of the Lax-Wendroff DG (LxW-DG)
 method \cite{article:Qiu05}. In particular, the starting point of this method is the LxW-DG formulation of
Gassner et al. \cite{article:GasDumHinMun2011}, in which every time-step is comprised of two distinct phases:
\begin{description}
\item[{\bf Prediction phase.}] This phase is akin to a single-step of a block-Jacobi iteration of a fully implicit spacetime DG approach \cite{article:KlaVegVen2006,article:Sudirham2006}, and is
completely local on each space-time element.
\item[{\bf Correction phase.}]  This phase is an Euler-like step used to advance the solution from the old to the new time, and requires the computation of temporal and spatiotemporal integrals of the predicted solution.
\end{description}
A novel feature of the proposed scheme is that the prediction step is done entirely using
primitive variables, which both simplifies the prediction step and subsequently 
allows a simple introduction of limiters.
We formulate the resulting scheme so that it can be made arbitrarily high-order,
but must, as always, confront the challenge that the scheme may breakdown at shocks,
rarefactions, and vacuum states. In order to overcome this difficulty, we introduce
four sets of limiters:
\begin{description}
\item[{\bf Prediction step positivity limiter.}] Using ideas similar to 
the celebrated Zhang and Shu \cite{article:ZhangShu11} limiter, 
we develop a completely local limiter that minimally damps the high-order corrections to the primitive variables in order to keep the numerical solution inside the region of hyperbolicity.
\item[{\bf Correction step positivity limiter I.}] Following the limiter developed by
Moe et al. \cite{article:MoRoSe17}, we introduce a limiter that blends the high-order time-averaged numerical fluxes used to update the cell averages with a low-order flux in such a way to obtain high-order cell averages that are inside the region of hyperbolicity.
\item[{\bf Correction step positivity limiter II.}] Similar to what was done in the prediction step, we use the Zhang and Shu \cite{article:ZhangShu11} limiter to minimally damp the high-order corrections to the conserved variables in order to keep the numerical solution inside the region of hyperbolicity.
\item[{\bf Correction step unphysical oscillation limiter.}] We develop a limiter based on the 
hierarchical minmod limiter of Krivodonova \cite{article:Kriv07} that is able to minimally damp high-order corrections to remove unphysical oscillations due to the
Gibbs phenomenon at shocks and rarefactions.
\end{description}
The resulting scheme is applied to several standard test cases for the shallow water and compressible Euler equations. All of the presented examples are written in a freely available open-source Python code.

The remainder of this paper is structured as follows.
After reviewing the specific hyperbolic conservation laws considered in this work in \cref{sec:equations}, we explain the full details of both the prediction and correction
steps in \cref{sec:LxW-DG}. The limiters are fully described in \cref{sec:limiters}; for each limiter we provide detailed pseudo-code algorithms. We also present a pseudo-code for a full time-step of the proposed LxW-DG scheme. The resulting scheme is implemented in a Python code that we are making freely available; a brief description of this code is presented in
\cref{sec:pycode}. In \cref{sec:numerical_examples} we apply the proposed algorithm to a series of numerical tests for the Burgers, shallow water, and compressible Euler equations. We clearly demonstrate the efficacy of both the non-oscillatory and positivity-preserving limiters.
We conclude in \cref{conclusions}.


\section{Model equations}
\label{sec:equations}
In this section we briefly review the mathematical properties of hyperbolic conservation laws
(\cref{subsec:hyp_cons_laws}) and the three equations of interest in this work: (1) Burgers equation (\cref{subsec:burgers}),
(2) the shallow water equations (\cref{subsec:shllw}), and (3) the compressible Euler equations
(\cref{subsec:euler}). For a full treatment of these
equations see for example the textbooks of LeVeque \cite{book:Le02}.

\subsection{\it Hyperbolic conservation laws}
\label{subsec:hyp_cons_laws}
We consider a class of 
partial differential equations in one spatial dimension known as {\it conservation laws}, which
can be written in the form:
\begin{equation}
\label{eqn:conslaw}
\vec{q}_{,t}+\vec{f}\left(\vec{q}\right)_{,x} = \vec{0},
\end{equation}
where $t \in \reals_{\ge 0}$ is time, $x \in \reals$ is the one-dimensional spatial coordinate, $\vec{q}(t,x): \reals^{+} \times \reals \mapsto \reals^\meq$ is the vector of $\meq$ conserved variables, which may include things such as mass, momentum, and energy, and $\vec{f}\left(\vec{q}(t,x)\right): \reals^\meq \mapsto \reals^\meq$ is the flux function.

We refer to \cref{eqn:conslaw} as the equation written in {\it conservative form}. This form is fundamental since it is directly connected to the {\it integral} conservation law:
\begin{equation}
\frac{d}{dt} \int_{x_1}^{x_2} \vec{q}(t,x) \, dx = \vec{f}\left(\vec{q}(t,x_1)\right) - 
\vec{f}\left(\vec{q}(t,x_2)\right),
\end{equation}
where $x_1$ and $x_2$ are arbitrary, 
which states that the total amount of $q$ on the domain $[x_1,x_2]$ can only be modified
by a flux at $x=x_1$ into the domain and a flux $x=x_2$ out of the domain. The
above integral form does not require smoothness on $q$, and is necessary to properly define the notion of {\it weak solutions} of \cref{eqn:conslaw}.
If the solution is smooth, we can use the chain rule to put the equation in {\it quasilinear form}:
\begin{equation}
\label{eqn:quasilinear_cons}
\vec{q}_{,t}+\mat{A}\left(\vec{q}\right) \vec{q}_{,x} = \vec{0}, \qquad \mat{A}\left(\vec{q} \right) = \vec{f}\left(\vec{q}\right)_{,\vec{q}},
\end{equation}
where $\mat{A}\left(\vec{q}\right)$ is the {\it flux Jacobian}.

It is sometimes useful to consider writing the quasilinear equation in terms of variables
other than the conservative variables; typical examples include the {\it primitive} and {\it entropy} variables.
For simplicity, we will just refer to these ``other'' variables as the primitive variables and denote them by the
symbol $\alpha$. More concretely, these variables are related to the conservative variables via the chain rule:
\begin{equation}
\vec{q}_{,t} = \vec{q}_{,\vec{\alpha}} \, \vec{\alpha}_{,t} \quad \text{and} \quad
\vec{q}_{,x} = \vec{q}_{,\vec{\alpha}} \, \vec{\alpha}_{,x},
\end{equation}
which we can then use to rewrite \cref{eqn:quasilinear_cons} as
\begin{equation}
\label{eqn:quasilinear_prim}
\vec{\alpha}_{,t}+\mat{B}\left(\vec{\alpha}\right) \vec{\alpha}_{,x} = \vec{0}, \qquad
\mat{B} = \vec{q}_{,\vec{\alpha}}^{-1} \, \mat{A} \, \vec{q}_{,\vec{\alpha}}.
\end{equation}
Note that the matrix $\mat{B}$ and the flux Jacobian, $\mat{A}$, are similar matrices, which means that
they have the same eigenvalues.

In this work we consider a subclass of conservation laws of the form \cref{eqn:conslaw} that are
{\it hyperbolic}. Hyperbolicity is connected to the concept of {\it causality}; and therefore, hyperbolic conservation laws model phenomena characterized by waves propagating at
finite speeds; examples include the shallow water (gravity waves), compressible Euler (sound waves), Maxwell (light waves), and Einstein (gravitational waves) equations.
Mathematically, {\it hyperbolicity} is defined as follows.

\begin{definition}
Conservation law \cref{eqn:conslaw} is {\it hyperbolic} on the convex set
$S \subset \reals^\meq$ if the {flux Jacobian} \cref{eqn:quasilinear_cons}
 is diagonalizable with only real eigenvalues for all $\vec{q} \in S \subset \reals^\meq$.
\end{definition}

\subsection{\it Burgers equation}
\label{subsec:burgers}
The inviscid Burgers equation on the real line $x \in \reals$ in conservation form can be written as follows:
\begin{equation}
\label{eqn:burgers}
q_{,t} + \left( \frac{1}{2} q^2 \right)_{,x} = 0, \qquad q(t=0,x) = q_0(x).
\end{equation}
For smooth solutions, we can put this equation in quasilinear form:
\begin{equation}
q_{,t} + q q_{,x} = 0,
\end{equation}
which implies that the flux Jacobian is $A(q) = q$. This equation is hyperbolic on the set
\begin{equation}
S = \Bigl\{ q \in \reals \Bigr\}.
\end{equation}
The system supports a single wave, $\lambda_1 = q$, which happens to be {\it genuinely nonlinear}:
\begin{equation}
\lambda_{1,q} = 1 \quad \Longrightarrow \quad \lambda_{1,q} \ne 0 \quad \forall q \in S.
\end{equation}

If the initial condition, $q_0$, is smooth on $x \in \reals$, then there exist a time interval, $[0,t_{\text{shock}})$,
over which the solution remains smooth and can be expressed as
\begin{equation}
\label{eqn:Burgers_char1}
 q(t,x) = q_0\left( \xi(t,x) \right), \quad x = \xi + t q_0(\xi),
\end{equation}
where $\xi(t,x)$ are the {\it characteristics}, which in the case of the Burgers equation are straight lines
in the $xt$-plane with slopes determined by the initial condition $q_0$. 
The breakdown of this solution, i.e., the shock-formation, is the first value of $t$ when the
mapping between $x$ and $\xi$, encoded by the characteristics, is singular:
\begin{equation}
\frac{dx}{d\xi} = 1 + t q'_0(\xi) = 0 \quad \Longrightarrow \quad
t_{\text{shock}} =  \left[ \max_{\xi \in \reals} \left\{ -q'_0(\xi) \right\} \right]^{-1}.
\end{equation}

\subsection{\it Shallow water equations}
\label{subsec:shllw}
The shallow water equations model the dynamics of a thin, constant density, fluid layer that is in vertical hydrostatic balance. In 1D, this system can be written in conservative form \cref{eqn:conslaw}, with
$\meq=2$ and
\begin{equation}
\label{eqn:shllw}
\vec{q} = \left( h, \, hu \right), \quad 
\vec{f}\left(\vec{q}\right) = \left( hu, \, hu^2 + \frac{1}{2} g h^2 \right),
 \quad \text{and} \quad
 \mat{A}\left(\vec{q}\right) = \begin{bmatrix}[1.5] 0 & \quad 1 \\ gh-u^2 & \quad 2u \end{bmatrix},
\end{equation}
where $h$ is the thickness of the fluid layer, $u$ is the vertically integrated horizontal fluid velocity,
$hu$ is the macroscopic horizontal momentum density, and $g$ is the gravitational constant.

The eigenvalues of the flux Jacobian are
\begin{equation}
\lambda_1 = u-\sqrt{gh}\quad \text{and} \quad
\lambda_2 = u+\sqrt{gh},
\end{equation}
while the matrices of right and left eigenvectors can be written as
\begin{equation}
\label{eqn:shllw_eignevectors}
\mat{R} =
\begin{bmatrix}[1.5]
1 &  \quad 1 \\
u-\sqrt{gh} & \quad u+\sqrt{gh}
\end{bmatrix} 
\quad \text{and} \quad
\mat{L} = \mat{R}^{-1} = \frac{1}{2\sqrt{gh}}
\begin{bmatrix}[1.5]
\sqrt{gh}+u &  \quad -1 \\
\sqrt{gh}-u &  \quad  \hspace{2.5mm} 1
\end{bmatrix}.
\end{equation}
The primitive variables, the matrix $\mat{B}$ in the primitive quasilinear system
\cref{eqn:quasilinear_prim}, and the matrix of right eigenvectors of
$\mat{B}$, are given by:
\begin{equation}
\label{eqn:shllw_prim_matrix}
\vec{\alpha} = \left( h, \, u \right), \quad
\mat{B} = \begin{bmatrix}[1.5]
u & \quad h \\ g & \quad u
\end{bmatrix},
 \quad \text{and} \quad
 \mat{R_B} = \begin{bmatrix}[1.5]
-\sqrt{h} & \quad \sqrt{h} \\ \sqrt{g} & \quad \sqrt{g}
\end{bmatrix}.
\end{equation} 

We infer from the above information that the shallow water equations are hyperbolic on the
convex set:
\begin{equation}
  S = \Bigl\{ \vec{q} = (h, \, hu) \in \reals^2: \, h > 0 \Bigr\}.
\end{equation}
Note that the convexity of $S$ follows from the fact that the Hessian of $h$ with respect to $\vec{q}$
is negative semi-definite for all $\vec{q} \in S$:
\begin{equation}
h_{,\vec{q},\vec{q}} = \left( 1, \, 0 \right)_{,\vec{q}}= \mat{0}.
\end{equation}
We also note that the two {\it wave families} represented by $\lambda_1$ and $\lambda_2$ are
both referred to as {\it gravity waves}; a quick calculation shows that these waves are
genuinely nonlinear:
\begin{equation}
\lambda_{1,\vec{\alpha}} \cdot \vec{r_B}_1 = \lambda_{2,\vec{\alpha}} \cdot \vec{r_B}_2
 = \frac{3}{2} \sqrt{g} \quad \Longrightarrow \quad \lambda_{1,\vec{\alpha}} \cdot \vec{r_B}_1 = \lambda_{2,\vec{\alpha}} \cdot \vec{r_B}_2 \ne 0 \quad \forall \vec{q} \in S,
\end{equation}
where $\vec{r_B}_1$ and $\vec{r_B}_2$ are the two columns of the matrix of right eigenvectors
$\mat{R_B}$ given by \cref{eqn:shllw_prim_matrix}.
The fact that the gravity waves are genuinely nonlinear
means that each wave family can form a shock or rarefaction.

\subsection{\it Compressible Euler equations}
\label{subsec:euler}
The compressible Euler equations model the dynamics of a monatomic gas in thermodynamic equilibrium.
In 1D, this system can be written in conservative form \cref{eqn:conslaw}, with
$\meq=3$ and
\begin{equation}
\label{eqn:euler}
\begin{gathered}
\vec{q} = \left( \rho, \, \rho u, \, {\mathcal E} \right), \quad 
\vec{f}\left(\vec{q}\right) = \left( \rho u, \, \rho u^2 + p, \, u\left({\mathcal E}+p \right) \right),
 \quad \text{and} \\
\mat{A} = 
\begin{bmatrix}[1.5]
0 & \quad 1 & \quad 0 \\
\frac{1}{2} \left(\gamma-3\right) u^2 & \quad (3 - \gamma) u & \quad \gamma - 1 \\
\frac{1}{2} u^3 (\gamma - 2) + \frac{\gamma p u}{(1 - \gamma) \rho} & \quad
\frac{1}{2} u^2 (3 - 2 \gamma) - \frac{\gamma p}{(1-\gamma)\rho} & \quad \gamma u
\end{bmatrix},
\end{gathered}
\end{equation}
where $\rho$ is the mass density, $u$ is the fluid velocity, $p$ is the pressure, $\rho u$ is the
momentum density, $\gamma>1$ is the {\it specific heat ratio}, sometimes referred to as
the {\it adiabatic index},
and ${\mathcal E}$ is the energy density:
\begin{equation}
{\mathcal E} = \frac{p}{\gamma - 1} + \frac{1}{2} \rho u^2.
\end{equation}

The eigenvalues of the flux Jacobian are
\begin{equation}
\lambda_1 = u-c, \quad
\lambda_2 = u, \quad \text{and} \quad
\lambda_3 = u+c,
\end{equation}
where the sound speed is $c = \sqrt{{\gamma p}/{\rho}}$,
while the matrices of right and left eigenvectors can be written as
\begin{gather}
\label{eqn:euler_right_eignevectors}
\mat{R} =
\begin{bmatrix}[1.5]
1 &  \quad 1 & \quad 1 \\
u-c & \quad u & \quad u+c \\
\frac{p\gamma}{(\gamma-1) \rho} + \frac{1}{2} u (u - 2 c) & \quad
\frac{1}{2} u^2 & \quad
\frac{p\gamma}{(\gamma-1) \rho} + \frac{1}{2} u (u + 2 c)
\end{bmatrix} 
\quad \text{and} \\
\mat{L} = \mat{R}^{-1} = \frac{1}{4 c \gamma p}
\begin{bmatrix}[1.5]
c \rho u^2 (\gamma - 1) + 2 \gamma p u & \quad 
2 c \rho u (1 - \gamma) - 2 \gamma p & \quad 
2 c \rho (\gamma - 1) \\
2 c ((1 - \gamma) \rho u^2 + 2 \gamma p) & \quad 
4 c \rho u (\gamma - 1) & \quad 
4 c \rho (1 - \gamma) \\
c \rho u^2 (\gamma - 1) - 2 p u \gamma & \quad 
2 p \gamma + 2 c \rho u (1 - \gamma) & \quad
2 c \rho (\gamma - 1)
\end{bmatrix}.
\end{gather}
The primitive variables and matrix $\mat{B}$ in the primitive quasilinear system
\cref{eqn:quasilinear_prim}, as well as the matrix of right eigenvectors of
$\mat{B}$, are given by:
\begin{equation}
\label{eqn:euler_prim_matrix}
\vec{\alpha} = \left( \rho, \, u, \, p \right), \quad
\mat{B} = \begin{bmatrix}[1.5]
u & \quad \rho & \quad 0 \\ 0 & \quad u & \quad \frac{1}{\rho} \\
0 & \quad \gamma p & \quad u
\end{bmatrix},
 \quad \text{and} \quad
 \mat{R_B} = \begin{bmatrix}[1.5]
 \frac{\rho}{c} & \quad 1 & \quad \frac{\rho}{c} \\
  -1 & \quad 0 & \quad 1 \\
  \rho c & \quad 0 & \quad \rho c
 \end{bmatrix}.
\end{equation} 

We infer from the above information that the compressible Euler equations are hyperbolic on the
convex set:
\begin{equation}
  S = \Bigl\{ \vec{q} = (\rho, \, \rho u, \, {\mathcal E}) \in \reals^3: \quad \gamma > 1, \quad
  \rho > 0, \quad
  \text{and} \quad p>0 \Bigr\}.
\end{equation}
Note that the convexity of $S$ follows from the fact that the Hessian of $\rho$ with respect to $\vec{q}$
and the Hessian of $p$ with respect to $\vec{q}$ are negative semi-definite for all $\vec{q} \in S$:
\begin{gather}
\rho_{,\vec{q},\vec{q}} = \left( 1, \, 0, \, 0 \right)_{,\vec{q}}= \mat{0}, \\
p_{,\vec{q},\vec{q}} = \frac{(\gamma-1)}{\rho} \begin{bmatrix}[1.5] -u^2 & \quad u & \quad 0 \\
u & \quad -1 & \quad 0 \\ 0 & \quad 0 & \quad 0\end{bmatrix} \quad \Longrightarrow \quad
\lambda\left( p_{,\vec{q},\vec{q}}\right) = 0, \, 0, \, \frac{\left(1 - \gamma \right) \left(1 + u^2
\right)}{\rho},
\end{gather}
where we note that $\lambda_k(p_{,\vec{q},\vec{q}}) \le 0$ for all $k=1,2,3$ if $\gamma>1$ and $\rho>0$.

We also note that the wave families represented by $\lambda_1$ and $\lambda_3$ are
both referred to as {\it sound waves}, while $\lambda_2$ is referred to as the {\it contact discontinuity};
 a quick calculation shows that the sound waves are
genuinely nonlinear, while the contact is linearly degenerate:
\begin{alignat}{2}
\lambda_{1,\vec{\alpha}} \cdot \vec{r_B}_1 &= -\frac{1}{2} (\gamma+1) 
   \quad & \Longrightarrow \quad \lambda_{1,\vec{\alpha}} \cdot \vec{r_B}_1 \ne 0 \quad \forall \vec{q} \in S, \\
\lambda_{2,\vec{\alpha}} \cdot \vec{r_B}_2 &= 0 \quad & \Longrightarrow \quad 
    \lambda_{2,\vec{\alpha}} \cdot \vec{r_B}_2 = 0 \quad \forall \vec{q} \in S, \\
\lambda_{3,\vec{\alpha}} \cdot \vec{r_B}_3 &= +\frac{1}{2} (\gamma+1) 
   \quad & \Longrightarrow \quad \lambda_{3,\vec{\alpha}} \cdot \vec{r_B}_3 \ne 0 \quad \forall \vec{q} \in S,
\end{alignat}
where $\vec{r_B}_1$, $\vec{r_B}_2$, and $\vec{r_B}_3$ are the three columns of the matrix of right eigenvectors
$\mat{R_B}$ given by \cref{eqn:euler_prim_matrix}.
The fact that the sound waves are genuinely nonlinear
means that each wave family can form a shock or rarefaction. The fact that the contact discontinuity is linearly degenerate means that this wave does not undergo compression or rarefaction.


\section{Locally-implicit Lax-Wendroff discontinuous Galerkin}
\label{sec:LxW-DG}
The Lax-Wendroff method \cite{article:LxW60} is a fully discrete numerical method for
hyperbolic conservation laws (e.g., of the form \cref{eqn:conslaw}) based on
the the Cauchy-Kovalevskaya 
\cite{article:Kovaleskaya1875} procedure to convert temporal derivatives into spatial derivatives.
For example, in the case of conservation law \cref{eqn:conslaw},
we begin with a Taylor series in time:
\begin{equation}
 \vec{q}(t+\Delta t, x) = \vec{q}(t,x) + \Delta t \vec{q}_{,t}(t,x) + \frac{1}{2} \Delta t^2 \vec{q}_{,t,t}(t,x) +
 \ldots,
\end{equation}
and then replace all time derivatives by spatial derivatives:
\begin{equation}
\label{eqn:LxW_time_derivs}
 \vec{q}_{,t} = - \vec{f}\left(\vec{q}\right)_{,x}, \quad
 \vec{q}_{,t,t} = - \vec{f}\left(\vec{q}\right)_{,t,x} = -\left[ \mat{f'}\left(\vec{q}\right) \, \vec{q}_{,t} \right]_{,x} = 
 \left[ \mat{f'}\left(\vec{q}\right) \vec{f}\left(\vec{q}\right)_{,x} \right]_{,x}, \quad \ldots,
\end{equation}
which results in the following:
\begin{equation}
\vec{q}(t+\Delta t ,x) = \vec{q} - \Delta t \vec{f}\left(\vec{q}\right)_{,x} + \frac{1}{2} \Delta t^2 \left[ \mat{f'}\left(\vec{q}\right) \vec{f}\left(\vec{q}\right)_{,x} \right]_{,x}
+ \ldots,
\end{equation}
where on the right-hand side we have suppressed the evaluation at $(t,x)$.
The final step is to truncate the Taylor series at some finite number of terms,
and then replace all spatial derivatives by some discrete spatial derivative operators.
The above Lax-Wendroff formalism \cite{article:LxW60} has been used in conjunction with
with a variety of spatial discretizations, including finite volume  \cite{article:Le97},
weighted essentially non-oscillatory (WENO) \cite{article:TitaToro02}, and discontinuous Galerkin \cite{article:Qiu05} operators.

In this work, we are concerned with the discontinuous Galerkin version 
of Lax-Wendroff \cite{article:Qiu05}; and in particular, we make use
of the reformulation of Gassner et al. \cite{article:GasDumHinMun2011} of the Lax-Wendroff discontinuous Galerkin (LxW-DG) scheme
in terms of a locally-implicit prediction step, followed by an explicit correction step.
The key advantage of this formulation is that we do not need to explicitly compute the
partial derivatives as shown in \cref{eqn:LxW_time_derivs}; and instead, the locally-implicit
solver automatically produces discrete versions of these derivatives.
The next challenge is to efficiently solve the nonlinear algebraic equations that arise from the 
locally-implicit prediction step; we solve these equations by again following Gassner et al. \cite{article:GasDumHinMun2011} and making use of a Picard fixed point iteration. 
One key difference in this work is that we formulate the prediction step in terms of primitive
variables, which pays dividends when we develop limiters (see \cref{sec:limiters}).
We will follow the notational conventions of Guthrey and Rossmanith \cite{article:GuRo2017}
developed for locally-implicit and regionally-implicit LxW-DG schemes.

In the remainder of this section we develop the basic locally-implicit Lax-Wendroff 
discontinuous Galerkin (LxW-DG) scheme. We begin with a quick review of 
the DG spatial discretization \cref{sec:dg-fem-space}, 
followed by a detailed description of the prediction (\cref{subsec:prediction_step}) and
correction (\cref{subsec:corr_step}) phases in one time-step of the proposed locally-implicit LxW-DG method.
The discussion in this section is on the scheme {\it without} limiters;
limiters will receive our full attention in \cref{sec:limiters}.

\subsection{\it DG-FEM spatial discretization}
\label{sec:dg-fem-space}
We discretize system \cref{eqn:conslaw} in space via the discontinuous Galerkin (DG) method, which
was first introduced by Reed and Hill \cite{article:ReedHill73} for neutron transport, and then fully developed for time-dependent hyperbolic conservation laws in a series of papers by Bernardo Cockburn, Chi-Wang Shu, and collaborators (see \cite{article:CoShu98} and references therein for details). 

The computational domain
is a finite interval on the real line: $\Omega = [x_{\text{low}},x_{\text{high}}] \subset \reals$.
Let ${\mathbb P}\left(\mdeg \right)$ denote the set of polynomials from $\reals$ to $\reals$ 
with maximal polynomial degree $\mdeg$.
On the mesh of $\melems$ elements we define the {\it broken} finite element space:
\begin{equation}
\label{eqn:broken_space}
    \WS^{h} := \left\{ \vec{w}^h \in \left[ L^{\infty}(\Omega) \right]^{\meq}: \,
    \vec{w}^h \bigl|_{\Tm_i} \in \left[ {\mathbb P} \left(\mdeg \right) \right]^{\meq} \, \, \forall \Tm_i \right\},
\end{equation}
where $h:=\Delta x = (x_{\text{high}}-x_{\text{low}})/\melems$ is the uniform grid spacing, 
$\meq$ is the number of conserved variables, and $\mdeg$ is the 
maximal polynomial degree in the finite element representation.
The above expression means that $\vec{w} \in \WS^h$ has
$\meq$ components, each of which when restricted to some element $\Tm_i$
is a polynomial in ${\mathbb P}\left(\mdeg \right)$, and no continuity is assumed
across element faces. 

The computational mesh is comprised of $\melems$ elements, which we denote by
\begin{equation}
{\Tm_i} = \left[ x_{i} - \frac{\Delta x}{2},  x_{i} + \frac{\Delta x}{2} \right] \quad \text{for} \quad
i=1,\ldots,\melems.
\end{equation}
For convenience, we define the canonical variable, $\xi$, which on each element is related
the physical coordinate, $x$, as follows:
\begin{equation}
\label{eqn:local_xi}
x = x_{i} + \left( \frac{\Delta x}{2} \right) \xi, \qquad \xi \in [-1,1].
\end{equation}
Using the canonical variable, we define on each element the following Legendre polynomial
basis:
\begin{equation}
\label{eqn:phi_basis_1d}
\vec{\Phi} = \left( 1, \, \sqrt{3} \xi, \, \frac{\sqrt{5}}{2} \left( 3 \xi^2 - 1 \right), \, \frac{\sqrt{7}}{2}\left(5 \xi^3 - 3 \xi \right), \, 
\frac{\sqrt{9}}{8} \left(35 \xi^4 - 30 \xi^2 + 3 \right), \, \cdots \right),
\end{equation}
which also satisfies the following three-term recurrence relationship
for $k\ge3$:
\begin{equation}
\Phi_k(\xi) = \left( \frac{\sqrt{(2k-3)(2k-1)}}{(k-1)} \right) \xi \Phi_{k-1}(\xi) 
       - \left( \frac{(k-2)\sqrt{2k-1}}{(k-1)\sqrt{2k-5}} \right) \Phi_{k-2}(\xi),
\end{equation}
where $\Phi_1 = 1$ and $\Phi_2 = \sqrt{3} \xi$.

The approximate solution on each element at time $t=t^n$:
\begin{equation}
\label{eqn:corrector_ansatz}
   \vec{q}^h\left(t^n, \, x_i +  \frac{\Delta x}{2} \, \xi \right) := \vec{\Phi}\left(\xi \right)^T \mat{Q_i^n},
   \quad \text{for} \quad \xi \in [-1,1],
\end{equation}
where
\begin{align}
\vec{\Phi}\left(\xi \right): [-1,1] \mapsto \reals^{\mcorr} \qquad \text{and} \qquad
\mat{Q_i^n} \in \reals^{\mcorr \times \meq}.
\end{align}
Note that we are denoting the number of basis functions by $\mcorr$ to signify that this refers to the number of basis functions in the {\it correction step}.
If the exact solution, $\vec{q}(t^n,x)$, is known {\it a priori}, then we can compute the corresponding coefficients, ${Q}$, via $L^2$-projection:
\begin{equation}
\mat{Q^n_i} = \frac{1}{2} \int_{-1}^{1} \vec{\Phi}(\xi) \left[ \vec{q}\left(t^n, x_i + \frac{\Delta x}{2} \, \xi \right) \right]^T \, d\xi 
\approx \frac{1}{2} \sum_{a=1}^{\morder} \omega_a \, \vec{\Phi}(\mu_a) \left[ \vec{q}\left(t^n, x_i + \frac{\Delta x}{2} \mu_a \right) \right]^T,
\end{equation}
where $\morder = \mcorr$ is the maximum achievable order of accuracy, and $\omega_a$ and $\mu_a$ for $a=1,\ldots,\morder$ are the weights and abscissas of the $\morder$-point Gauss-Legendre quadrature rule.
In practice, the only solution that is known {\it a priori} is the initial condition: $\vec{q}(0,x)$.
For all subsequent time, the coefficients ${Q}$ must be computed by a numerical procedure, which in
in this work will be the locally-implicit Lax-Wendroff scheme described in detail below.

\subsection{\it Prediction step}
\label{subsec:prediction_step}
The numerical update of the proposed scheme is divided into two distinct parts: (1) the prediction step and (2)
the correction step. In the prediction step we will {\it not} enforce consistency of the numerical method
with the underlying conservation laws \cref{eqn:conslaw}; and therefore, we have a significant amount of freedom
in how this portion of the update can be accomplished. In particular, one freedom which we will exercise
is the choice of variables used in the prediction step (e.g., conservative, primitive, or entropy variables). For simplicity of discussion, we will simply refer to the choice of variables as the {\it primitive variables}. We denote these variables by $\alpha$ and assume (without loss of generality) that they satisfy quasilinear equation \cref{eqn:quasilinear_prim}.

The prediction step is entirely local on each element; and therefore, without loss of generality, we focus our attention on element $\Tm_i$ over time interval $[t^n, t^{n+1}]$, where $t^{n+1} = t^n + \Delta t$. On this
{\it space-time} element, we introduce the local spatial variable $\xi$ as defined by \cref{eqn:local_xi}, the local temporal variable
$\tau$ as defined by
\begin{equation}
t = t^{n}+\frac{\Delta t}{2} (1+\tau), \qquad \tau \in [-1,1],
\end{equation}
and rewrite \cref{eqn:quasilinear_prim} as follows:
\begin{equation}
\label{eqn:alpha_theta}
\vec{\alpha}_{,\tau} = \vec{\Theta}\left(\vec{\alpha} \right) :=  - \nu \mat{B}\left(\vec{\alpha} \right) \, \vec{\alpha}_{,\xi}, \qquad \nu:=\frac{\Delta t}{\Delta x}.
\end{equation}
This equation can also be written in component form as
\begin{equation}
\label{eqn:alpha_theta_comp}
{\alpha}_{m,\tau} = {\Theta}_m\left(\vec{\alpha} \right) =
- \nu \sum_{k=1}^{\meq} {B}_{mk} \left(\vec{\alpha} \right) \, {\alpha}_{k,\xi} \quad \text{for} \quad m=1,\ldots,\meq.
\end{equation}
We introduce a space-time Legendre basis on each element:
\begin{equation}
\label{eqn:spacetime_basis}
\vec{\Psi}(\tau, \xi): [-1,1]^2 \mapsto \reals^{\mpred}, \quad \mpred = \frac{\mcorr(\mcorr+1)}{2}, \quad {\Psi}_{\ell}(\tau, \xi) = {\Phi}_{\ell_\tau}(\tau) \, {\Phi}_{\ell_\xi}(\xi),
\end{equation}
which is orthonormal on $[-1,1]^2$:
\begin{equation}
\label{eqn:psi_orthonormality}
\frac{1}{4} \int_{-1}^{1} \int_{-1}^{1} \vec{\Psi} \, \vec{\Psi}^T \, d\tau \, d\xi = 
  \mat{\mathbb I} \in \reals^{\mpred\times\mpred}.
\end{equation}
Note that we are denoting the number of space-time basis functions by $\mpred$ to signify that this refers to the number of basis functions in the {\it prediction step}.
We catalog, at least up to fifth-order accuracy, how the one-dimensional indices $\ell_{\tau}$ and $\ell_{\xi}$ vary with
the index $\ell$ in \cref{table:spacetime_index}.
Using these basis functions, we write the predicted solution as follows:
\begin{equation}
\label{eqn:pred_ansatz}
\vec{\alpha}^{\text{ST}}\left(t^{n}+\frac{\Delta t}{2} (1+\tau), \, x_i +\frac{\Delta x}{2} \, \xi \right)   := \vec{\Psi}\left(\tau,\xi\right)^T \mat{W^{n+\half}_i},
\qquad \mat{W_i^{n+\half}} \in \reals^{\mpred \times \meq},
\end{equation}
for $(\tau,\xi) \in [-1,1]^2$, where $\mat{W}$ represents the matrix of unknown  coefficients.

\begin{table}[!t]
\begin{center}
\begin{Large}
\begin{tabular}{|cc||cc||cc||cc||cc|}
\hline
{\normalsize $\ell$} & {\normalsize $\left(\ell_\tau, \ell_\xi \right)$} &
{\normalsize $\ell$} & {\normalsize $\left(\ell_\tau, \ell_\xi \right)$} &
{\normalsize $\ell$} & {\normalsize $\left(\ell_\tau, \ell_\xi \right)$} &
{\normalsize $\ell$} & {\normalsize $\left(\ell_\tau, \ell_\xi \right)$} &
{\normalsize $\ell$} & {\normalsize $\left(\ell_\tau, \ell_\xi \right)$} \\
\hline\hline
{\normalsize 1} & {\normalsize (1,1)} & {\normalsize 4} & {\normalsize (1,3)} & {\normalsize 7} & {\normalsize (1,4)} & {\normalsize 10} & {\normalsize (4,1)} & {\normalsize 13} & {\normalsize (3,3)} \\
{\normalsize 2} & {\normalsize (1,2)} & {\normalsize 5} & {\normalsize (2,2)} & {\normalsize 8} & {\normalsize (2,3)} & {\normalsize 11} & {\normalsize (1,5)} & {\normalsize 14} & {\normalsize (4,2)} \\
{\normalsize 3} & {\normalsize (2,1)} & {\normalsize 6} & {\normalsize (3,1)} & {\normalsize 9} & {\normalsize (3,2)} & {\normalsize 12} & {\normalsize (2,4)} & {\normalsize 15} & {\normalsize (5,1)} \\
\hline
\end{tabular}
\caption{Index conversion table for the space-time basis, $\psi$, as defined in 
\cref{eqn:spacetime_basis}. Shown here are the basis element up to fifth-order of accuracy.
\label{table:spacetime_index}}
\end{Large}
\end{center}
\end{table}

Before describing how to compute the space-time coefficients, $
\mat{W}$, for the primitive variables, $\vec{\alpha}$,
we need to address one small issue: before the prediction step, the solution at time $t=t^n$ is given
only in terms of conservative variables (see equation \cref{eqn:corrector_ansatz}). 
In order to convert the conservative variable coefficients  from \cref{eqn:corrector_ansatz} 
to primitive variable coefficients,
\begin{equation}
\label{eqn:prim_ic}
   \vec{\alpha}^h\left(x_i +  \frac{\Delta x}{2} \, \xi \right) := \vec{\Phi}\left(\xi \right)^T \mat{A_i^n}
   \quad \text{for} \quad \xi \in [-1,1], \quad i=1,\ldots,\melems,
\end{equation}
we apply a simple $L^2[-1,1]$ projection:
\begin{equation}
\mat{A_i^n} = \frac{1}{2} \sum_{a=1}^{\morder} \omega_a \, \vec{\Phi}\left(\mu_a \right)  \left[ \vec{\alpha}\left( 
\vec{\Phi}\left(\mu_a \right)^T \mat{Q_i^n} \right) \right]^T \in \reals^{\mcorr \times \meq},
\end{equation}
where $\vec{\alpha}(\vec{q}): \reals^\meq \mapsto \reals^\meq$ gives the relationship between conservative
and primitive variables.

Next, an algebraic equation for the solution of the unknown coefficients in ansatz \cref{eqn:pred_ansatz} is 
obtained by multiplying \cref{eqn:alpha_theta} by $\vec{\Psi}$, integrating over 
$(\tau,\xi) \in [-1,1]^2$,  integrating-parts {\it only} in $\tau$ and {\it not} in $\xi$, and
making use of ansatz \cref{eqn:pred_ansatz}:
\begin{equation}
\label{eqn:nonlinear_pred}
\begin{split}
\mat{L} \, \, \vec{W^{n+\half}_{i \, (:,m)}} &= \frac{1}{4}\int_{-1}^{1} \int_{-1}^{1} 
  {\Theta}_m\left( \vec{\Psi}^T \mat{W^{n+\half}_i} \right) \vec{\Psi} \, d\tau \, d\xi \\ &+ 
  \left[ \frac{1}{4} \int_{-1}^{1} \, \vec{{\Psi}} \left(-1,\xi \right)  \, \vec{\Phi} \left(\xi \right)^T \, d\xi
 \right] \, \vec{A^n_{i \, (:,m)}},
 \end{split}
\end{equation}
for $m=1,\ldots,\meq$, where
\begin{gather}
\vec{W^{n+\half}_{i \, (:,m)}} \in \reals^\mpred : \, \text{$m^{\text{th}}$
column of $\mat{W^{n+\half}_i}$} \quad \text{(all coefficients for equation $m$)}, \\
\vec{A^{n}_{i \, (:,m)}} \in \reals^\mcorr : \, \text{$m^{\text{th}}$
column of $\mat{A^{n}_i}$} \quad \text{(all coefficients for equation $m$)}, \\
\mat{L} := \frac{1}{4} \int_{-1}^{1} \int_{-1}^{1} \vec{\Psi} \, \vec{\Psi}^T_{,\tau} \, d\tau \, d\xi
+ \frac{1}{4} \int_{-1}^{1} \vec{\Psi}_{|_{\tau=-1}} \vec{\Psi}_{|_{\tau=-1}}^T \, d\xi \in \reals^{\mpred \times \mpred}.
\end{gather}

We note that system \cref{eqn:nonlinear_pred} represents something akin to a single block-Jacobi iteration of a fully implicit spacetime DG approach \cite{article:KlaVegVen2006,article:Sudirham2006}, and is a set of  nonlinear algebraic equations that must be solved
independently on each space-time element. There are several techniques, including Newton's method for systems, that could be used solve these equations. However, following Gassner et al. \cite{article:GasDumHinMun2011}, we make use of an even simpler fixed-point iteration: the {\it Picard} iteration. 

After replacement of the
space-time integration with Gauss-Legendre quadrature, we can write the Picard iteration as
\begin{equation}
\label{eqn:picard}
\begin{split}
\vec{W^{n+\half}_{i \, (:,m)}} \leftarrow & \frac{1}{4} \sum_{a=1}^{\morder}
\sum_{b=1}^{\morder} \, \omega_a  \, \omega_b \, 
  \vec{\hat{\Psi}}\left(\mu_b, \, \mu_a \right)
   \, {\Theta}_m\left( \vec{\Psi}\left(\mu_b, \mu_a \right)^T \mat{W^{n+\half}_i} \right) \\
   + & \frac{1}{4} \sum_{b=1}^{\morder} \omega_{b} \, \vec{\hat{\Psi}} \left(-1,\xi_b \right)  \, \vec{\Phi} \left(\xi_b \right)^T \vec{A^{n}_{i \, (:,m)}},
   \end{split}
\end{equation}
for $m=1,\ldots,\meq$, where $\vec{\hat{\Psi}} = \mat{L}^{-1} \vec{\Psi}$ and
 $\omega_a$ and $\mu_a$ for $a=1,\ldots,\morder$ are the weights and abscissas of the $\morder$-point Gauss-Legendre quadrature rule.
This iteration -- like all fixed point iterations -- requires some appropriate initial guess; we explain how this
is done in \cref{subsec:full_algorithm}.

The two main advantages of the Picard iteration over Newton's method are: 
(1) it is Jacobian-free (the only inverse that must be 
computed is of $\mat{L}$, which is independent of the solution); and (2) the iteration converges to sufficient high-order accuracy after exactly $\morder$ iterations, obviating the need to compute residuals.
The basic principle of this approach is that each iteration improves the quality of the guess by one order of accuarcy: one iteration gives first-order accuracy, two iterations gives second-order 
accuracy, etc$\ldots$,
up to the maximum possible order of accuracy\footnote{Technically, the predicted solution is not even consistent with the underlying partial differential equation. What we mean by 
{\it high-order accuracy} in this context is what happens when the predicted solution is fed to the correction step. See \cref{subsec:corr_step} for more details.}: $\morder$ \cite{article:GasDumHinMun2011}.

\subsection{\it Correction step}
\label{subsec:corr_step}
The prediction step as outlined above is clearly not sufficient to produce a consistent numerical approximation
of hyperbolic conservation law \cref{eqn:conslaw}. Without having done a proper integration-by-parts in the spatial variable, the predicted solution on each space-time element is completely decoupled from all other space-time elements, which is inconsistent with the underlying partial differential equation.
Fortunately, there is a simple remedy that makes the solution not only {\it consistent} with conservation law \cref{eqn:conslaw}, but in fact  {\it high-order accurate} (under the assumption of sufficiently smooth solutions).
We refer to this remedy as the {\it correction step}, which is a single forward Euler-like step that makes use of the predicted solution.

To enact the correction step, we take hyperbolic conservation law \cref{eqn:conslaw}, multiply
by the spatial basis functions $\vec{\Phi}$ (see \cref{eqn:phi_basis_1d}), 
integrate over $(\tau,\xi) \in [-1,1]^2$, make use of ansatz \cref{eqn:corrector_ansatz},
apply integration-by-parts on the spatial variable, and replace all exact integration
by Gauss-Legendre quadrature:
\begin{gather}
\label{eqn:correction_update}
\begin{split}
\mat{Q^{n+1}_i} = \mat{Q^{n}_i} &+  \frac{\nu}{2} \sum_{a=1}^{\morder}
\sum_{b=1}^{\morder} \, \omega_a  \, \omega_b \, \vec{\Phi}_{,\xi}\left(\mu_a \right) \left[
 \vec{f} \left( \vec{\Psi}\left(\mu_b,\mu_a\right)^T \mat{W^{n+\half}_{i}} \right) \right]^T \\
&- {\nu}\left( \, \vec{\Phi}(1) \left[ \vec{{\mathcal F}^{n+\half}_{i+\half}} \right]^T
-  \vec{\Phi}(-1) \left[ \vec{{\mathcal F}^{n+\half}_{i-\half}} \right]^T \, \right),
\end{split}
\end{gather}
where $\omega_a$ and $\mu_a$ for $a=1,\ldots,\morder$ are the weights and abscissas of the $\morder$-point Gauss-Legendre quadrature rule.
The time-integrated numerical fluxes are defined using the predicted solution and the Rusanov \cite{article:Ru61} time-averaged flux: 
\begin{gather}
\label{eqn:correction_update_flux_1}
  \vec{{\mathcal F}^{n+\half}_{i-\half}} := \frac{1}{2} \sum_{a=1}^{\morder}
  \omega_a \, \vec{\mathcal F} \left( \mu_a \right),
\end{gather}
where the numerical flux at each temporal quadrature point is given by
\begin{gather}
  \label{eqn:correction_update_flux_2}
   \vec{\mathcal F} \left( \tau \right) := \frac{1}{2} \left( \vec{f}\left( \vec{W_{\text{R}}}(\tau) \right) + 
  \vec{f} \left( \vec{W_{\text{L}}}(\tau) \right) \right)  -  
       \frac{1}{2} \bigl|\lambda(\tau)\bigr| \Bigl( \vec{q} \left( \vec{W_{\text{R}}}(\tau) \right) 
         - \vec{q} \left( \vec{W_{\text{L}}}(\tau) \right) \Bigr),
\end{gather}
where
\begin{gather}
         \label{eqn:correction_update_flux_3}
  \vec{W_{\text{L}}}(\tau) := \vec{\Psi}\left(\tau,1\right)^T \mat{W^{n+\half}_{i-1}}, \qquad
  \vec{W_{\text{R}}}(\tau) := \vec{\Psi}\left(\tau,-1\right)^T \mat{W^{n+\half}_{i}},
\end{gather}
and $\bigl|\lambda(\tau)\bigr|$ is a local bound on the spectral radius of $\mat{A}\left(\vec{q} \right)$ in the neighborhood of interface $x=x_{i-\half}$ and at time $\tau$.


\section{Limiters for positivity-preservation and oscillation-control}
\label{sec:limiters}

In this section we give full details of the proposed limiting strategy.
In order to achieve discrete positivity-preservation and
non-oscillatory behavior in the presence of shocks and rarefactions, we need to apply 
limiters at various steps in the full algorithm. After making a few clarifying
definitions in \cref{subsec:positivity_defns}, we develop the proposed positivity limiter
in the prediction step in \cref{subsec:positivity_pred}, and the correction step in \cref{subsec:pos_corr_1,subsec:pos_corr_2}. We then develop the non-oscillatory limiter in
\cref{subsec:nonoscillatory}. Finally, we put all the pieces together and write out the full algorithm in \cref{subsec:full_algorithm}.

\subsection{\it Definition of the discrete positivity constraints}
\label{subsec:positivity_defns}
Before proceeding to the details of the prediction and correction step limiting strategies, it is useful to first define some notation and what we mean by the 
discrete positivity constraints.

\subsubsection{Choice of positivity points}
In Runge-Kutta discontinuous Galerkin schemes, the optimal points on which to enforce positivity (in the sense of achieving the minimal number of positivity points that allow the largest possible time-step) are the Gauss-Lobatto points \cite{article:ZhangShu11}.
For the positivity-preserving limiting strategy we propose in this work, the maximum allowable stable time-step is not directly tied to the choice of positivity points.
For this reason, we choose as our positivity points the Gauss-Legendre points augmented with the end points:
\begin{equation}
\label{eqn:space_pos_points}
{\mathbb X}_{\morder} := \Bigl\{ -1, 1 \Bigr\} \cup \Bigl\{ \text{roots of the $\morderth$ degree Legendre polynomial} \Bigr\},
\end{equation}
where $\morder$ is the desired order of accuracy. Note that ${\mathbb X}_{\morder}$
contains a total of $\morder+2$ points.
The reason for this choice is simple: for a fixed order of accuracy, $\morder$, all  purely spatial quadrature in the numerical scheme, both internally on the element and on the element faces, will only involve points taken from
${\mathbb X}_{\morder}$.

For the prediction step we require the two-dimensional version of \cref{eqn:space_pos_points}, which is the Cartesian product of ${\mathbb X}_{\morder}$
with itself:
\begin{equation}
\label{eqn:spacetime_pos_points}
{\mathbb X}^2_{\morder} := {\mathbb X}_{\morder} \otimes {\mathbb X}_{\morder}.
\end{equation}
Note that ${\mathbb X}^2_{\morder}$
contains a total of $(\morder+2)^2$ points.
Again, the reason for this choice is that all space-time quadrature in the numerical scheme,
both internally on the space-time element and on the space-time element faces,
will involve only points taken from ${\mathbb X}^2_{\morder}$.

\subsubsection{Discrete positivity constraints for the shallow water equations}
Strong hyperbolicity of the shallow water system is guaranteed 
if the height, $h(t,x)$, remains bounded away from zero for all $t\ge0$ and $x \in \left[
x_{\text{xlow}}, x_{\text{high}} \right]$. The discrete version of this positivity
constraint at time $t=t^n$ is defined separately for each element $\Tm_i$, 
and involves both the average and pointwise heights:
\begin{align}
\overline{h}^{\, n}_i := Q^n_{i \, \left(1,1\right)} \qquad \text{and} \qquad
h^n_i\left(\xi \right) := \vec{\Phi}\left(\xi \right)^T \vec{Q^n_{i \, \left(:,1 \right)}},
\end{align}
where we have used the following conventions:
\begin{alignat}{3}
\label{eqn:corr_access_1}
\mat{Q^n_i} \in \reals^{\mcorr \times \meq} &: \, \text{matrix of Legendre} &&\text{ coefficients in element $\Tm_i$ at time $t^n$}, \\
\label{eqn:corr_access_2}
Q^n_{i \, \left(\ell,k\right)} \in \reals &: \, \text{$(\ell,k)$ entry of $\mat{Q^n_i}$}  &&\text{($\ell = $ polynomial index, $k$ = equation index)}, \\
\label{eqn:corr_access_3}
\vec{Q^n_{i \, \left(:,k\right)}} \in \reals^\mcorr &: \, \text{$k^{\text{th}}$
column of $\mat{Q^n_i}$}  &&\text{(all polynomial coefficients for equation $k$)}, \\
\label{eqn:corr_access_4}
\vec{Q^n_{i \, \left(\ell,:\right)}} \in \reals^\meq &: \, \text{$\ell^{\text{th}}$
row of $\mat{Q^n_i}$}  &&\text{(polynomial coefficient $\ell$ for all equations)}.
\end{alignat}
In particular, we define two notions of discrete positivity: (1) positivity-in-the-mean, and (2) positivity at the points ${\mathbb X}_{\morder}$
defined in \cref{eqn:space_pos_points}. For the shallow water equations,
discrete positivity is defined via the following two sets:
\begin{alignat}{2}
\label{eqn:shllw_pos_set_1}
\overline{\mathcal P}^n_i  &:= \biggl\{ \mat{Q_i^n} \in \reals^{\mcorr \times 2}: \, \,
\overline{h}^{\, n}_i  \ge \varepsilon
\biggr\} \qquad &&(\text{positivity-in-the-mean}), \\
\label{eqn:shllw_pos_set_2}
{\mathcal P}^{ \, n}_i &:= \biggl\{ \mat{Q_i^n} \in \reals^{\mcorr \times 2}: \, \,
\min_{\xi \in {\mathbb X}_{\morder}} \Bigl\{
 h^n_i\left(\xi \right) \Bigr\}
\ge \varepsilon
\biggr\} \qquad &&(\text{positivity at \cref{eqn:space_pos_points}}),
\end{alignat}
for some $\varepsilon>0$.

\subsubsection{Discrete positivity constraints for the compressible Euler equations}
Well-posedness for compressible Euler equations is guaranteed 
if the density, $\rho(t,x)$, and pressure, $p(t,x)$, remain bounded away from zero for all $t\ge0$ and $x \in \left[
x_{\text{xlow}}, x_{\text{high}} \right]$. The discrete version of this positivity
constraint at time $t=t^n$ is defined separately for each element, $\Tm_i$, 
and involves the average and pointwise densities and pressures:
\begin{equation}
\label{eqn:corr_vars_on_element}
\begin{split}
\overline{\rho}^{\, n}_i := Q^n_{i \, \left(1,1\right)}, \quad
\overline{p}^{\, n}_i &:= \left(\gamma-1 \right) 
\left( Q^n_{i \, \left(1,3\right)} - \frac{\left[Q^n_{i \, \left(1,2\right)}\right]^2}{2\overline{\rho}^{\, n}_i} \right), \quad
\rho^n_i\left(\xi \right) := \vec{\Phi}\left(\xi \right)^T \vec{Q^n_{i \, \left(:,1 \right)}},
\\
p^n_i\left(\xi \right) &:= \left( \gamma - 1 \right) \left( \, \vec{\Phi}\left(\xi \right)^T \vec{Q^n_{i \, \left(:,3 \right)}} -  \frac{\left[\vec{\Phi}\left(\xi \right)^T \vec{Q^n_{i \, \left(:,2 \right)}} \right]^2}{2\rho^n_i\left(\xi \right)} \, \right),
\end{split}
\end{equation}
where again we have made use of the conventions from \cref{eqn:corr_access_1}--\cref{eqn:corr_access_4}.
In particular, we define two notions of discrete positivity: (1) positivity-in-the-mean, and (2) positivity at the points ${\mathbb X}_{\morder}$
defined in \cref{eqn:space_pos_points}. For the compressible Euler equations 
discrete positivity is defined via the following two sets:
\begin{align}
\label{eqn:euler_pos_set_1}
\overline{\mathcal P}^{ \, n}_i  &:=
 \biggl\{ \mat{Q_i^n} \in \reals^{\mcorr \times 3}: \, \,
\min\Bigl\{ \overline{\rho}^{\, n}_i, \, \overline{p}^{\, n}_i \Bigr\}
 \ge \varepsilon
\biggr\},\\
 \label{eqn:euler_pos_set_2}
{\mathcal P}^{ \, n}_i &:= \left\{ \mat{Q_i^n} \in \reals^{\mcorr \times 3}: \, \,
\min \left\{ \min_{\xi \in {\mathbb X}_{\morder}} \Bigl\{
 \rho^n_i\left(\xi \right) \Bigr\}, \, 
 \min_{\xi \in {\mathbb X}_{\morder}} \Bigl\{ p^n_i\left(\xi \right) \Bigr\} \right\}
\ge \varepsilon
\right\},
\end{align}
for some $\varepsilon>0$.

\subsection{\it Positivity-preservation in the prediction step}
\label{subsec:positivity_pred}
We stated in \cref{subsec:prediction_step} that an important flexibility in
the prediction step is the choice of the variables $\alpha$. For the sake of the simplest possible scheme for positivity-preservation, we now make specific choice
of using the true primitive variables: \cref{eqn:shllw_prim_matrix} for the shallow water equations and \cref{eqn:euler_prim_matrix} for the compressible Euler equations.

Let ${\mathbb I}_{\text{PrimPos}}$ be the set of 
equation indices of the predicted solution, \cref{eqn:pred_ansatz},
 for which positivity is required; for example,
 ${\mathbb I}_{\text{PrimPos}} = \{ 1 \}$ for shallow water and 
 ${\mathbb I}_{\text{PrimPos}} = \{ 1, 3 \}$ for compressible Euler.
Let ${\mathbb X}^2_{\morder}$ be the set of space-time positivity points 
defined by \cref{eqn:spacetime_pos_points}. Following the philosophy developed by
Zhang and Shu \cite{article:ZhangShu11} for the Runge-Kutta discontinuous Galerkin scheme, 
we seek the maximum value of $\theta \in [0,1]$
such that the space-time solution,
\begin{equation}
\label{eqn:pred_step_limiter_damping}
\begin{split}
\vec{\alpha}^{\text{ST}}\left(t^{n}+\frac{\Delta t}{2} (1+\tau), \, x_i +\frac{\Delta x}{2} \, \xi; \, \theta \right) &:= 
\vec{W^{n+\half}_{i \, (1, :)}} + \theta \, \sum_{\ell=2}^{\mpred}  \Psi_{ \, \ell} (\tau,\xi) \, \vec{W^{n+\half}_{i \, (\ell, :)}}, \\
&= (1-\theta) \, \vec{W^{n+\half}_{i \, (1, :)}} + \theta \, \vec{\Psi}\left(\tau,\xi\right)^T \mat{W^{n+\half}_i},
\end{split}
\end{equation}
is positive at all the space-time points ${\mathbb X}^2_{\morder}$ for all variables with index in ${\mathbb I}_{\text{PrimPos}}$. $\theta=0$ means that the solution is limited down to its cell
average (i.e., full limiting), while $\theta=1$ means that the full high-order approximation
can be used (i.e., no limiting).
Finding the optimal $\theta$ involves sampling the unlimited 
$\alpha^{\text{ST}}$ (i.e., \cref{eqn:pred_ansatz})  at all
the points in ${\mathbb X}^2_{\morder}$, computing the minimum over these point evaluations,
and then solving a linear scalar equation to find the parameter of $\theta$ that produces the
minimal damping to achieve positivity.

For example, if $k \in {\mathbb I}_{\text{PrimPos}}$, then we seek a value of $\theta_k$ such that:
\begin{equation}
   \left(1-\theta_k\right) {W^{n+\half}_{i \, (1, k)}} + \theta_k \alpha^k_{\text{min}} = \varepsilon > 0, 
   \quad \alpha^k_{\text{min}} := \min_{(\tau, \xi) \in {\mathbb X}^2_{\morder}}
   \left\{ \vec{\Psi}\left(\tau,\xi\right)^T \vec{W^{n+\half}_{i \, (:,k)}} \right\},
\end{equation}
which is a scalar linear equation that can be easily solved for $\theta_k$:
\begin{equation}
\theta_k = \min\left\{ 1, \frac{{W^{n+\half}_{i \, (1, k)}} - \varepsilon}{{W^{n+\half}_{i \, (1, k)}} - \alpha^k_{\text{min}}} \right\}.
\end{equation}
We need to do this for every index in ${\mathbb I}_{\text{PrimPos}}$, compute the minimum over all of these
$\theta_k$ values, and finally damp {\it all} variables with this minimum $\theta$ using
definition \cref{eqn:pred_step_limiter_damping}.
The full prediction step limiting process is summarized in \cref{alg:pred_limiter}. 

With this limiting procedure, we guarantee that all the conservative variable and flux function evaluations required in the
correction step, \cref{eqn:correction_update}--\cref{eqn:correction_update_flux_3}, 
involve only discrete primitive variables that satisfy the correct positivity constraint.
Mathematically, this means that we only ever evaluate conserved variables and fluxes
inside the convex set $S\in \reals^\meq$ over which the conservation law is hyperbolic.
Practically, this means we avoid computing square roots of negative numbers or dividing by zero.
It turns out, however, this simple limiting is insufficient
 to guarantee that the solution at the next step, $Q^{n+1}$,
satisfies the positivity constraint; in order to achieve positivity of $Q^{n+1}$, we also need
to apply positivity-preserving limiters in the correction step. 

\begin{algorithm}[!t]
\caption{Prediction step limiter. \label{alg:pred_limiter}}
\begin{algorithmic}[1]
 \Require \,
      $\melems$; \quad $\meq$; \quad $\mpred$;
       \quad $\varepsilon$;
       \quad ${\mathbb I}_{\text{PrimPos}}$;
       \quad ${\mathbb X}^2_{\morder}$;
       \quad $\mat{W^{n+\half}_i} \in \reals^{\mpred \times \meq}$ 
       \, $\forall i=1,\ldots,\melems$;
      
 \medskip
      
 \Ensure \, limited \, $\mat{W^{n+\half}_i} \in \reals^{\mpred \times \meq}$ \, 
 $\forall i=1,\ldots,\melems$;
 
 \medskip

\For{$i=1,\ldots,\melems$} \texttt{\quad \# loop over all elements}
\State
\State \texttt{\# find minimum value at augmented quadrature points and compute theta}
\State $\theta = 1$;
\For{$k\in {\mathbb I}_{\text{PrimPos}}$} \texttt{\quad \# loop over positivity variables}
\State $\alpha_{\text{min}} = \underset{(\tau, \xi) \in {\mathbb X}^2_{\morder}}{\min} 
\left\{ \vec{\Psi}\left(\tau,\xi\right)^T \vec{W^{n+\half}_{i \, (:,k)}} \right\};$ \quad
 $\theta = \min \left\{ \theta, \,  
\left({W^{n+\half}_{i \, (1,k)} - \varepsilon}\right)\bigg/\left({W^{n+\half}_{i \, (1,k)} - \alpha_{\text{min}}}
\right)
 \right\};$
\EndFor
\State
\State \texttt{\# if needed, limit all high-order space-time coefficients for all equations}
\If{$\theta<1$}
\For{$\ell=2,\ldots,\mpred$} \texttt{\quad \# loop over high-order space-time coefficients}
\For{$k=1,\ldots,\meq$} \texttt{\quad \# loop over all equations}
\State $W^{n+\half}_{i \, (\ell, k)}
  = \theta \, W^{n+\half}_{i \, (\ell, k)}$;
\EndFor
\EndFor
\EndIf
\State
\EndFor
\end{algorithmic}
\end{algorithm}

\subsection{\it Positivity-preservation in the correction step I: positivity-in-the-mean}
\label{subsec:pos_corr_1}
As described in \cref{subsec:positivity_defns}, there are two notions of discrete positivity:
(1) positivity-in-the-mean and (2) positivity at the augmented quadrature points
\cref{eqn:space_pos_points}. For the prediction step limiter described in \cref{subsec:positivity_pred}, we were able to ignore the {\it positivity-in-the-mean}
portion due the simplicity of the update (i.e., use of primitive variables and
an update that is completely local to the current element); and instead, it sufficed to
enforce positivity at the augmented space-time quadrature points \cref{eqn:spacetime_pos_points}.
For the correction step, we can no longer ignore the positivity-in-the-mean condition.

In order to achieve positivity-in-the-mean (e.g., \cref{eqn:shllw_pos_set_1} for shallow water and
\cref{eqn:euler_pos_set_1} for Euler), we employ a strategy that will
compare the cell average solution as computed by
our high-order scheme against a low-order scheme that is guaranteed to satisfy
the positivity-in-the-mean condition. In particular, if our scheme violates positivity-in-the-mean,
we will minimally limit the high-order fluxes so that the resulting cell average satisfies positivity. 

The idea of comparing high-order and low-order fluxes for the sake of limiting the high-order fluxes has a long history.  Harten and Zwas \cite{harten1972self} used such an idea
in their self-adjusting hybrid scheme. The flux-corrected transport (FCT) method developed
 by Boris, Book, and collaborators  \cite{boris1973flux,book1975flux,boris1976flux,book1981finite}
is also based on this idea. 
In the context of positivity-preservation for weighted essentially non-oscillatory (WENO)
schemes, this idea has been used by several authors in recent papers
\cite{xu2013,liang2014parametrized,tang14,christlieb2015high,article:seal2014explicit,ChFeSeTa2015,article:Xiong2013}.
In the context of discontinuous Galerkin schemes,
Xiong, Qiu, and Xu \cite{article:Xiong15} developed such an approach for scalar convection-diffusion equations. In this work we closely follow the flux limiting strategy developed by Moe, Rossmanith, and Seal \cite{article:MoRoSe17} for Lax-Wendroff discontinuous Galerkin schemes.

The basic idea of the proposed limiter is as follows. First, we compute
the high-order time-averaged numerical fluxes according to equations
\cref{eqn:correction_update_flux_1}--\cref{eqn:correction_update_flux_3}.
Next, we compute the Rusanov \cite{article:Ru61} (often called local Lax-Friedrichs) update
from $t=t^n$ to $t=t^n + \Delta t$, using as initial data
the cell averages, $Q^n_{i \, (1,:)}$:
\begin{align}
  \vec{Q^{\text{LxF}}_i} := \vec{Q^n_{i \, (1,:)}} - \nu
  \left( \, \vec{{\mathcal F}^{\text{LxF}}_{i+\half}} - 
  \vec{{\mathcal F}^{\text{LxF}}_{i-\half}} \, \right),
\end{align}
where $\nu = \Delta t/\Delta x$, the numerical flux is
\begin{align}
\label{eqn:LxF_fluxes}
\vec{{\mathcal F}^{\text{LxF}}_{i-\half}} := \frac{1}{2} \left[ \vec{f}\left(
  \, \vec{Q^n_{i \, (1,:)}} \, \right) + \vec{f}\left(
  \, \vec{Q^n_{i-1 \, (1,:)}} \, \right) \right] - \frac{1}{2} \bigl| \lambda \bigr|
  \left( \, \vec{Q^n_{i \, (1,:)}} - \vec{Q^n_{i-1 \, (1,:)}} \, \right),
\end{align}
and $\bigl|\lambda \bigr|$ is a local bound on the spectral radius of the flux Jacobian, 
$\mat{A}\left(\vec{q} \right)$, in the neighborhood of interface $x=x_{i-1/2}$
at time $t=t^n$. Given initial coefficients, ${Q^n_{i \, (1,:)}}$, that satisfy
the positivity constraints, we are guaranteed that ${Q^{\text{LxF}}_i}$ also satisfy
the positivity constraints under a suitable time-step restriction 
(see Perthame and Shu \cite{article:PerShu96} for an elegant proof).
Next, we update the cell averages via a limited flux:
\begin{equation}
\label{eqn:pos_in_mean_updated}
  \vec{Q^{n+1}_{i \, (1,:)}} = \vec{Q^\text{LxF}_{i \, (1,:)}} - \nu 
   \left( \, \theta_{i+\half} \, \vec{\Delta {\mathcal F}_{i+\half}}  
  - \theta_{i-\half} \, \vec{\Delta {\mathcal F}_{i-\half}} \, \right),
\end{equation}
where
\begin{equation}
\label{eqn:flux_differences}
\vec{\Delta {\mathcal F}_{i-\half}} := \vec{{\mathcal F}^{n+\half}_{i-\half}} 
   - \vec{{\mathcal F}^{\text{LxF}}_{i-\half}},
\end{equation}
and the maximum $\theta\in[0,1]$ on each face is 
chosen so that the updated solution 
satisfies the positivity constraints.
$\theta=0$ means that the solution is limited down to the positive local Lax-Friedrichs
cell average (i.e., full limiting), while $\theta=1$ means that the full high-order flux
can be used (i.e., no limiting).
Note that the high-order coefficients, $Q^{n+1}_{i \, (k,:)}$ for $k\ge 2$, are still updated using the full high-order flux as shown described by 
\cref{eqn:correction_update}--\cref{eqn:correction_update_flux_3}.

The final ingredient for obtaining positivity-in-the-mean is to determine a formula for
computing the values of $\theta_{i-1/2}\in[0,1]$ and $\theta_{i+1/2}\in[0,1]$ in \cref{eqn:pos_in_mean_updated}. We closely follow the methodology developed by
Moe, Rossmanith, and Seal \cite{article:MoRoSe17}. We summarize the process for computing
the optimal $\theta$ values in \cref{alg:corr_limiter_1}.

\begin{algorithm}[!t]
\caption{Correction step limiter I: positivity-in-the-mean. This algorithm is used
to determine the amount of damping on the high-order
fluxes in the correction step cell-average update \cref{eqn:pos_in_mean_updated}. \label{alg:corr_limiter_1}}
\begin{algorithmic}[1]
 \Require \,
      $\melems$; \quad $\nu$;
       \quad $\varepsilon$;
       \quad $\vec{\overline{Q}_i}:=\vec{Q^{n}_{i \, (1,:)}}
        \in \reals^{\meq}$ 
       \, $\forall i=1,\ldots,\melems$; \quad
       \quad $\vec{Q^{\text{LxF}}_{i}}
        \in \reals^{\meq}$ 
       \, $\forall i=1,\ldots,\melems$; 
       
       ${\mathcal F}^{n+\half}_{i-\half}
       \in \reals^{\meq}$
       \, $\forall i=1,\ldots,\melems+1$; \quad
       $\Delta{\mathcal F}_{i-\half}
       \in \reals^{\meq}$
       \, $\forall i=1,\ldots,\melems+1$;
      
 \medskip
      
 \Ensure \, $\theta_{i-\half}$ \, 
 $\forall i=1,\ldots,\melems+1$;
 
 \medskip

\State $\vec{\theta} = \vec{1}$; \texttt{\quad \# initialize all theta values to one}
\State
\For{$i=1,\ldots,\melems$} \texttt{\quad \# loop over all elements}
\State 
\State \texttt{\# Part I: limit for positive height (shallow water) or density (Euler)}
\State $\Gamma = \left( Q^{\text{LxF}}_{i (1)} - \varepsilon \right)\big/ \nu$;
\quad $\Delta{\mathcal F}_{\text{left}} = \Delta{\mathcal F}_{i-\half \, (1)}$;
\quad $\Delta{\mathcal F}_{\text{right}} = \Delta{\mathcal F}_{i+\half \, (1)}$;
\quad $\Lambda_{\text{left}} = 1$; \quad $\Lambda_{\text{right}} = 1$;
\State
\If{$\left( \Delta{\mathcal F}_{\text{left}}>0 \right) \, \text{\bf and} \,
    \left( \Delta{\mathcal F}_{\text{right}}<0 \right)$}
\State  $\Lambda_{\text{left}} = \Lambda_{\text{right}} = \min\Bigl\{ 
               1, \, \Gamma\big/\left( \big| \Delta {\mathcal F}_{\text{left}} \big| + 
               \big| \Delta {\mathcal F}_{\text{right}} \big| \right) \Bigr\}$;
               \texttt{\quad \# outflow on both faces}
\ElsIf {$\left( \Delta{\mathcal F}_{\text{left}}>0 \right)$}
\State $\Lambda_{\text{left}} = \min\Bigl\{ 
               1, \, \Gamma\big/\big| \Delta {\mathcal F}_{\text{left}} \big|  \Bigr\}$;
               \texttt{\quad \# outflow only on left face}
\ElsIf {$\left( \Delta{\mathcal F}_{\text{right}}<0 \right)$}
\State $\Lambda_{\text{right}} = \min\Bigl\{ 
               1, \, \Gamma\big/ 
               \big| \Delta {\mathcal F}_{\text{right}} \big| \Bigr\}$;
                \texttt{\quad \# outflow only on right face}
\EndIf
\State 
\State \texttt{\# Part II: limit for positive pressure (Euler case only)}
\State $p^{\text{LxF}} = (\gamma-1) \left(Q^{\text{LxF}}_{i \, (3)}
            - \left( Q^{\text{LxF}}_{i \, (2)} \right)^2 \Big/ \left( 2 \, Q^{\text{LxF}}_{i \, (1)} \right) \right)$; \quad $\mu_{11} = 1$; \quad $\mu_{10} = 1$; \quad $\mu_{01} = 1$; 
\State
\State \texttt{\# case 1:  limit based on both left and right face fluxes}
\State $\vec{Q^{\star}} = \vec{\overline{Q}_i} - \nu \left( \Lambda_{\text{right}} 
          \, \vec{{\mathcal F}_{i+\half}^{n+\half}} 
          - \Lambda_{\text{left}} \, \vec{{\mathcal F}_{i-\half}^{n-\half}} \right)$;
\quad   $p^{\star} = (\gamma-1) \left( Q^{\star}_{(3)} - \left(Q^{\star}_{(2)} \right)^2 \Big/\left(2 \, 
         Q^{\star}_{(1)} \right) \right)$;
\State {{\bf if} \, $\left(p^{\star}<\varepsilon\right)$:} \quad
$\mu_{11} = \left(p^{\text{LxF}}-\varepsilon \right)\Big/\left(p^{\text{LxF}}-p^{\star}\right)$; 
 \quad \text{\bf end if}          
\State
\State \texttt{\# case 2:  limit based on zero right face flux}
\State $\vec{Q^{\star}} = \vec{\overline{Q}_i} + \nu \,
      \Lambda_{\text{left}} \, \vec{{\mathcal F}_{i-\half}^{n-\half}}$;
\quad   $p^{\star} = (\gamma-1) \left( Q^{\star}_{(3)} - \left(Q^{\star}_{(2)} \right)^2 \Big/\left(2 \, 
         Q^{\star}_{(1)} \right) \right)$;
\State {{\bf if} \, $\left(p^{\star}<\varepsilon\right)$:} \quad
$\mu_{10} = \left(p^{\text{LxF}}-\varepsilon \right)\Big/\left(p^{\text{LxF}}-p^{\star}\right)$; 
 \quad \text{\bf end if}
\State
\State \texttt{\# case 3:  limit based on zero left face flux}
\State $\vec{Q^{\star}} = \vec{\overline{Q}_i} - \nu \, \Lambda_{\text{right}} 
          \, \vec{{\mathcal F}_{i+\half}^{n+\half}}$;
\quad   $p^{\star} = (\gamma-1) \left( Q^{\star}_{(3)} - \left(Q^{\star}_{(2)} \right)^2 \Big/\left(2 \, 
         Q^{\star}_{(1)} \right) \right)$;
\State {{\bf if} \, $\left(p^{\star}<\varepsilon\right)$:} \quad
$\mu_{01} = \left(p^{\text{LxF}}-\varepsilon \right)\Big/\left(p^{\text{LxF}}-p^{\star}\right)$; 
 \quad \text{\bf end if}
\State
\State $\mu = \text{min}\left\{ \mu_{11}, \, \mu_{10}, \, \mu_{01} \right\}$;
\quad $\Lambda_{\text{left}} = \mu \, \Lambda_{\text{left}}$; \quad
$\Lambda_{\text{right}} = \mu \, \Lambda_{\text{right}}$;
\State $\theta_{i-\half} = \min\left\{ \theta_{i-\half}, \, \Lambda_{\text{left}} \right\}$;
\quad $\theta_{i+\half} = \min\left\{ \theta_{i+\half}, \, \Lambda_{\text{right}} \right\}$;
\State
\EndFor
\end{algorithmic}
\end{algorithm}

\subsection{\it Positivity-preservation in the correction step II: positivity at quadrature points}
\label{subsec:pos_corr_2}
Once we have ensured that the new solution, $Q^{n+1}$, satisfies the positive-in-the-mean condition, 
we now seek to enforce positivity at the augmented quadrature points \cref{eqn:space_pos_points}
 (e.g., \cref{eqn:shllw_pos_set_2} for shallow water and
\cref{eqn:euler_pos_set_2} for Euler).
This limiting step is similar to what was done for prediction step limiting (see \cref{subsec:positivity_pred}), but with the added complication that we are now working
with conserved variables. 

Following the philosophy developed by
Zhang and Shu \cite{article:ZhangShu11} for the Runge-Kutta discontinuous Galerkin scheme,
the idea is to find the maximum $\theta\in[0,1]$ such that
\begin{equation}
\begin{split}
\vec{q}^{h}\left( t^{n+1}, x_i +\frac{\Delta x}{2} \, \xi; \theta \right) :=& \, \,
\vec{Q^{n+1}_{i \, (:, 1)}} + \theta \, \sum_{k=2}^{\mcorr}  \Phi_{ \, k} (\xi) \, \vec{Q^{n+1}_{i \, (:, k)}} \\
=& \left(1-\theta \right) \, \vec{Q^{n+1}_{i \, (:, 1)}} + \theta \,  \vec{\Phi} (\xi)^T \, 
\mat{Q^{n+1}_{i}}
\end{split}
\end{equation}
satisfies the appropriate positivity constraints at all points in ${\mathbb X}_{\morder}$.
$\theta=0$ means that the solution is limited down to its cell
average (i.e., full limiting), while $\theta=1$ means that the full high-order approximation
can be used (i.e., no limiting).
For all variables that are both conservative and required to positive (e.g., the height in the shallow water equations and the density in the Euler equations), finding the optimal $\theta$ involves solving a scalar linear equation (just as in the prediction step limiter from
\cref{subsec:positivity_pred}). However,
some variables that are required to be positive may not be conservative variables (e.g., the pressure in the Euler equations), and thus finding the optimal $\theta$ requires solving nonlinear
equations. However, if we give up on finding the exact optimizer, we can linearize this process by
invoking convexity of the pressure:
\begin{equation}
\begin{split}
 p^{\star}(\theta) :=& \, \, \left( \gamma - 1 \right) \min_{\xi \in {\mathbb X}_{\morder}} \left\{  q^h_3\left(t^{n+1},x_i + \frac{\Delta x}{2} \xi; \, \theta \right) -  \frac{1}{2} \frac{q^h_2\left(t^{n+1},x_i +  \frac{\Delta x}{2} \xi; \, \theta \right)^2}{q^h_1\left(t^{n+1},x_i + \frac{\Delta x}{2} \xi; \, \theta \right)}  \right\} \\
 \ge& \, \, p^{\star}(1) + (1-\theta) \left( p^{\star}(0) - p^{\star}(1) \right) = 
 p_{\text{min}} + (1-\theta) \left( \overline{p}^{ \, n+1}_i - p_{\text{min}} \right),
 \end{split}
\end{equation}
for $\theta \in [0,1]$, where $p_{\text{min}}$ is the minimum pressure over all the points in ${\mathbb X}_{\morder}$
of the unlimited solution, and $\overline{p}$ is the cell average pressure defined in
\cref{eqn:corr_vars_on_element}. Finding the near-optimal $\theta \in [0,1]$ according to the
above linearization is straightfoward:
\begin{equation}
\theta = \min\left\{ 1, \frac{p_{\text{min}}-\varepsilon}{p_{\text{min}}-\overline{p}^{ \, n+1}_i}\right\},
\end{equation}
for some $\varepsilon>0$.
We summarize the full limiting procedure for both the shallow water and compressible Euler
equations in \cref{alg:corr_limiter_2}.

\begin{algorithm}[!t]
\caption{Correction step limiter II: positivity at augmented 
quadrature points ${\mathbb X}_{\morder}$. \label{alg:corr_limiter_2}}
\begin{algorithmic}[1]
 \Require \,
      $\melems$; \quad $\meq$; \quad $\mcorr$;
       \quad $\varepsilon$;
       \quad ${\mathbb X}_{\morder}$;
       \quad $\mat{Q^{n+1}_i} \in \reals^{\mcorr \times \meq}$ 
       \, $\forall i=1,\ldots,\melems$;
      
 \medskip
      
 \Ensure \, limited \, $\mat{Q^{n+1}_i} \in \reals^{\mcorr \times \meq}$ \, 
 $\forall i=1,\ldots,\melems$;
 
 \medskip

\For{$i=1,\ldots,\melems$} \texttt{\quad \# loop over all elements}
\State
\State \texttt{\# Part I: limit for positive height (shallow water) or density (Euler)}
\State $\overline{\rho}^{ \, n+1}_{i} = Q^{n+1}_{i \, (1,1)}$; \quad $\rho_{\text{min}} 
= \underset{\xi \in {\mathbb X}_{\morder}}{\min}  
\left\{ \vec{\Phi}\left(\xi\right)^T \vec{Q^{n+1}_{i \, (:,1)}} \right\};$ \quad
 $\theta = \min \left\{ 1, \,  
\left(\overline{\rho}^{ \, n+1}_{i} - \varepsilon\right)\bigg/\left(\overline{\rho}^{ \, n+1}_{i} - \rho_{\text{min}}
\right)
 \right\};$
\If{$\theta<1$}
\For{$\ell=2,\ldots,\mcorr$} \texttt{\quad \# loop over high-order coefficients}
\For{$k=1,\ldots,\meq$} \texttt{\quad \# loop over all equations}
\State $Q^{n+1}_{i \, (\ell, k)}
  = \theta \, Q^{n+1}_{i \, (\ell, k)}$;
\EndFor
\EndFor
\EndIf
\State
\State \texttt{\# Part II: limit for positive pressure (Euler case only)}
\State $\overline{p}^{\, n+1}_i = \left(\gamma-1 \right) 
\left( Q^{n+1}_{i \, \left(1,3\right)} - {\left(Q^{n+1}_{i \, \left(1,2\right)}\right)^2}
\Big/\left(2 Q^{n+1}_{i \, \left(1,1\right)} \right) \right)$; 
\State $p^{n+1}_i\left(\xi \right) = \left( \gamma - 1 \right) \left( \, \vec{\Phi}\left(\xi \right)^T \vec{Q^{n+1}_{i \, \left(:,3 \right)}} -  {\left(\vec{\Phi}\left(\xi \right)^T \vec{Q^{n+1}_{i \, \left(:,2 \right)}} \right)^2}\Big/\left({2\vec{\Phi}\left(\xi \right)^T \vec{Q^{n+1}_{i \, \left(:,1 \right)}}}\right) \, \right)$;
\State $p_{\text{min}} = \underset{\xi \in {\mathbb X}_{\morder}}{\min} 
\left\{ p^{n+1}_i\left(\xi\right) \right\};$ \quad
 $\theta = \min \left\{ 1, \,  
\left(\overline{p}^{ \, n+1}_i - \varepsilon\right)\bigg/\left(\overline{p}^{ \, n+1}_i - p_{\text{min}}
\right)
 \right\};$
\If{$\theta<1$}
\For{$\ell=2,\ldots,\mcorr$} \texttt{\quad \# loop over high-order coefficients}
\For{$k=1,\ldots,\meq$} \texttt{\quad \# loop over all equations}
\State $Q^{n+1}_{i \, (\ell, k)}
  = \theta \, Q^{n+1}_{i \, (\ell, k)}$;
\EndFor
\EndFor
\EndIf
\State
\EndFor
\end{algorithmic}
\end{algorithm}

\subsection{\it Controlling unphysical oscillations}
\label{subsec:nonoscillatory}
The limiters described in \cref{subsec:positivity_pred,subsec:pos_corr_1,subsec:pos_corr_2} 
 guarantee positivity, but they are generally not sufficient to damp out all unphysical oscillations at shocks and rarefactions. 
In order to eliminate these oscillations we augment the method with one final limiter.
Through numerical experiments we have found that applying a limiting strategy similar to  the one developed by Krivodonova \cite{article:Kriv07}, {\it once per time-step} after the correction step update \cref{eqn:correction_update}, provides the necessary limiting to remove unphysical oscillations without unduly diffusing the numerical solution.

This Krivodonova \cite{article:Kriv07} limiter is applied
 on the {\it characteristic variables}:
\begin{equation}
 \vec{q} = \mat{R} \, \vec{c} \quad \Longleftrightarrow \quad 
 \vec{c} = \mat{L} \, \vec{q},
\end{equation}
where $\mat{R}$ and $\mat{L}$ are the matrices of right and left-eigenvectors of the flux Jacobian
\cref{eqn:quasilinear_cons}, respectively.
The limiter is applied in a hierarchical manner starting from the highest degree Legendre coefficient, $Q_{(\mcorr,:)}$, down to the second lowest coefficient,
$Q_{(2,:)}$. The lowest coefficient, $Q_{(1,:)}$, which is the cell average, is never limited in order to maintain the conservative property of the scheme.

In each element, $i$, for each characteristic variable, $m$, and for each of the Legendre coefficients from the highest, $k=\mcorr$, down to the second lowest, $k=2$, we compare the current coefficient to two one-sided finite differences
of coefficients of one lower order:
\begin{equation}
\label{eqn:minmod_limiter}
\begin{split}
\vec{\widehat{L}_{(m,:)}} \cdot \vec{Q^{n+1}_{i \, (k,:)}} \leftarrow 
\text{minmod}\Biggl( \vec{\widehat{L}_{(m,:)}} \cdot \vec{Q^{n+1}_{i \, (k,:)}}, \, \,
        &a_k \, \vec{\widehat{L}_{(m,:)}} \cdot \left(
           \vec{Q^{n+1}_{i (k-1,:)}} - \vec{Q^{n+1}_{i-1 (k-1,:)}} \right), \\
        &a_k \, \vec{\widehat{L}_{(m,:)}} \cdot \left(
           \vec{Q^{n+1}_{i+1 (k-1,:)}} - \vec{Q^{n+1}_{i (k-1,:)}} \right) \Biggr),
\end{split}
\end{equation}
where $a_k = \sqrt{(2k-1)/(2k+1)}$ is the largest possible constant allowed in the 
Krivodonova \cite{article:Kriv07} limiter, which results in the least aggressive limiter
possible in this framework.
The minmod function with three arguments is defined as follows:
\begin{equation}
\text{minmod}(a,b,c) = \begin{cases}
0 & \quad \text{if} \quad \text{min}\Bigl\{ a  b, b  c, a  c \Bigr\} \le 0, \\
\text{sign}(a) \cdot \text{min}\Bigl\{ |a|,|b|,|c| \Bigr\}
& \quad \text{otherwise}.
\end{cases}
\end{equation}
We note that the matrices of right and left-eigenvectors the flux Jacobian \cref{eqn:quasilinear_cons}, $\mat{R}$ and $\mat{L}$, depend on the solution; 
in this work we evaluate both of these matrices at the cell averages:
\begin{equation}
\mat{\widehat{R}}= \mat{R}\left( \, \vec{Q^{n}_{i (1,:)}} \, \right)
\qquad \text{and} \qquad
\mat{\widehat{L}} = \left( \, \mat{\widehat{R}} \, \right)^{-1},
\end{equation}
and denote the $m^{\text{th}}$ row of $\mat{\widehat{L}}$ by $\vec{\widehat{L}_{(m,:)}}$.
Our version of the Krivodonova \cite{article:Kriv07} limiting procedure is detailed in
\cref{alg:kriv_limiter}.

\begin{algorithm}[!t]
\caption{Limiter to remove any unphysical oscillations not handled by
the positivity limiters. This is a modified version of the Krivodonova \cite{article:Kriv07} limiter. \label{alg:kriv_limiter}}
\begin{algorithmic}[1]
 \Require \,
      $\melems$; \quad $\meq$; \quad $\mcorr$; \quad $\varepsilon$;
       \quad $\mat{Q^{n+1}_i} \in \reals^{\mcorr \times \meq}$ 
       \, $\forall i=1,\ldots,\melems$; \quad 
       $\mat{R}(\vec{q}): \reals^\meq \mapsto \reals^{\meq\times\meq}$;
      
 \medskip
      
 \Ensure \, limited \, $\mat{Q^{n+1}_i} \in \reals^{\mcorr \times \meq}$ \, 
 $\forall i=1,\ldots,\melems$;
 
 \medskip
 
 \State \texttt{\# convert variables from conservative to characteristic}
 \For{$i=1,\ldots,\melems$} \texttt{\quad \# loop over all elements}
 \State $\mat{\widehat{L}} = \mat{R}^{-1}\left( \vec{Q^{n+1}_{i (1,:)}} \right)$;
   \For{$k=1,\ldots,\mcorr-1$} \texttt{\quad \# loop over Legendre moments}
     \State 
     $\vec{C_{i (k,:)}} =  \mat{\widehat{L}} \, \vec{Q^{n+1}_{i (k+1,:)}}$
     \quad $\vec{\Delta_{-} C_{i (k,:)}} =  \mat{\widehat{L}} \left(
           \vec{Q^{n+1}_{i (k,:)}} - \vec{Q^{n+1}_{i-1 (k,:)}} \right)$, \quad
      $\vec{\Delta_{+} C_{i (k,:)}} = \mat{\widehat{L}} \left(
           \vec{Q^{n+1}_{i+1 (k,:)}} - \vec{Q^{n+1}_{i (k,:)}} \right)$;
     \State \texttt{\# Boundary conditions are needed to define $i=0$ and $i=\melems+1$ cases}
   \EndFor
 \EndFor
 \State
 \State {\tt \# hierarchical limiting procedure}
 \For{$i=1,\ldots,\melems$} \texttt{\quad \# loop over all elements}
   \For{$\ell=1,\ldots,\meq$} \texttt{\quad \# loop over all characteristic variables}
        \State $\text{mstop} = 0$, \, $k = \morder-1$
        \While{$\text{mstop} = 0$}
        \State $C^{\star} = C_{i (k,\ell)}$, \quad
         $C_{i (k,\ell)} = \text{minmod}\left( C^{\star}, \,
        \sqrt{\frac{2k-1}{2k+1}} \, \Delta_{+} C_{i (k,\ell)}, \, 
        \sqrt{\frac{2k-1}{2k+1}} \, \Delta_{-} C_{i (k,\ell)} \right)$;
        \If{ \, $\left( k>1 \right)$ \, \text{\bf and} \,
        $\Bigl( \Bigl| C_{i (k,\ell)} - C^{\star} \Bigr| > \varepsilon
        \, \, \, \, \text{\bf or} \, \, \, \, \Bigl| C_{i (k,\ell)} \Bigr| \le \varepsilon \Bigr)$ \, }
        \State $k=k-1$;  \texttt{\quad \# move down to next Legendre moment}
        \Else 
        \State $\text{mstop} = 1$;  \texttt{\quad \# reached lowest
        moment or solution is smooth enough}
        \EndIf
        \EndWhile
   \EndFor
 \EndFor
 \State
 \State \texttt{\# convert variables back to conservative from characteristic}
 \For{$i=1,\ldots,\melems$} \texttt{\quad \# loop over all elements}
   \State $\mat{\widehat{R}}= \mat{R}\left( \vec{Q^{n+1}_{i (1,:)}} \right)$;
   \For{$k=2,\ldots,\mcorr$} \texttt{\quad \# loop over Legendre moments}
     \State $\vec{Q^{n+1}_{i (k,:)}} =  \mat{\widehat{R}} \, \vec{C_{i (k-1,:)}}$;
   \EndFor
 \EndFor
\end{algorithmic}
\end{algorithm}

\subsection{\it Full algorithm: one complete time-step with limiters}
\label{subsec:full_algorithm}
Finally, in order to  clearly demonstrate where each limiter is applied in the course of a single time-step, we have summarized the full scheme over one time-step in 
\cref{alg:full_method}.

\begin{algorithm}[!th]
\caption{One full time-step of the locally-implicit Lax-Wendroff DG scheme.\label{alg:full_method}}
\begin{algorithmic}[1]
 \Require \,
      $\nu = \frac{\Delta t}{\Delta x}$; \quad
      $\melems$; \quad $\meq$; \quad $\morder$; \quad $\mcorr$; \quad $\mpred$; \quad $\varepsilon$;
       \quad $\mat{Q^{n}_i} \in \reals^{\mcorr \times \meq}$ 
       \, $\forall i=1,\ldots,\melems$;
       
       $\vec{\omega}, \, \vec{\mu} \in \reals^{\morder}$; \, \texttt{\#Gaussian quad\#}; \quad
       $\vec{\alpha}(\vec{q}): \reals^\meq \mapsto \reals^\meq$;  \, \texttt{\#Cons-to-prim\#};
       
       $\vec{\Phi}(\xi): [-1,1] \mapsto \reals^{\mcorr}$; \quad
       $\vec{\Psi}(\tau,\xi): [-1,1]^2 \mapsto \reals^{\mpred}$; \, \texttt{\#Legendre polys\#};
       
       $\vec{f}\left(\vec{q}\right): \reals^\meq \mapsto \reals^\meq$; \, \texttt{\#flux\#};
        \quad
       $\mat{R}\left(\vec{q}\right): \reals^\meq \mapsto \reals^{\meq \times \meq}$;
        \, \texttt{\#right e-vecs flux J\#};
       
       $\big|\lambda\left(\vec{q}\right)\bigr|: \reals^\meq \mapsto \reals$; \, \texttt{\#spectral
       rad of flux J\#};
       \quad
       $\vec{\Theta}(\vec{q}):\reals^\meq \mapsto \reals^\meq$;
        \, \texttt{\#eqn \cref{eqn:alpha_theta_comp}\#};
      
 \medskip
      
 \Ensure \, $\mat{Q^{n+1}_i} \in \reals^{\mcorr \times \meq}$ \, 
 $\forall i=1,\ldots,\melems$;
 
 \medskip
 
 \State \texttt{\# initial guess for Picard iteration}
 \For{$i=1,\ldots,\melems$} \texttt{\quad \# loop over all elements}
 \State $\mat{A_i^n} = \frac{1}{2} \sum\limits_{a=1}^{\morder} \omega_a \, \vec{\Phi}\left(\mu_a \right)  \left[ \vec{\alpha}\left( 
\vec{\Phi}\left(\mu_a \right)^T \mat{Q_i^n} \right) \right]^T \in \reals^{\mcorr \times \meq}$; \,
\texttt{\# convert cons to prim vars}
\State $\mat{W^{n+\half}_{i}} = \left[ 
 \frac{1}{4} \sum\limits_{a=1}^{\morder} \sum\limits_{b=1}^{\morder}
 \omega_a \, \omega_b \, \vec{\Psi} \left(\tau_a, \xi_b\right)  \vec{\Phi} \left(\xi_b\right)^T
 \right] \, \mat{A^{n}_{i}}$; \,
 \texttt{\# extend prim vars to const in time}
 \EndFor
 \State
 \State \texttt{\# Picard iteration}
 \State $\vec{\hat{\Psi}}(\tau,\xi) = \left[ \sum\limits_{a=1}^{\morder}
 \sum\limits_{b=1}^{\morder} \frac{\omega_a\omega_b}{4} \vec{\Psi}\left(\mu_a, \mu_b \right)
  \, \vec{\Psi}_{,\tau}\left(\mu_a, \mu_b \right)^T
+ \sum\limits_{a=1}^{\morder} \frac{\omega_a}{4} \vec{\Psi}\left(-1,\mu_a\right)
 \vec{\Psi}\left(-1,\mu_a \right)^T \right]^{-1} \, \vec{\Psi}(\tau,\xi)$;
 \For{$I=1,\ldots,\morder$}
 \For{$i=1,\ldots,\melems$} \texttt{\quad \# loop over all elements}
 \State $\mat{W^{n+\half}_{i}} \leftarrow   \sum\limits_{a=1}^{\morder}
\sum\limits_{b=1}^{\morder}  \frac{\omega_a  \, \omega_b}{4} 
  \vec{\hat{\Psi}}\left(\mu_b, \, \mu_a \right) 
   \, \left[\vec{\Theta}\left( \vec{\Psi}\left(\mu_b, \mu_a \right)^T \mat{W^{n+\half}_i} \right)
   \right]^T
   +  \sum\limits_{b=1}^{\morder} \frac{\omega_{b}}{4} \vec{\hat{\Psi}} \left(-1,\xi_b \right)   \vec{\Phi} \left(\xi_b \right)^T \mat{A^{n}_{i}}$;
   \EndFor
   \State $\mat{W^{n+\half}_{i}} \,  \leftarrow$ $\left\{ \text{apply prediction step limiter \cref{alg:pred_limiter} to $\mat{W^{n+\half}_{i}}$} \right\}$;
   \EndFor
   \State
   \State \texttt{\# compute numerical fluxes and Rusanov update}
   \State  $\vec{Q^{\text{LxF}}_i} = \vec{Q^n_i} \quad \forall i=1,\ldots,\melems$;
   \quad \texttt{\# initialize Rusanov update}
   \State {use BCs to set:} \quad $\vec{{\mathcal F}^{n+\half}_{i-\half}}$,
   \, $\vec{{\mathcal F}^{\text{LxF}}_{i-\half}}$, 
   \, $\vec{\Delta {\mathcal F}_{i-\half}} = \vec{{\mathcal F}^{n+\half}_{i-\half}}-
   \vec{{\mathcal F}^{\text{LxF}}_{i-\half}}$ \, for $i=1$ and $i=\melems+1$;
   \State $\vec{Q^{\text{LxF}}_1} \leftarrow \vec{Q^{\text{LxF}}_1} +
    \nu \vec{{\mathcal F}^{\text{LxF}}_{\half}}$; \quad
     $\vec{Q^{\text{LxF}}_\melems} \leftarrow \vec{Q^{\text{LxF}}_\melems} -  \nu \vec{{\mathcal F}^{\text{LxF}}_{\melems+\half}}$;
    \quad \texttt{\# partial LxF udpate}
   \For {$i=2,\ldots,\melems$} \texttt{\quad \# loop over all interior faces}
   \State compute $\vec{{\mathcal F}^{n+\half}_{i-\half}}$ via 
   \cref{eqn:correction_update_flux_1}--\cref{eqn:correction_update_flux_3}; \quad
   compute $\vec{{\mathcal F}^{\text{LxF}}_{i-\half}}$ via 
   \cref{eqn:LxF_fluxes}; \quad
   $\vec{\Delta {\mathcal F}_{i-\half}} = \vec{{\mathcal F}^{n+\half}_{i-\half}}-
   \vec{{\mathcal F}^{\text{LxF}}_{i-\half}}$;
   \State $\vec{Q^{\text{LxF}}_{i-1}} \leftarrow \vec{Q^{\text{LxF}}_{i-1}} -  \nu \vec{{\mathcal F}^{\text{LxF}}_{i-\half}}$ \quad \text{and} \quad
   $\vec{Q^{\text{LxF}}_{i}} \leftarrow \vec{Q^{\text{LxF}}_{i}} +  \nu \vec{{\mathcal F}^{\text{LxF}}_{i-\half}}$;
   \EndFor 
   \State $\vec{\theta} \,  \leftarrow \left\{ \text{apply correction step limiter I \cref{alg:corr_limiter_1}} \right\}$;
   \State
   \State \texttt{\# correction step update}
    \For{$i=1,\ldots,\melems$} \texttt{\quad \# loop over all elements}
   \State $\vec{Q^{n+1}_{i \, (1,:)}} = \vec{Q^\text{LxF}_{i \, (1,:)}} - \nu 
   \left( \, \theta_{i+\half} \, \vec{\Delta {\mathcal F}_{i+\half}}  
  - \theta_{i-\half} \, \vec{\Delta {\mathcal F}_{i-\half}} \, \right)$;
  \For {$k=2,\ldots,\mcorr$} \texttt{\quad \# loop over high-order moments}
  \State $N_k = \frac{1}{2} \sum\limits_{a=1}^{\morder}
\sum\limits_{b=1}^{\morder} \, \omega_a  \, \omega_b \, {\Phi}_{k,\xi}\left(\mu_a \right)
 \vec{f} \left( \vec{\Psi}\left(\mu_b,\mu_a\right)^T \mat{W^{n+\half}_{i}} \right)$;
  \State $\vec{Q^{n+1}_{i \, (k,:)}} = \vec{Q^{n}_{i \, (k,:)}} 
  +  \nu N_k
- {\nu}\left( \, {\Phi}_k(1) \vec{{\mathcal F}^{n+\half}_{i+\half}}
-  {\Phi}_k(-1)  \vec{{\mathcal F}^{n+\half}_{i-\half}} \, \right)$
  \EndFor
   \EndFor
   \State $\mat{Q^{n+1}_i} \leftarrow \left\{ \text{apply non-oscillatory limiter \cref{alg:kriv_limiter} to } \mat{Q^{n+1}_i} \right\}$;
   \State $\mat{Q^{n+1}_i} \leftarrow \left\{ \text{apply correction step limiter II \cref{alg:corr_limiter_2} to } \mat{Q^{n+1}_i} \right\}$;
   \end{algorithmic}
\end{algorithm}


\section{Freely available Python code}
\label{sec:pycode}

The pseudocode described in \cref{alg:pred_limiter,alg:corr_limiter_1,alg:corr_limiter_2,alg:kriv_limiter,alg:full_method}
has been implemented into source code using the Python (version 2.7) programming language   \cite{article:python}, with use of the NumPy  library for basic mathematical functions and numerical
linear algebra tools \cite{book:numpy}, the PyPy  implementation of Python that makes use of just-in-time compilers to drastically improve runtime \cite{code:pypy}, and Matplotlib  for
visualization \cite{article:matplotlib}. The resulting code is publicly available for download \cite{code:Pythonic_LxW_DG}.

The numerical time-step as described in \cref{alg:full_method}, at least in principle, can be made arbitrarily high-order. For practical reasons, the implementation in \cite{code:Pythonic_LxW_DG}
is currently limited to orders of accuracy from $\morder=1$ to $\morder=5$. One limitation with
very high-order schemes is that the maximum Courant-Friedrichs-Lewy number
for which the scheme is still stable,
\begin{equation}
\text{CFL} = |\lambda|_{\text{max}} \, \frac{\Delta t}{\Delta x}, 
\end{equation}
{\it decreases} with increasing $\morder$. 
In the above expression, $|\lambda|_{\text{max}}$ is a bound on the maximum spectral radius of the flux Jacobian over the entire mesh and over the current time-step, $[t^n,t^n+\Delta t]$. The maximum allowable CFL number decreases roughly
as the inverse of $\morder$ (e.g., see \cite{article:GuRo2017} and references therein). 
We catalog the the CFL numbers used in the implementation of
\cite{code:Pythonic_LxW_DG} in \cref{table:cfl_numbers}.

The Pythonic Lax-Wendroff DG code \cite{code:Pythonic_LxW_DG} we developed 
is structured so that the top-level directory contains all of 
the application-specific sub-directories
(Burgers, shallow water, and Euler), as well as the {\tt lib} sub-directory that contains 
the main LxW-DG and plotting functions. In each application-specific sub-directory there 
are several specific numerical example sub-directories, as well as another
{\tt lib} directory that contains information about the fluxes and limiters for that particular
equation. Each numerical example sub-directory contains the following
three files that are required to run the main code:
\begin{enumerate}
\item \texttt{parameters.py}: \, set parameters for all values needed in the simulation;
\item \texttt{run\_example.py}: \, set initial conditions for simulation and execute main routine;
\item \texttt{plot\_example.py}: \, using Matplotlib, create plots of desired variables.
\end{enumerate}
Inside each numerical example sub-directory, the Pythonic Lax-Wendroff DG code can be run
by executing the following Makefile commands:
\begin{enumerate}
\item \texttt{make run}: \, executes code with parameters set in \texttt{parameters.py} and
initial conditions set in \texttt{run\_example.py};
\item \texttt{make plot}: \, executes main plotting routine with options set in \texttt{plot\_example.py}.
\end{enumerate}

\begin{table}[!t]
\begin{center}
\begin{Large}
\begin{tabular}{|c||c|c|c|c|c|}
\hline
{\normalsize $\morder$} & {\normalsize 1} & {\normalsize 2} & {\normalsize 3} & {\normalsize 4} & {\normalsize 5} \\
\hline\hline
{\normalsize $\text{CFL} = |\lambda|_{\text{max}} \frac{\Delta t}{\Delta x}$} & 
{\normalsize 0.90} & {\normalsize 0.30} & {\normalsize 0.14} & {\normalsize 0.10} & {\normalsize 0.06} \\
\hline
\end{tabular}
\caption{List of CFL numbers as used in actual simulations
 for locally-implicit Lax-Wendroff DG schemes of
various orders as implemented in \cite{code:Pythonic_LxW_DG}.
These CFL numbers are slightly reduced from the values on the linear stability boundary. Note that these CFL numbers decrease roughly as the inverse of the method order $\morder$.\label{table:cfl_numbers}}
\end{Large}
\end{center}
\end{table}


\section{Numerical examples}
\label{sec:numerical_examples}

In this section we demonstrate the accuracy and robustness of the proposed numerical scheme on several standard
test cases for the Burgers equation \cref{eqn:burgers}, the shallow water equations \cref{eqn:shllw}, and
the compressible Euler equations \cref{eqn:euler}.
    
\subsection{\it Burgers equation}
\label{subsec:burgers_example}
\begin{figure}[!t]
\centering
\begin{tabular}{cc}
(a)\includegraphics[width=0.43\textwidth,trim={1.0cm 1.8cm 1.3cm 1.8cm},clip]{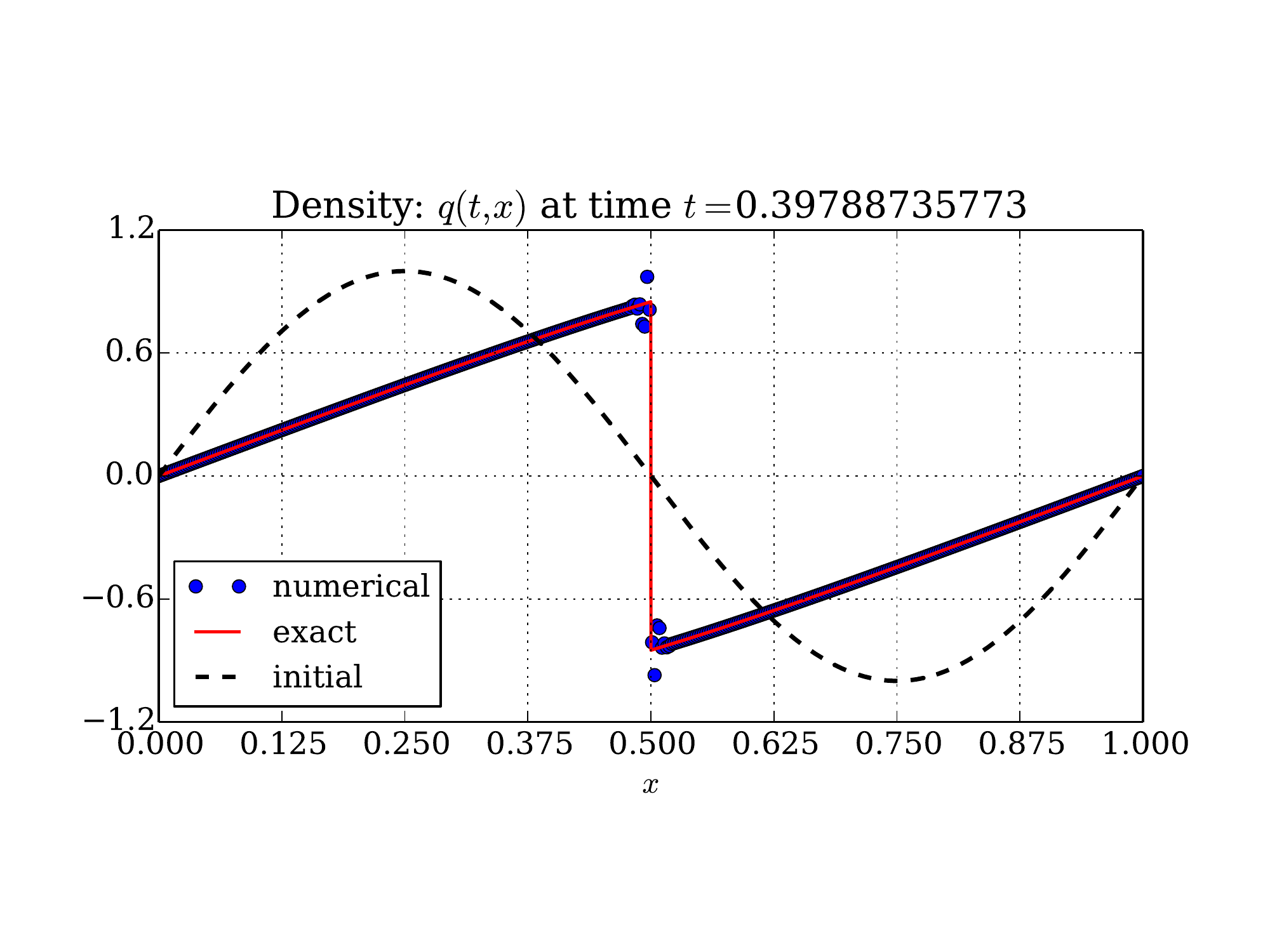} &
(b)\includegraphics[width=0.43\textwidth,trim={1.0cm 1.8cm 1.3cm 1.8cm},clip]{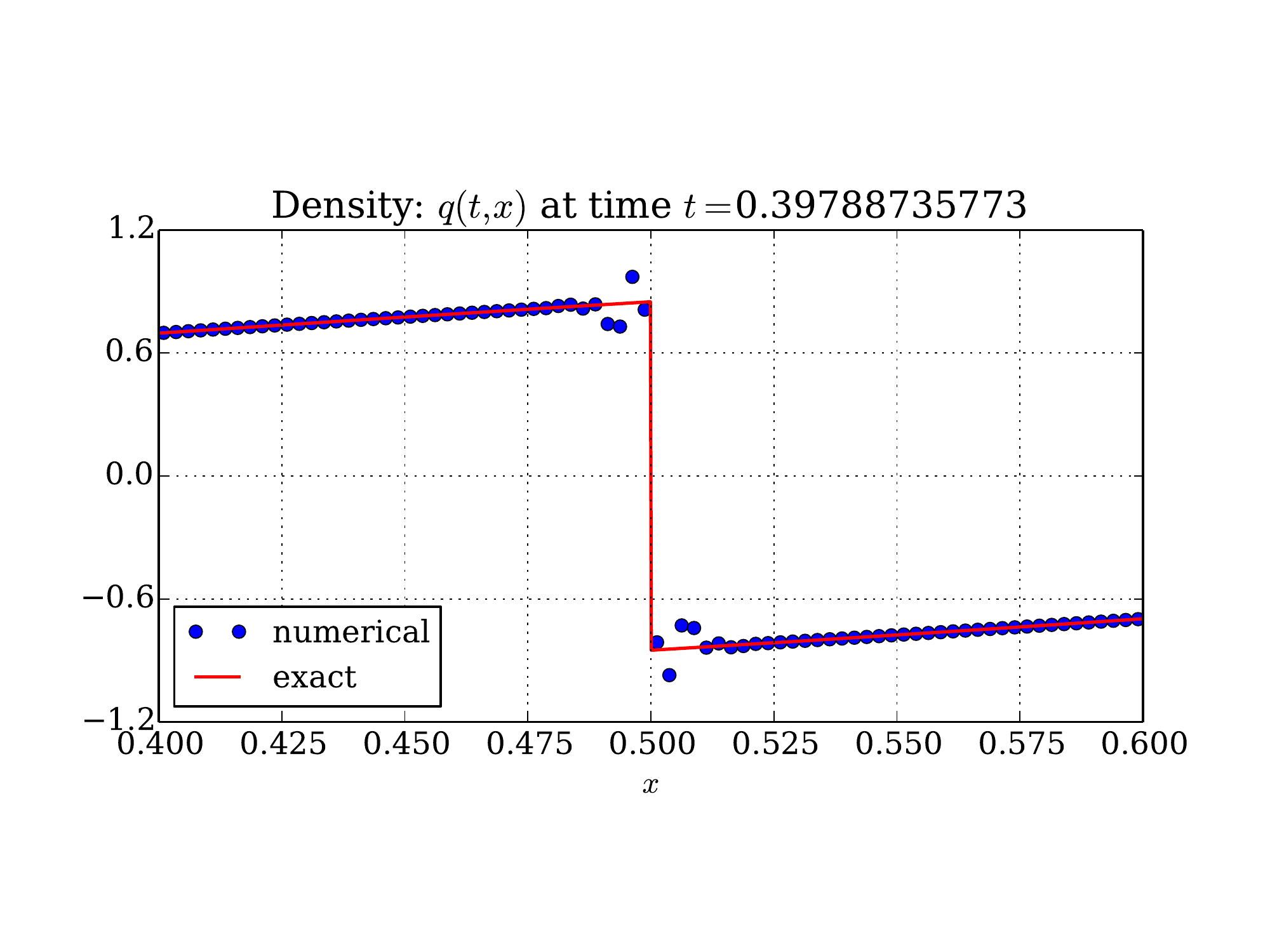} \\
(c)\includegraphics[width=0.43\textwidth,trim={1.0cm 1.8cm 1.3cm 1.8cm},clip]{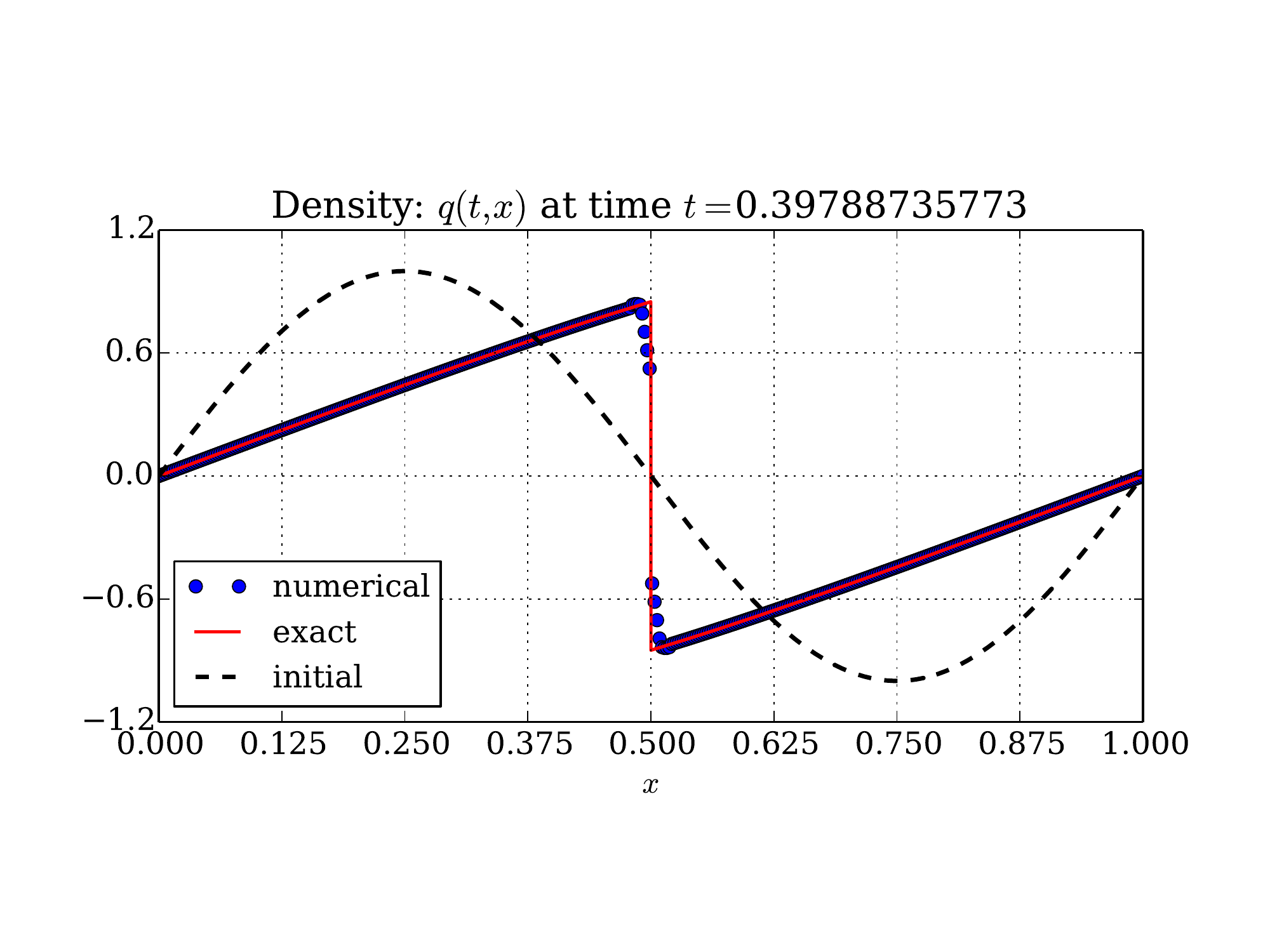} &
(d)\includegraphics[width=0.43\textwidth,trim={1.0cm 1.8cm 1.3cm 1.8cm},clip]{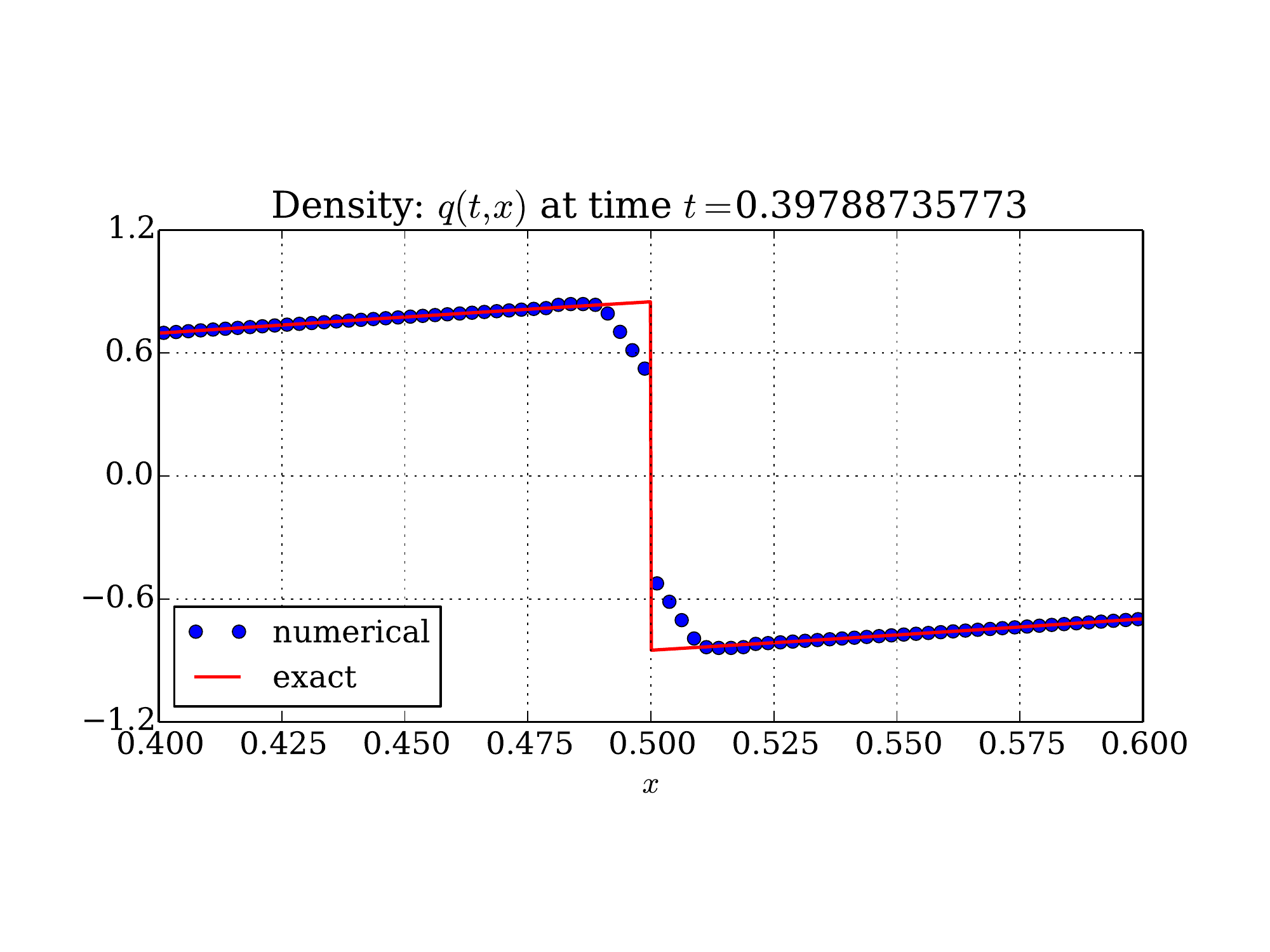}
\end{tabular}
\caption{Numerical and exact solution of Burgers equation \cref{eqn:burgers}. 
Shown are solutions computed with a smooth initial condition \cref{eqn:burgers_example} that forms a shock at
time \cref{eqn:burgers_example_shock_time}, periodic boundary conditions,
$\melems = 100$, and $\morder=4$ (i.e., fourth-order scheme with polynomial degree $\mdeg=3$). The solutions are shown at time $t=5(4\pi)^{-1} \approx 0.39788735773$, well after the shock has formed. 
In each panel we are plotting four points per element in order to clearly show
the subcell structure of the numerical solution. The individual panels show (a) the solution with no limiters,
(b) a zoomed-in version of this solution, (c) the solution with the correction step limiter, and 
(d) a zoomed-in version of this solution.
 \label{fig:burgers}}
\end{figure}

We begin by applying the proposed method to the Burgers equation \cref{eqn:burgers} with periodic boundary conditions
on $x \in [0,1]$ and the smooth initial data:
\begin{equation}
\label{eqn:burgers_example}
  q(t=0,x) = q_0(x) = \sin\left(2 \pi x \right).
\end{equation}
This initial condition forms a stationary shock at $x=0.5$ at time
\begin{equation}
\label{eqn:burgers_example_shock_time}
t_{\text{shock}} =  \left[ \max_{\xi \in \reals} \left\{ -2\pi\cos\left( 2 \pi \xi \right) \right\} \right]^{-1}
 = \left(2\pi \right)^{-1} \approx 0.15915494309.
\end{equation}
In this example, the only limiter that is active is our variant of the Krivodonova \cite{article:Kriv07} correction step limiter 
(see \cref{subsec:nonoscillatory}); and therefore, we can use this example to isolate its efficacy.  
We compare two simulations that are identical in every way except one: (1) 
in the first simulation no limiter is applied, and (2) in the second simulation the correction step limiter is used.
The results with $\melems = 100$, $\morder=4$, and at time $t=5(4\pi)^{-1} \approx 0.39788735773$ are shown
in \cref{fig:burgers}. The individual panels show (a) the solution with no limiters,
(b) a zoomed-in version of this solution, (c) the solution with the correction step limiter, and 
(d) a zoomed-in version of this solution. In each panel we are plotting four points per element in order to clearly show
the subcell structure of the numerical solution.

From these results we can draw some conclusions: (1) unsurprisingly, without a limiter the numerical solution is oscillatory,
producing unphysical overshoots and undershoots near the shock location; (2) the correction step limiter, which is applied only once per time-step (i.e., it is not applied in any of the prediction step iterations), successfully damps out the unphysical oscillations; and (3) the correction step limiter is not overly aggressive in that the damping effect is apparent only in the two elements on either side of the stationary shock.

\subsection{\it Shallow water equations}
We next consider solving the shallow water equations \cref{eqn:shllw} with various initial conditions that
demonstrate the accuracy of the scheme and the efficacy of the proposed limiting strategies.
In all examples the gravitational constant is taken as $g=1$ and and the limiter parameter is taken as
$\varepsilon = 10^{-14}$.

\subsubsection{Convergence test}
In order to verify the order of accuracy of the proposed scheme we make use of the so-called 
{\it method of manufactured solutions} (e.g., see \cite{book:KnSa02}), where we prescribe a ``solution'',
and then add a source term to our original PDE that guarantees that this ``solution'' satisfies the PDE with 
the additional source term. The additional source term will depend explicitly on time and space, but not on
the conserved variables $q$.

We consider an example where the manufactured solution is given by 
\begin{equation}
\label{eqn:meth_man_solns}
h_{\text{MS}}(t,x) := 1 + \frac{1}{2} \sin \left(\pi (x-t) \right), \quad
u_{\text{MS}}(t,x) := \cos\left(2 \pi (x-2t) \right),
\end{equation}
on the domain $[-1,1]$ with periodic boundary conditions.
The shallow water equations, with an additional source term to guarantee that \cref{eqn:meth_man_solns}
is indeed a solution, can be written as follows:
\begin{equation}
\label{eqn:forced_shllw}
\begin{bmatrix}[1.5] h \\ hu \end{bmatrix}_{,t} + 
\begin{bmatrix}[1.5] hu \\ hu^2 + \frac{1}{2} g h^2 \end{bmatrix}_{,x} = \begin{bmatrix}[1.5]
\left( h_{\text{MS}} \right)_{,t} + \left( h_{\text{MS}} u_{\text{MS}} \right)_{,x} \\ 
\left(h_{\text{MS}} u_{\text{MS}} \right)_{,t} + \left( h_{\text{MS}} u^2_{\text{MS}} + \frac{1}{2} g h_{\text{MS}}^2\right)_{,x}\end{bmatrix}.
\end{equation}

We compute numerical solutions to \cref{eqn:forced_shllw} with \cref{eqn:meth_man_solns},
on $[-1,1]$, with periodic boundary conditions, with $\melems = N = 10 \times 2^k$ elements for $k=0,1,2,\ldots,5$,
to time $t=0.5$, and with three different orders of accuracy: $\morder=3,4,5$. 
For each simulation, we aim to compute the relative  $L^2[-1,1]$ error:
\begin{equation}
\text{rel.} \, \,  L^2[-1,1] \, \text{error} := \sum_{\ell=1}^{\meq}
\sqrt{\frac{\int_{-1}^{1} \left| q^h_{\ell}(0.5,x) - q^{\star}_{\ell}(0.5,x) \right|^2 \, dx}{\int_{-1}^{1} \left| q^{\star}_{\ell}(0.5,x) \right|^2 \, dx}},
\end{equation}
where $q^h$ and $q^{\star}$ are the numerical and exact solutions, respectively.
In practice, we replace the exact solution by a piecewise Legendre polynomial approximation of
degree $\morder+1$, which allows us to obtain the following approximate relative $L^2[-1,1]$ error:
\begin{equation}
\label{eqn:rel_error}
e_N := \sum_{\ell=1}^{\meq} \sqrt{ \frac{\sum\limits_{i=1}^{N} \left( \sum\limits_{k=1}^{\mcorr} \left( Q_{i \, (k,\ell)} - 
Q^{\star}_{i \, (k,\ell)} \right)^2 + \left( Q^{\star}_{i \, (\mcorr+1,\ell)} \right)^2 \right)}{
\sum\limits_{i=1}^{N} \sum\limits_{k=1}^{\mcorr+1} \left( Q^{\star}_{i \, (k,\ell)} \right)^2 } },
\end{equation}
where $Q^{\star}$ are the Legendre coefficients of the exact solution \cref{eqn:meth_man_solns}, written
as conservative variables at
time $t=0.5$. The exact solution coefficients are computed via Gaussian quadrature
with $\morder+1$ points:
\begin{equation}
\vec{Q^{\star}_{i \, (k,:)}} := \frac{1}{2} \sum_{a=1}^{\morder+1} \omega^{\star}_a \, \vec{q_{\text{MS}}}\left(0.5,x_i
+ \frac{\Delta x}{2} \mu^{\star}_a \right),
\end{equation}
where $\omega^{\star}_a$ and $\mu_a^{\star}$ for $a=1,\ldots,\morder+1$ are the weights and abscissas
of the $(\morder+1)$-point Gaussian quadrature rule. The resulting errors and error ratios are catalogued
in \cref{table:shllw_1d_error}.

\begin{table}
\begin{center}
\begin{Large}
\begin{tabular}{|c||c|c||c|c||c|c|}
\hline
{\normalsize $N$} & {\normalsize $\text{e}_{N} (\morder=3)$} & {\normalsize $\log_2\frac{\text{e}_{N/2}}{\text{e}_{N}}$}  & {\normalsize $\text{e}_N (\morder=4)$} & {\normalsize $\log_2\frac{\text{e}_{N/2}}{\text{e}_{N}}$}  & {\normalsize $\text{e}_N (\morder=5)$} & {\normalsize $\log_2\frac{\text{e}_{N/2}}{\text{e}_{N}}$} \\
\hline\hline
{\normalsize 10} & {\normalsize 1.990e-02} & -- & {\normalsize 1.324e-03} & -- & {\normalsize 1.104e-04} & -- \\\hline
{\normalsize 20} & {\normalsize 2.409e-03} & {\normalsize $3.047$} & {\normalsize 7.650e-05} & {\normalsize $4.114$} & {\normalsize 3.216e-06} & {\normalsize $5.102$} \\\hline
{\normalsize 40} & {\normalsize 3.183e-04} & {\normalsize $2.920$} & {\normalsize 4.497e-06} & {\normalsize $4.088$} & {\normalsize 1.090e-07} & {\normalsize $4.883$} \\\hline
{\normalsize 80} & {\normalsize 4.210e-05} & {\normalsize $2.918$} & {\normalsize 2.792e-07} & {\normalsize $4.010$} & {\normalsize 3.679e-09} & {\normalsize $4.889$} \\\hline
{\normalsize 160} & {\normalsize 5.563e-06} & {\normalsize $2.920$} & {\normalsize 1.768e-08} & {\normalsize $3.981$} & {\normalsize 1.199e-10} & {\normalsize $4.939$} \\\hline
{\normalsize 320} & {\normalsize 7.341e-07} & {\normalsize $2.922$} & {\normalsize 1.126e-09} & {\normalsize $3.972$} & {\normalsize 3.846e-12} & {\normalsize $4.963$} \\\hline
\end{tabular} 
\caption{Relative $L^2$ errors for the forced shallow water equations with periodic boundary conditions.
\label{table:shllw_1d_error}}
\end{Large}
\end{center}
\end{table}

\subsubsection{Shock formation in finite time}
In \cref{subsec:burgers_example} we considered an example for the Burgers equation that started as a smooth initial condition and then shocked in finite time. We now attempt similar example, this time for the shallow
water equations. We take the following initial conditions on the computational domain $[-1,1]$ with outflow
boundary conditions at $x=\pm 1$:
\begin{equation}
\label{eqn:shllw_gauss_example}
h(t=0,x) = 1 + e^{-100x^2}, \quad u(t=0,x) = 0.
\end{equation}
We note that the subsequent dynamics will result in the initial Gaussian bump splitting into
two smaller height disturbances, one propagating to the left, the other to the right. The tops of these height disturbances will propagate faster than the rest of the profile, which will lead to wave steepening and eventually to
 the formation of a shock in both the left and right propagating disturbances.

\begin{figure}[!t]
\centering
\begin{tabular}{cc}
(a)\includegraphics[width=0.43\textwidth,trim={1.0cm 1.8cm 1.3cm 1.8cm},clip]{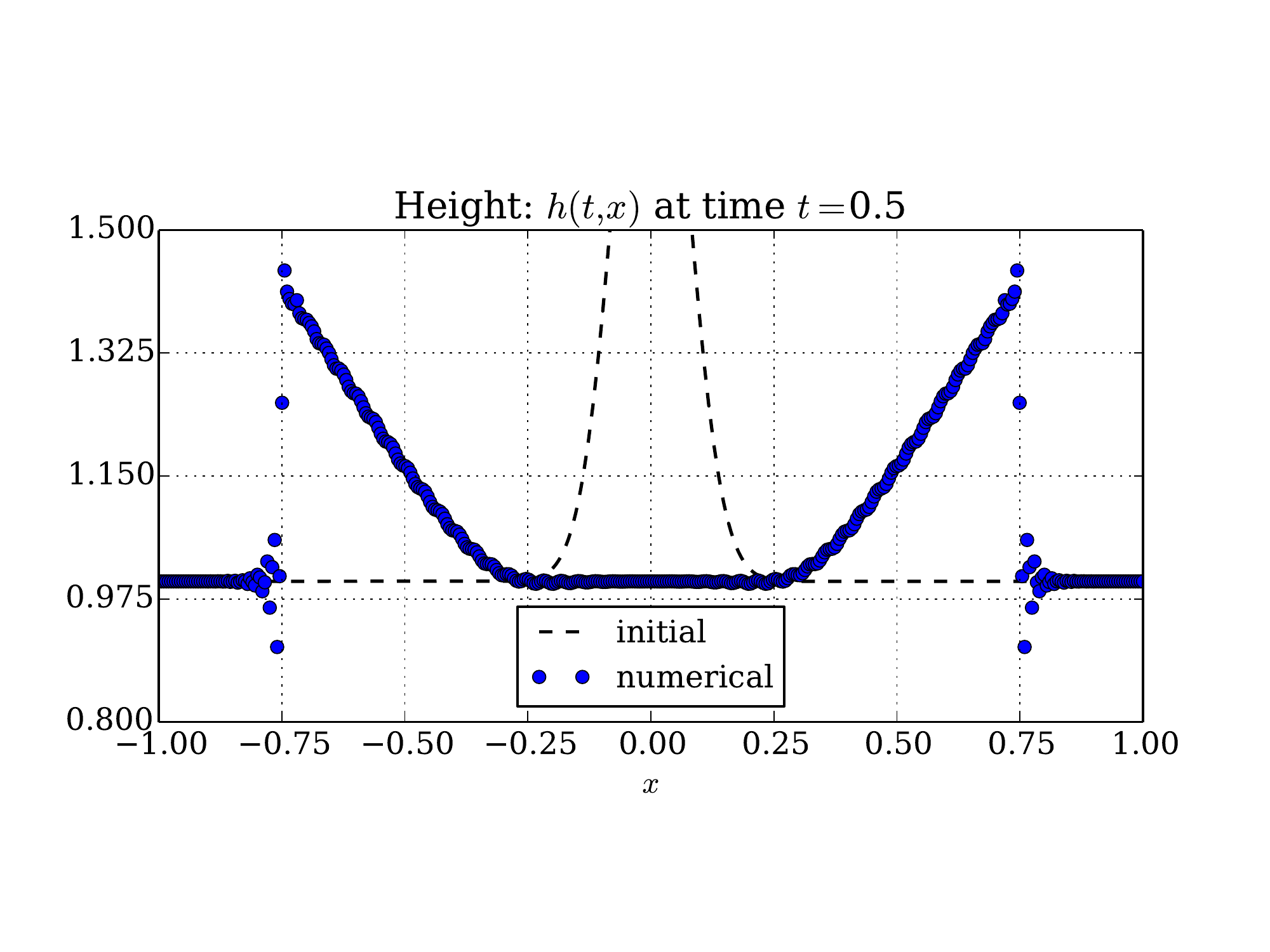} &
(b)\includegraphics[width=0.43\textwidth,trim={1.0cm 1.8cm 1.3cm 1.8cm},clip]{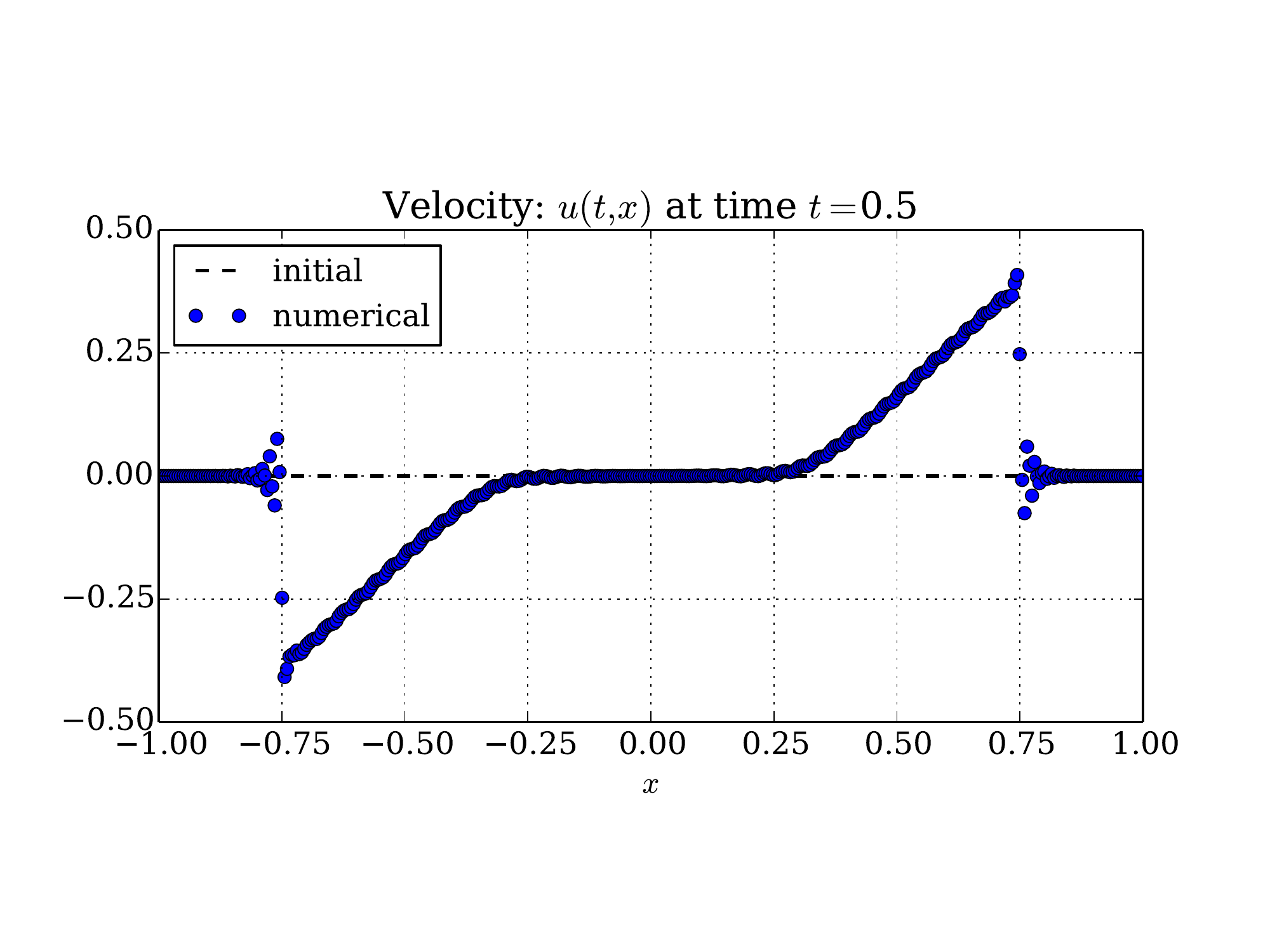} \\
(c)\includegraphics[width=0.43\textwidth,trim={1.0cm 1.8cm 1.3cm 1.8cm},clip]{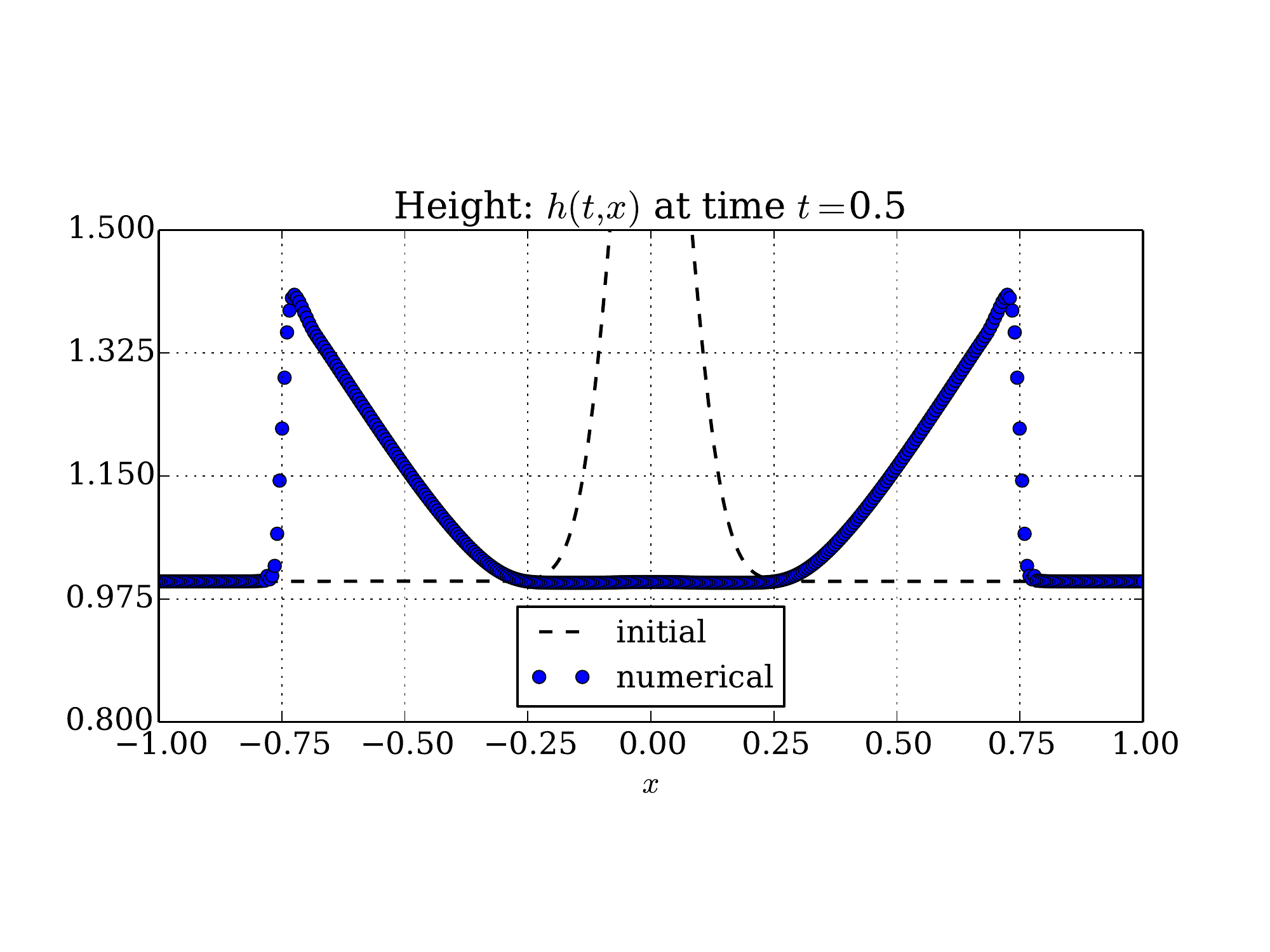} &
(d)\includegraphics[width=0.43\textwidth,trim={1.0cm 1.8cm 1.3cm 1.8cm},clip]{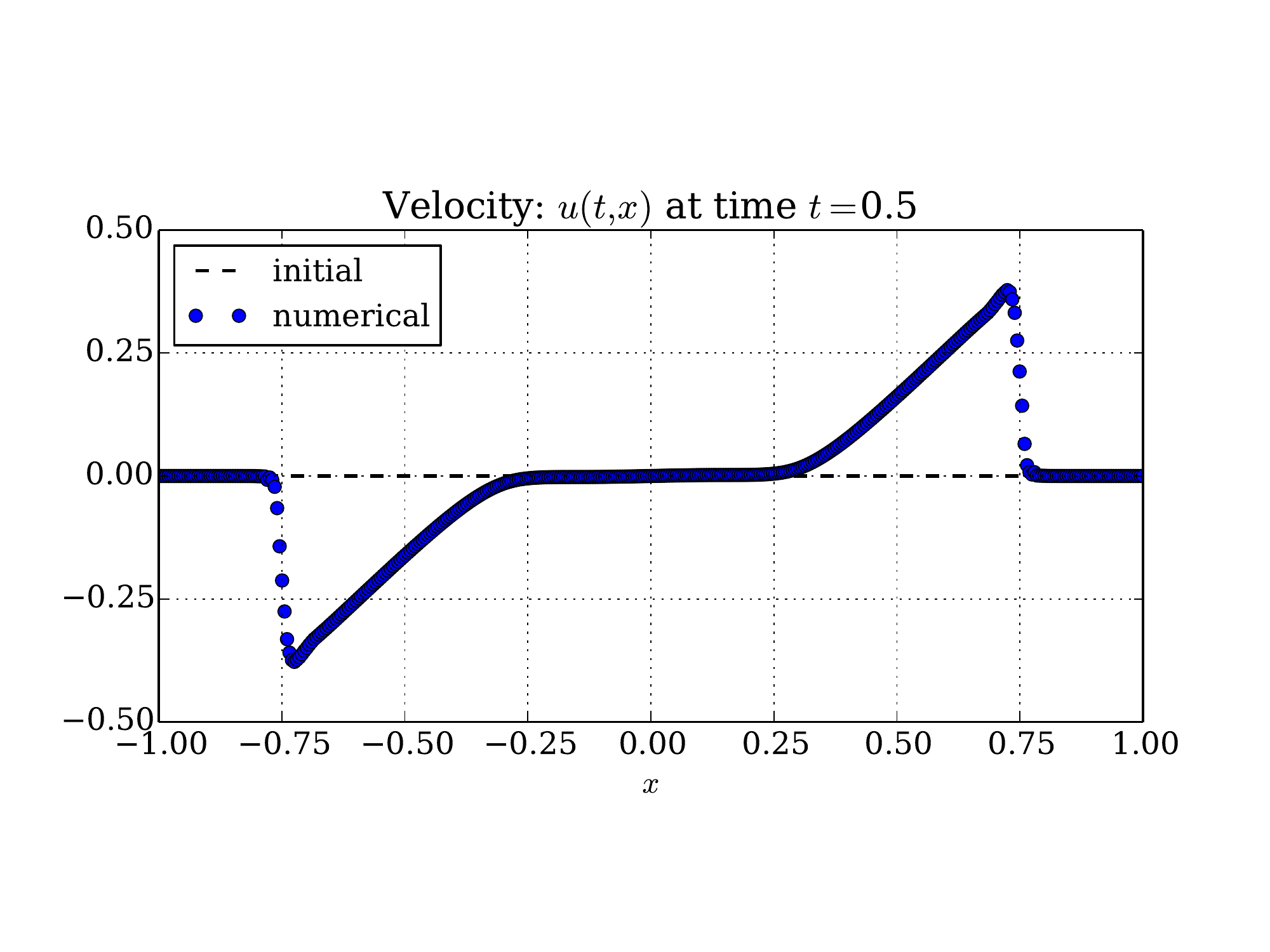}
\end{tabular}
\caption{Comparison of the unlimited and limited schemes on the shallow water equations \cref{eqn:shllw}. 
Shown are solutions computed with an initially smooth height profile and zero velocity,
\cref{eqn:shllw_gauss_example}, that forms shocks in finite time, outflow boundary conditions,
$\melems = 100$, and $\morder=4$ (i.e., fourth-order scheme with polynomial degree $\mdeg=3$). The solutions are shown at time $t=0.5$, after the shocks have formed. 
In each panel we are plotting four points per element in order to clearly show
the subcell structure of the numerical solution. The individual panels show the (a) height from 
the unlimited scheme, (b) velocity from the unlimited scheme, (c) height from the scheme with limiters,
and (d) velocity from the scheme with limiters. 
 \label{fig:shllw_gaussian}}
\end{figure}

The results of solving the shallow water equations with initial conditions \cref{eqn:shllw_gauss_example}, with the 
$\morder=4$ version of the scheme, and with $\melems=100$ elements, is shown in \cref{fig:shllw_gaussian}.
 As in \cref{subsec:burgers_example}, we run this problem
with and without limiters in order to demonstrate efficacy of the limiting procedure.
The individual panels show at time $t=0.5$: (a) the height and (b) velocity as computed with the unlimited scheme,
and (c) the height and (d) velocity as computed with the limited scheme.
In each panel we are plotting four points per element in order to clearly show
the subcell structure of the numerical solution. In both cases, the numerical solution is far from
exhibiting positivity violations; and therefore, the only limiter that is active is the 
correction step limiter (see \cref{subsec:nonoscillatory}). We again see the ability of this once-per-time-step
limiter to simultaneously remove unphysical oscillations and to maintain the sharpness of the solution.

\subsubsection{Dambreak problem}
Next we consider an example of a Riemann problem, which for shallow water equations is also called the
dambreak problem. The initial conditions are piecewise constant:
\begin{equation}
\label{eqn:shllw_dambreak}
\left( h, \, u \right)(t=0,x) = 
\begin{cases}
   \left( 1.0, \, 0 \right) & \quad x<0, \\
   \left( 0.1, \, 0 \right) & \quad x>0.
\end{cases}
\end{equation}
The entropy-satisfying solution with this initial data is a left-propagating rarefaction wave and a right-propagating
shock. A full mathematical explanation of the solution of this Riemann problem can be found in several textbooks, including Chapter 13 of LeVeque \cite{book:Le02}.

\begin{figure}[!t]
\centering
\begin{tabular}{cc}
(a)\includegraphics[width=0.43\textwidth,trim={0.7cm 1.8cm 1.3cm 1.8cm},clip]{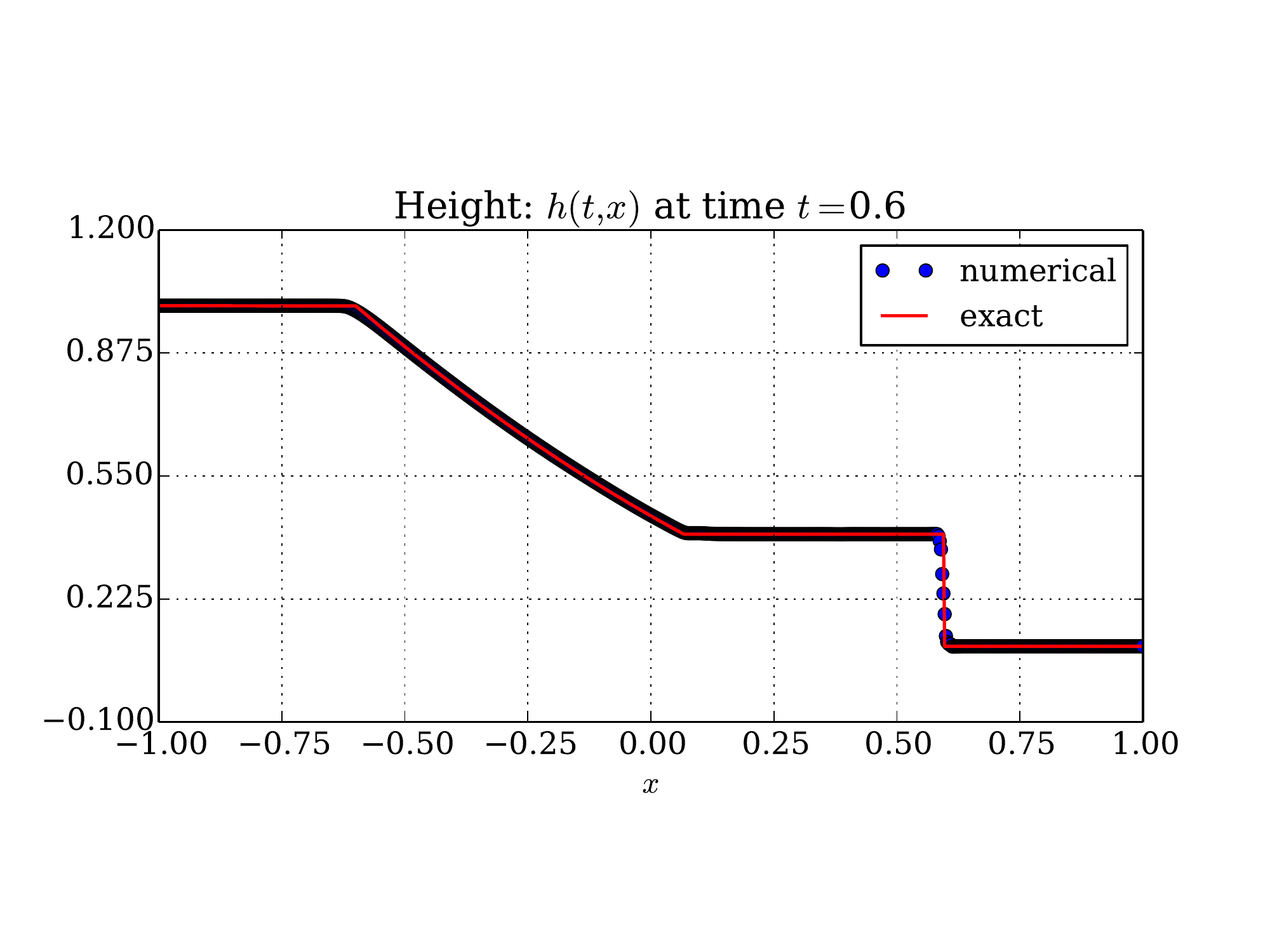} &
(b)\includegraphics[width=0.43\textwidth,trim={0.7cm 1.8cm 1.3cm 1.8cm},clip]{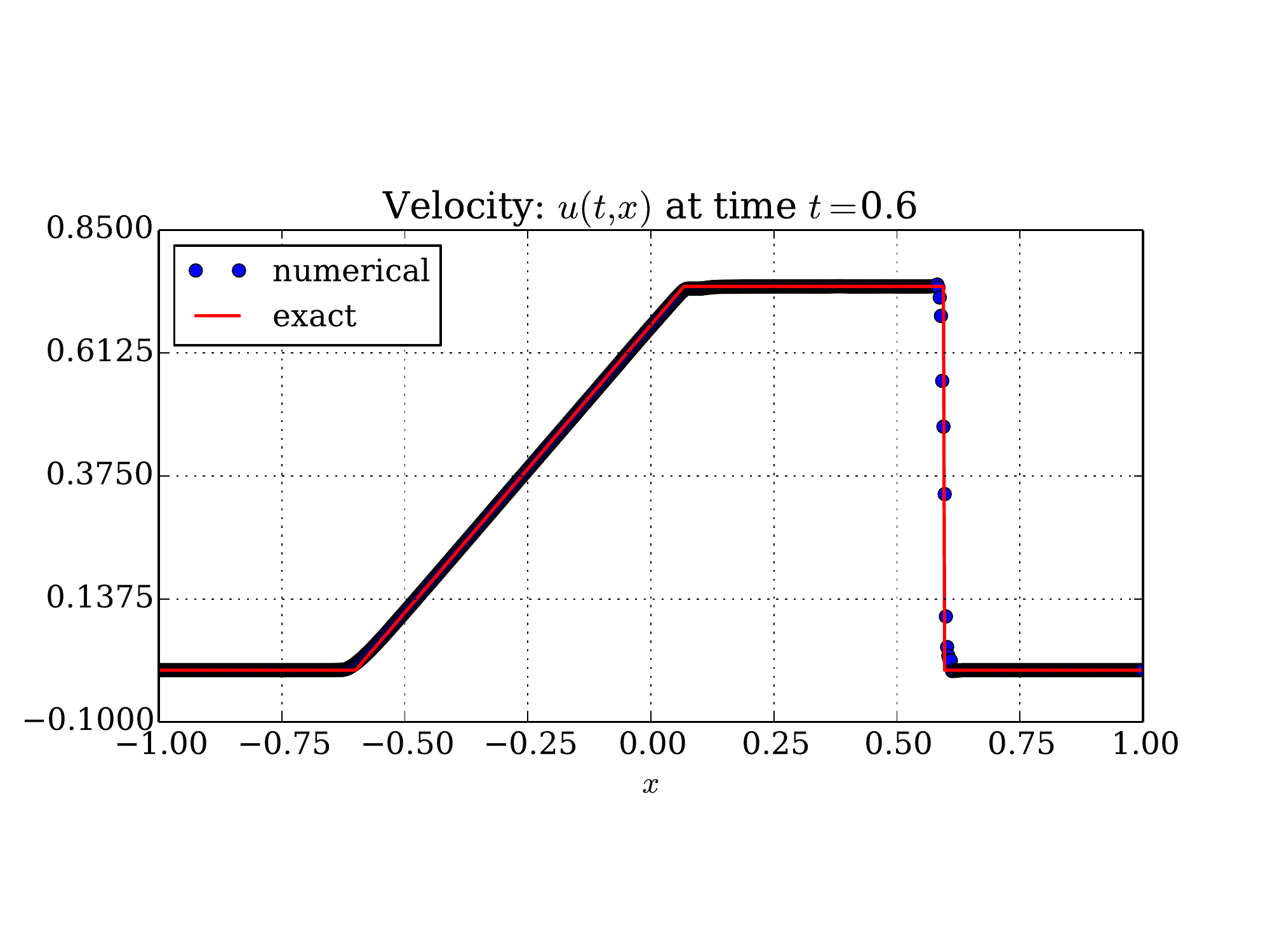}
\end{tabular}
\caption{
Dambreak Riemann problem for the shallow water equations \cref{eqn:shllw}. 
Shown are solutions computed with initial conditions
\cref{eqn:shllw_dambreak}, outflow boundary conditions,
$\melems = 200$, and $\morder=4$ (i.e., fourth-order scheme with polynomial degree $\mdeg=3$). The solutions are shown at time $t=0.4$. 
In each panel we are plotting four points per element in order to clearly show
the subcell structure of the numerical solution. The individual panels show the numerical solution and the superimposed exact Riemann solution for the (a) height and
(b) fluid velocity. 
\label{fig:shllw_dambreak}}
\end{figure}

The results of solving the shallow water equations with initial conditions \cref{eqn:shllw_dambreak}, on
the domain $[-1,1]$ with outflow boundary conditions, with the 
$\morder=4$ version of the scheme, and with $\melems=200$ elements is shown in \cref{fig:shllw_dambreak}.
The individual panels show at time $t=0.6$: (a) the height and (b) the velocity,
 with the exact Riemann solution superimposed.
In each panel we are plotting four points per element in order to clearly show
the subcell structure of the numerical solution.
As in the previous example, the numerical solution is far from
exhibiting positivity violations; and therefore, the only limiter that is active is the 
correction step limiter (see \cref{subsec:nonoscillatory}). We again see the ability of this once-per-time-step
limiter to simultaneously remove unphysical oscillations and to maintain the sharpness of the solution.

\subsubsection{Double rarefaction}
\label{subsubsec:shllw_double_raref}
In order to test the positivity limiters on the shallow water equations, we 
attempt a different Riemann problem, this time with initial conditions
that result in two counter-propagating rarefactions that leave in their wake a near-vacuum
state. The initial conditions are as follows:
\begin{equation}
\label{eqn:shllw_double_raref}
\left( h, \, u  \right)(t=0,x) = 
\begin{cases}
   \left( 1, \, -2 \right) & \quad x<0, \\
   \left( 1, \, +2 \right) & \quad x>0.
\end{cases}
\end{equation}
A full mathematical explanation of the solution of this Riemann problem can be found in several textbooks, including in Chapter 13 of LeVeque \cite{book:Le02}.

The results of solving the  shallow water equations with initial conditions \cref{eqn:shllw_double_raref}, on
the domain $[-1,1]$ with outflow boundary conditions, with the 
$\morder=4$ version of the scheme, and with $\melems=200$ elements is shown in \cref{fig:euler_doublerarefaction}.
The individual panels show at time $t=0.25$: (a) the full height profile, (b) a zoomed-in view of the height in the
near-vacuum region, (c) the full momentum profile, and (d) a zoomed-in view of the 
momentum. The exact Riemann solution is superimposed.
In each panel we are plotting four points per element in order to clearly show
the subcell structure of the numerical solution.
The proposed limiters are able to handle the near-vacuum solution in the center of the computational domain,
and the numerical solution remains stable and positivity-preserving in the sense of \cref{eqn:shllw_pos_set_1}
and \cref{eqn:shllw_pos_set_2}.

\begin{figure}[!t]
\centering
\begin{tabular}{cc}
(a)\includegraphics[width=0.43\textwidth,trim={1.0cm 1.8cm 1.3cm 1.8cm},clip]{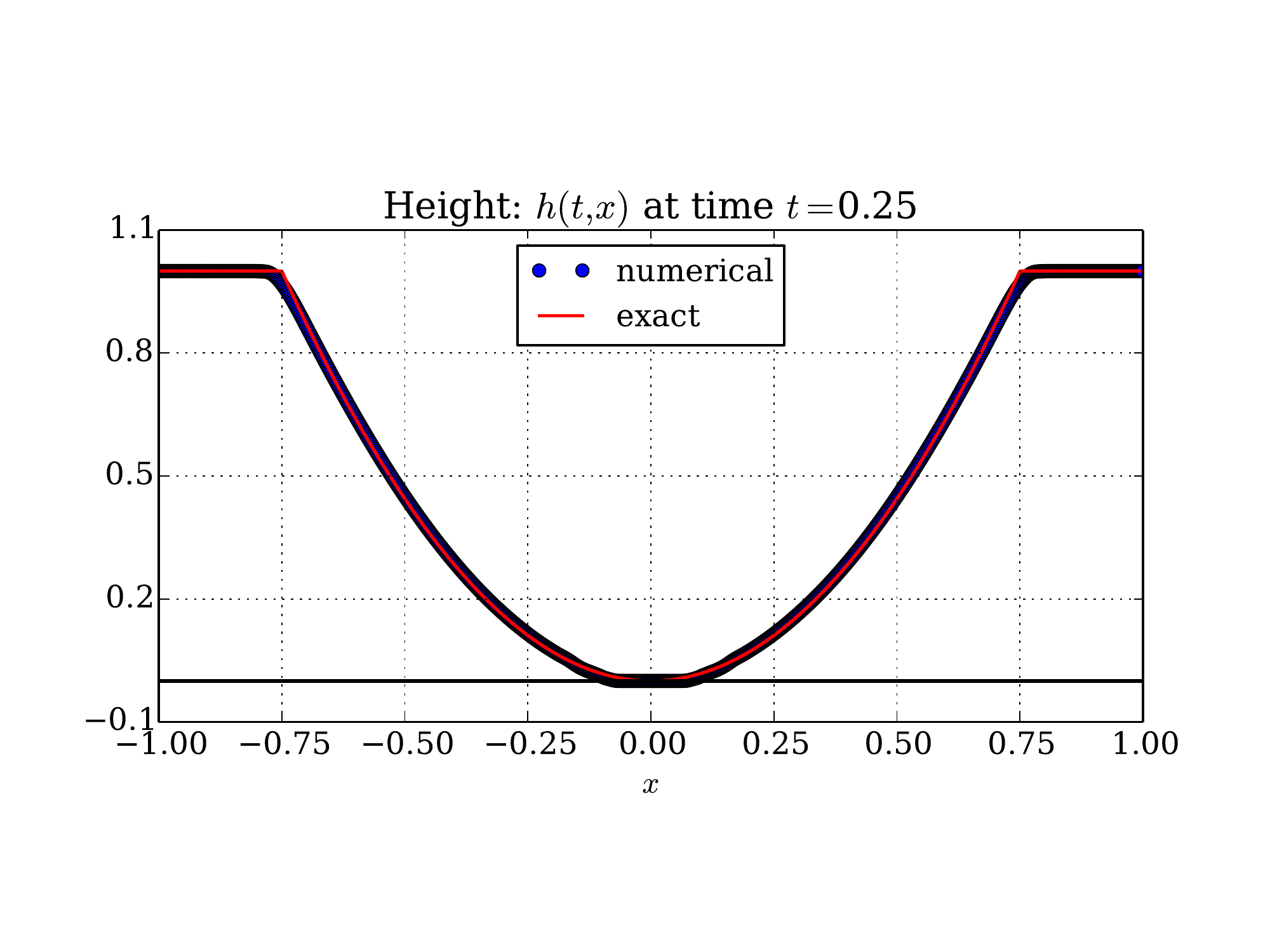} &
(b)\includegraphics[width=0.43\textwidth,trim={1.0cm 1.8cm 1.3cm 1.8cm},clip]{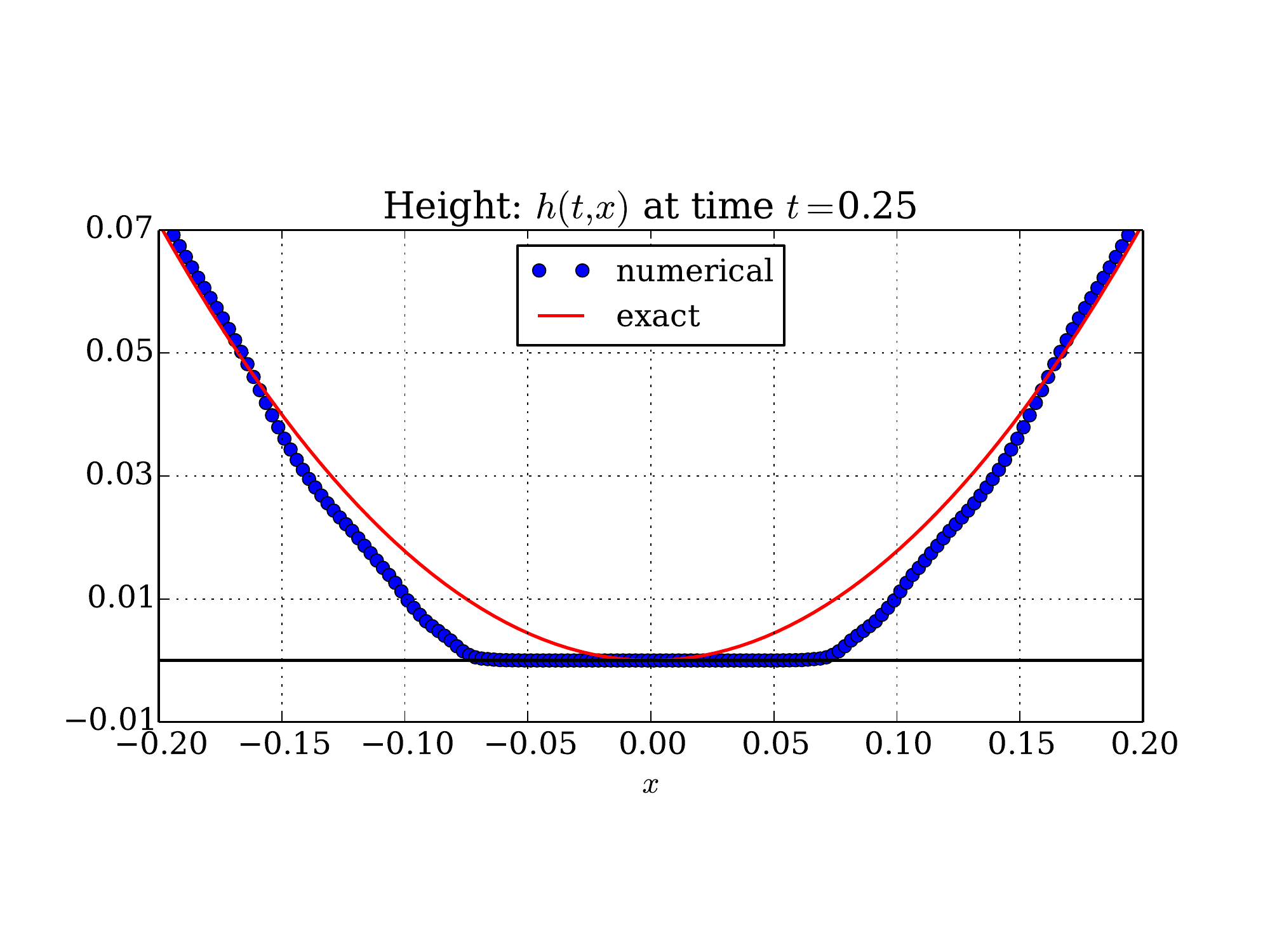} \\
(c)\includegraphics[width=0.43\textwidth,trim={1.0cm 1.8cm 1.3cm 1.8cm},clip]{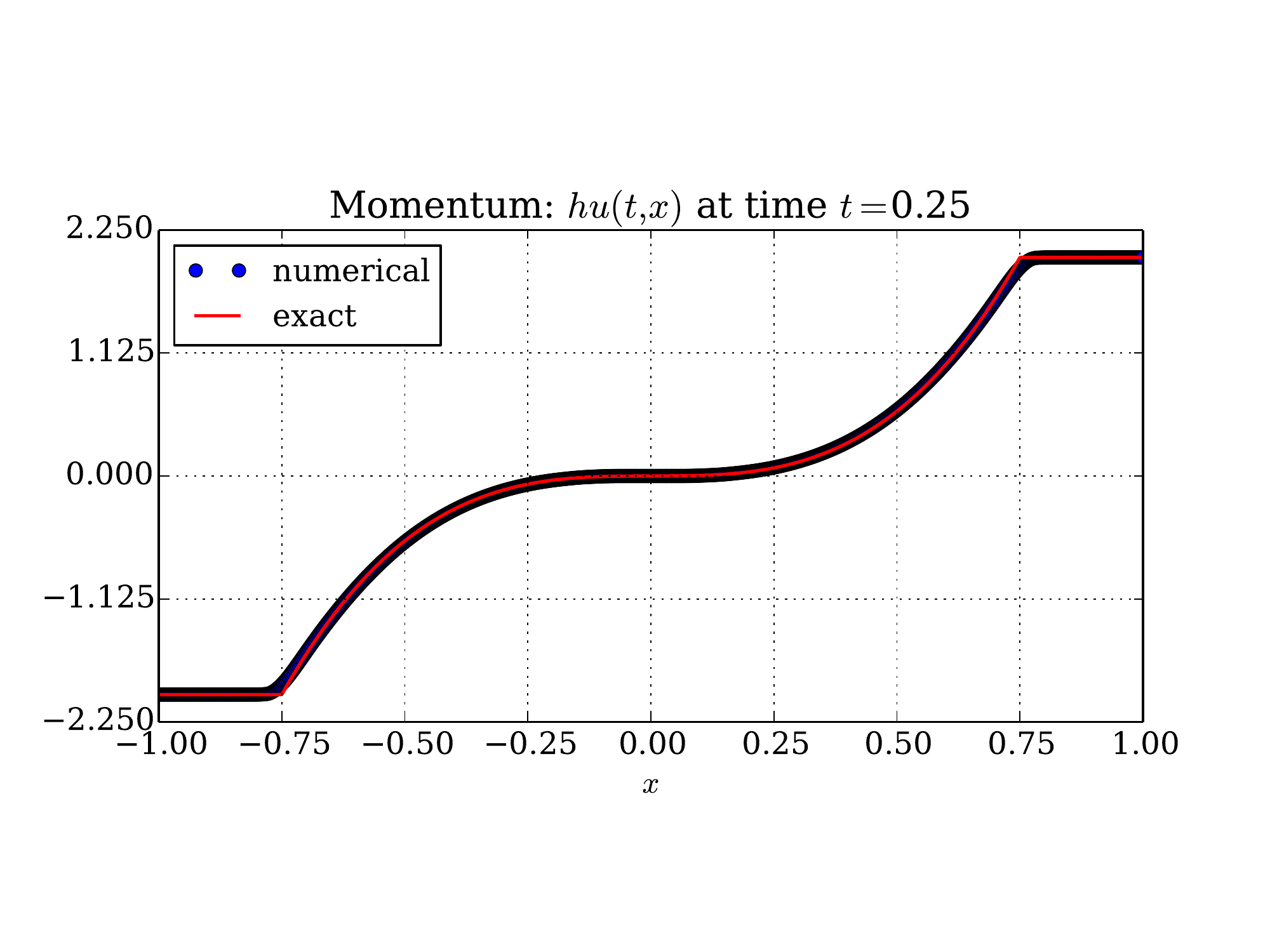} &
(d)\includegraphics[width=0.43\textwidth,trim={1.0cm 1.8cm 1.3cm 1.8cm},clip]{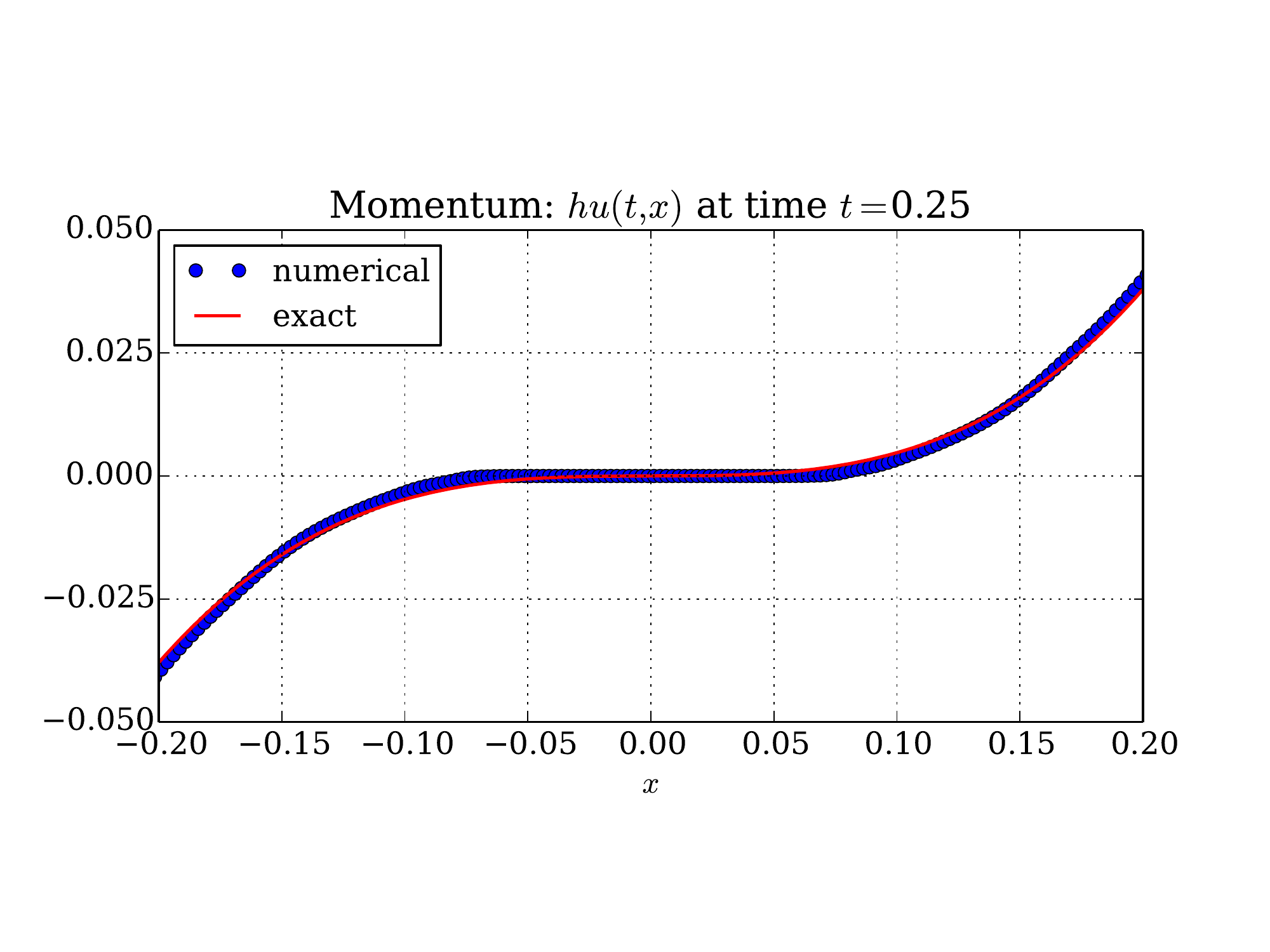}
\end{tabular}
\caption{Double rarefaction Riemann problem for the shallow water equation \cref{eqn:shllw}.
Shown are solutions computed with initial conditions
\cref{eqn:shllw_double_raref}, outflow boundary conditions,
$\melems = 200$, and $\morder=4$ (i.e., fourth-order scheme with polynomial degree $\mdeg=3$). 
The individual panels show at time $t=0.25$: (a) the full height profile, (b) a zoomed-in view of the height in the
near-vacuum region, (c) the full momentum profile, and (d) a zoomed-in view of the 
momentum. The exact Riemann solution is superimposed.
In each panel we are plotting four points per element in order to clearly show
the subcell structure of the numerical solution.
\label{fig:shllw_doublerarefaction}}
\end{figure}

\subsection{\it Compressible Euler equations}
Finally we consider solving the compressible Euler equations \cref{eqn:euler} with various initial conditions that
demonstrate the accuracy of the scheme and the efficacy of the proposed limiting strategies.
In all examples the adiabatic constant is taken as $\gamma=1.4$ and the limiter parameter is taken as
$\varepsilon = 10^{-14}$.

\subsubsection{Convergence test}
The compressible Euler equations, along with suitable boundary conditions,
admit a class of non-trivial and non-stationary
exact solutions in which the fluid velocity and pressure remain globally constant, but
an arbitrary density profile is advected by the fluid velocity.
We will use such a solution to verify the order of accuracy of the proposed scheme.
In particular, we consider the following exact solution on $[-1,1]$ with periodic boundary conditions:
\begin{equation}
\label{eqn:euler_smooth_ic}
\rho(t,x) = 1 + 0.5 \sin\left(3 \pi \left( x - 0.5 t \right) \right), \quad
u(t,x) = 0.5, \quad p(t,x) = 0.75.
\end{equation}

We compute numerical solutions to \cref{eqn:euler} with initial conditions
obtained from \cref{eqn:euler_smooth_ic},
on $[-1,1]$, with periodic boundary conditions, with $\melems = N = 10 \times 2^k$ elements for $k=0,1,2,\ldots,5$, to time $t=1$, and with three different orders of accuracy: $\morder=3,4,5$. For each simulation we compute
the relative $L^2[-1,1]$ error via \cref{eqn:rel_error}, where the $Q^{\star}$ coefficients
are computed from the conservative variable version of \cref{eqn:euler_smooth_ic} at time $t=1$.
 The resulting errors and error ratios are catalogued
in \cref{table:euler_1d_error}.

\begin{table}
\begin{center}
\begin{Large}
\begin{tabular}{|c||c|c||c|c||c|c|}
\hline
{\normalsize $N$} & {\normalsize $\text{e}_N (\morder=3)$} & {\normalsize $\log_2\frac{\text{e}_{N/2}}{\text{e}_{N}}$}  & {\normalsize $\text{e}_N (\morder=4)$} & {\normalsize $\log_2\frac{\text{e}_{N/2}}{\text{e}_{N}}$}  & {\normalsize $\text{e}_N (\morder=5)$} & {\normalsize $\log_2\frac{\text{e}_{N/2}}{\text{e}_{N}}$} \\
\hline\hline
{\normalsize 10} & {\normalsize 2.161e-02} & -- & {\normalsize 3.109e-03} & -- & {\normalsize 2.179e-04} & -- \\\hline
{\normalsize 20} & {\normalsize 3.742e-03} & {\normalsize $2.530$} & {\normalsize 1.225e-04} & {\normalsize $4.665$} & {\normalsize 1.010e-05} & {\normalsize $4.431$} \\\hline
{\normalsize 40} & {\normalsize 6.540e-04} & {\normalsize $2.517$} & {\normalsize 7.182e-06} & {\normalsize $4.093$} & {\normalsize 4.438e-07} & {\normalsize $4.509$} \\\hline
{\normalsize 80} & {\normalsize 9.633e-05} & {\normalsize $2.763$} & {\normalsize 4.398e-07} & {\normalsize $4.029$} & {\normalsize 1.623e-08} & {\normalsize $4.773$} \\\hline
{\normalsize 160} & {\normalsize 1.279e-05} & {\normalsize $2.913$} & {\normalsize 2.728e-08} & {\normalsize $4.011$} & {\normalsize 5.343e-10} & {\normalsize $4.925$} \\\hline
{\normalsize 320} & {\normalsize 1.629e-06} & {\normalsize $2.973$} & {\normalsize 1.706e-09} & {\normalsize $3.999$} & {\normalsize 1.695e-11} & {\normalsize $4.979$} \\\hline
\end{tabular} 
\caption{Relative $L^2$ errors for the compressible Euler equations with constant pressure and fluid velocity and with periodic boundary conditions.
\label{table:euler_1d_error}}
\end{Large}
\end{center}
\end{table}

\subsubsection{Shock tube problem}
Next we consider an example of a Riemann problem, which for compressible Euler equations is also called the
shock tube problem. We consider the celebrated Sod shock tube problem \cite{article:Sod78}, for
which the initial conditions are
\begin{equation}
\label{eqn:euler_shocktube}
\left( \rho, \, u, \, p \right)(t=0,x) = 
\begin{cases}
   \left( 1.000, \, 0, \, 1.0 \right) & \quad x<0, \\
   \left( 0.125, \, 0, \, 0.1 \right) & \quad x>0.
\end{cases}
\end{equation}
The entropy-satisfying solution with this initial data is a left-propagating rarefaction wave,
a right-propagating contact wave, and a (faster) right-propagating
shock. A full mathematical explanation of the solution of this Riemann problem can be found in several textbooks, including Chapter 14 of LeVeque \cite{book:Le02}.

\begin{figure}[!t]
\centering
\begin{tabular}{cc}
(a)\includegraphics[width=0.43\textwidth,trim={1.0cm 1.8cm 1.3cm 1.8cm},clip]{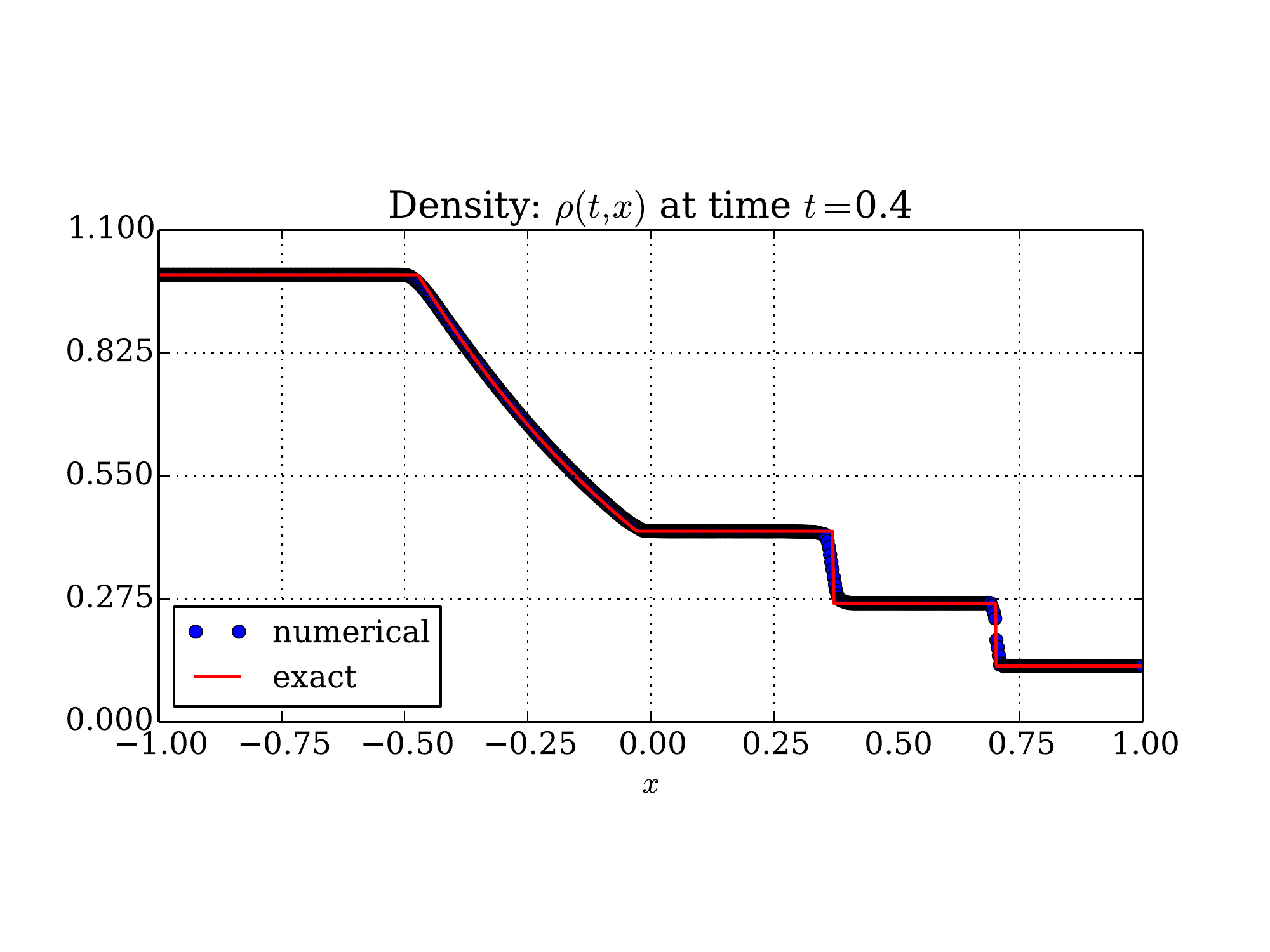} &
(b)\includegraphics[width=0.43\textwidth,trim={1.0cm 1.8cm 1.3cm 1.8cm},clip]{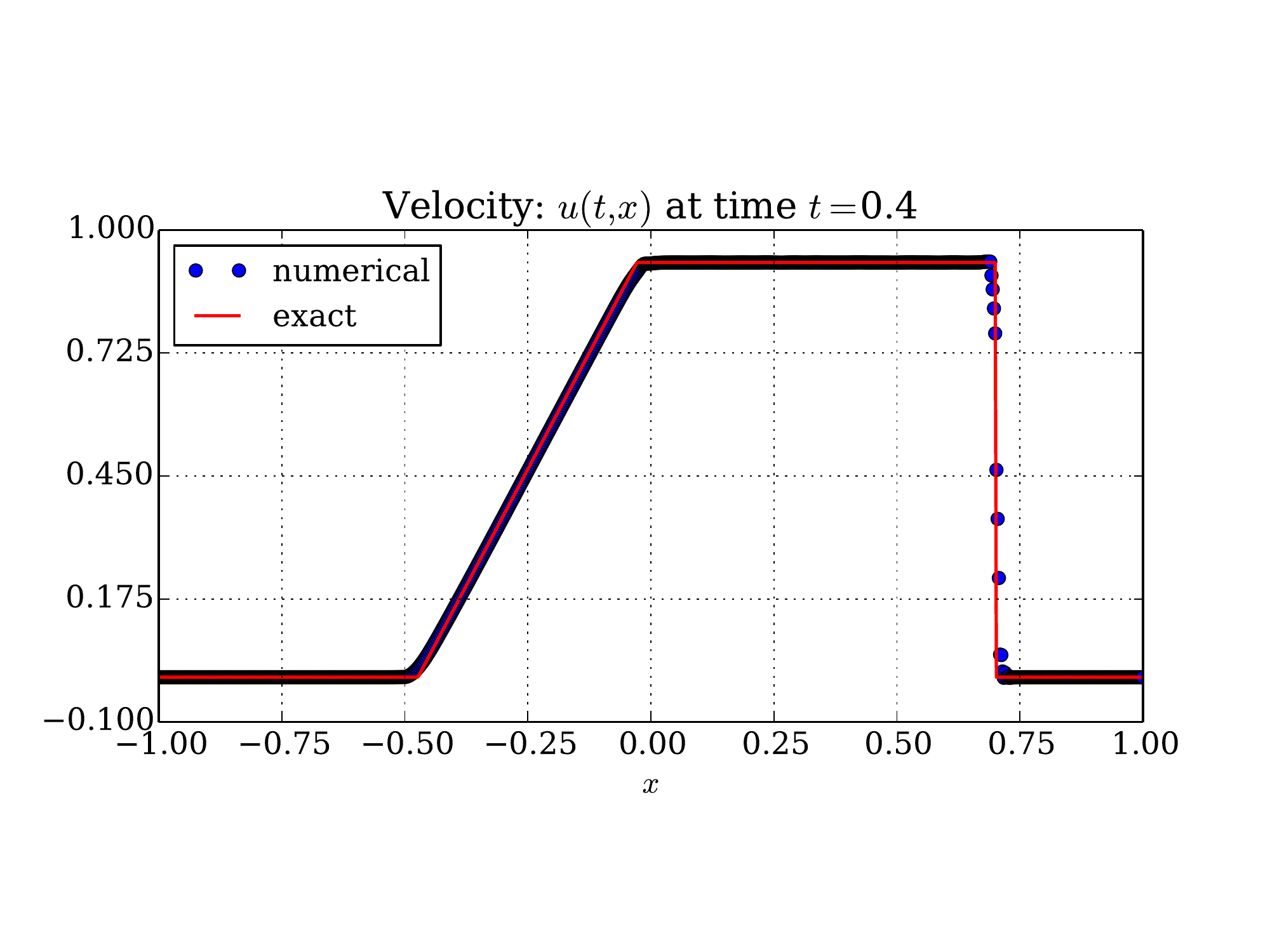} \\
\multicolumn{2}{c}{(c)\includegraphics[width=0.43\textwidth,trim={1.0cm 1.8cm 1.3cm 1.8cm},clip]{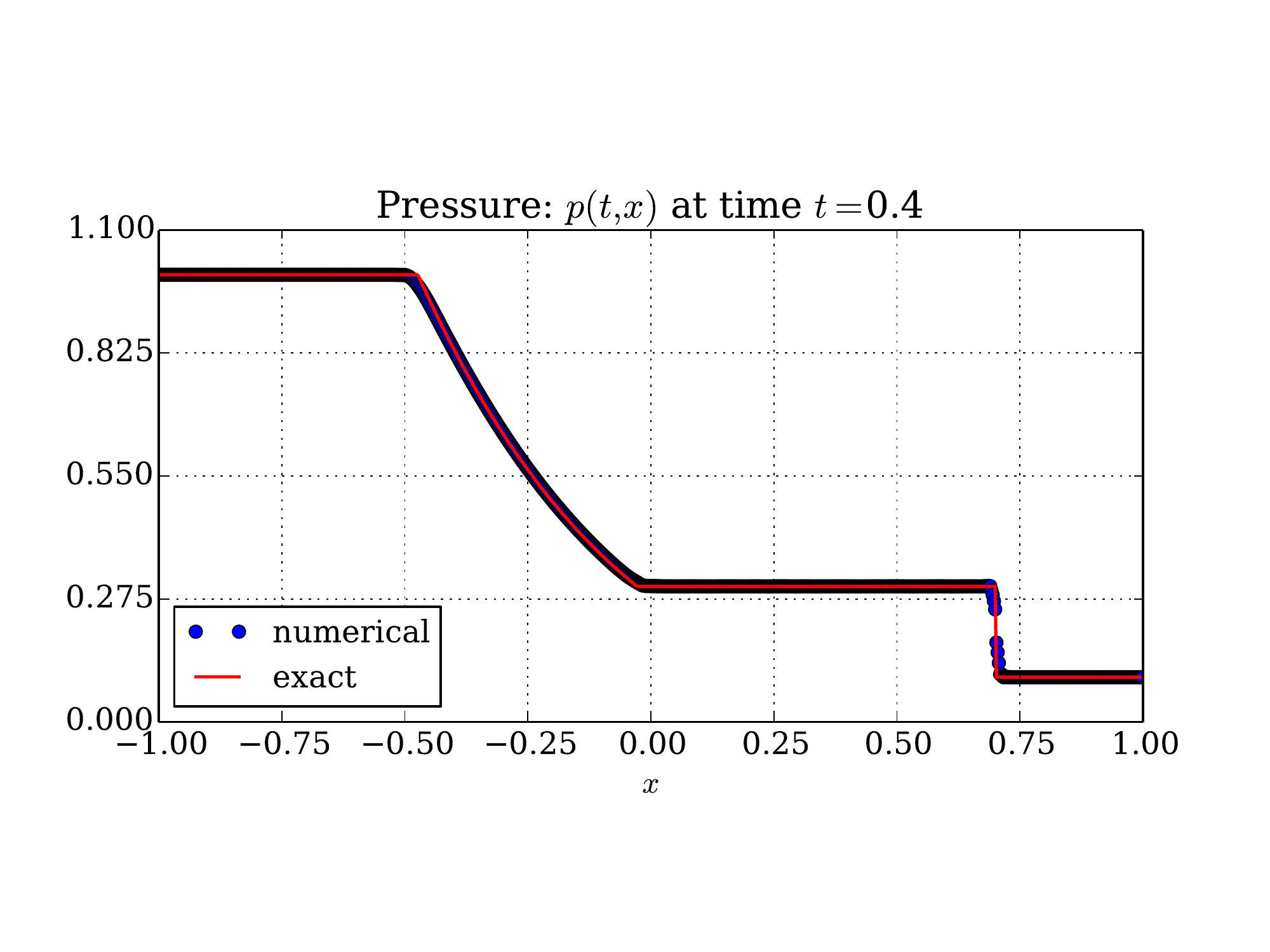}}
\end{tabular}
\caption{
Sod shocktube Riemann problem for the compressible Euler equations \cref{eqn:euler}. 
Shown are solutions computed with initial conditions
\cref{eqn:euler_shocktube}, outflow boundary conditions,
$\melems = 200$, and $\morder=4$ (i.e., fourth-order scheme with polynomial degree $\mdeg=3$). The solutions are shown at time $t=0.4$. 
In each panel we are plotting four points per element in order to clearly show
the subcell structure of the numerical solution. The individual panels show the numerical solution and the superimposed exact Riemann solution for the (a) density,
(b) fluid velocity, and (c) pressure. 
\label{fig:euler_shocktube}}
\end{figure}

The results of solving the shallow water equations with initial conditions \cref{eqn:euler_shocktube}, on
the domain $[-1,1]$ with outflow boundary conditions, with the 
$\morder=4$ version of the scheme, and with $\melems=200$ elements is shown in \cref{fig:euler_shocktube}.
The individual panels show at time $t=0.4$: (a) the density, (b) velocity, and (b) the pressure,
 with the exact Riemann solution superimposed.
In each panel we are plotting four points per element in order to clearly show
the subcell structure of the numerical solution.
The numerical solution is far from
exhibiting positivity violations; and therefore, the only limiter that is active is the 
correction step limiter (see \cref{subsec:nonoscillatory}). This example
demonstrates the ability of this once-per-time-step
limiter to simultaneously avoid producing unphysical oscillations and 
to maintain the sharpness of the solution.

\subsubsection{Double rarefaction}
In order to test the positivity limiters on the compressible Euler equations, we again 
attempt a Riemann problem that results in two counter-propagating rarefactions that leave in their wake a near-vacuum
state (see \cref{subsubsec:shllw_double_raref} for the shallow water version). 
This particular test case can be found in several papers, e.g., \cite{article:ZhangShu10rectmesh,article:ZhangShu11,article:seal2014explicit,article:MoRoSe17}.
The initial conditions are as follows:
\begin{equation}
\label{eqn:euler_double_raref}
\left( \rho, \, u, \, p \right)(t=0,x) = 
\begin{cases}
   \left( 7, \, -1, \, 0.2 \right) & \quad x<0, \\
   \left( 7, \, +1, \, 0.2 \right) & \quad x>0.
\end{cases}
\end{equation}
A full mathematical explanation of the solution of this Riemann problem can be found in several textbooks, including Chapter 14 of LeVeque.\cite{book:Le02}.

\begin{figure}[!t]
\centering
\begin{tabular}{cc}
(a)\includegraphics[width=0.43\textwidth,trim={1.0cm 1.8cm 1.3cm 1.8cm},clip]{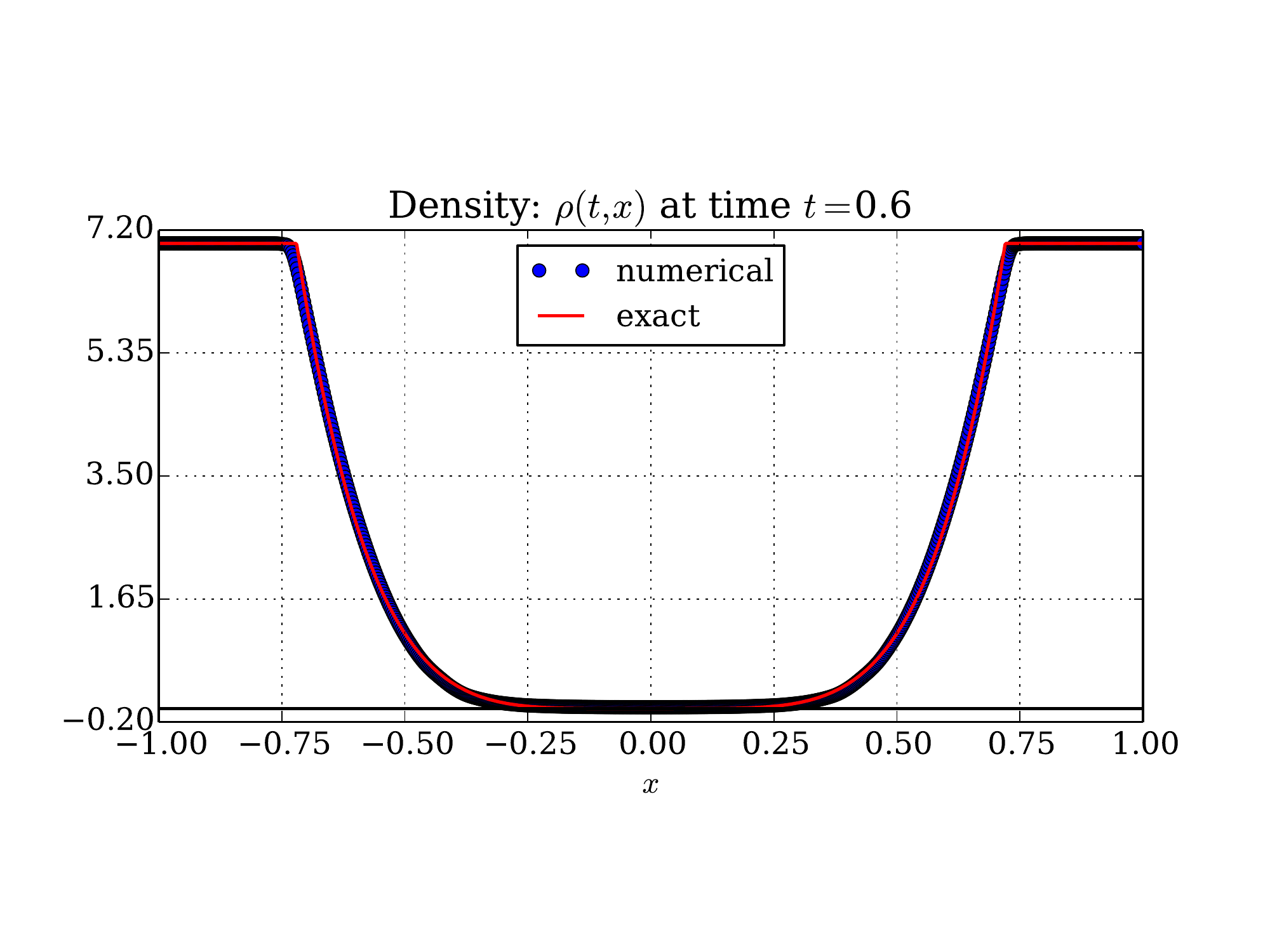} &
(b)\includegraphics[width=0.43\textwidth,trim={1.0cm 1.8cm 1.3cm 1.8cm},clip]{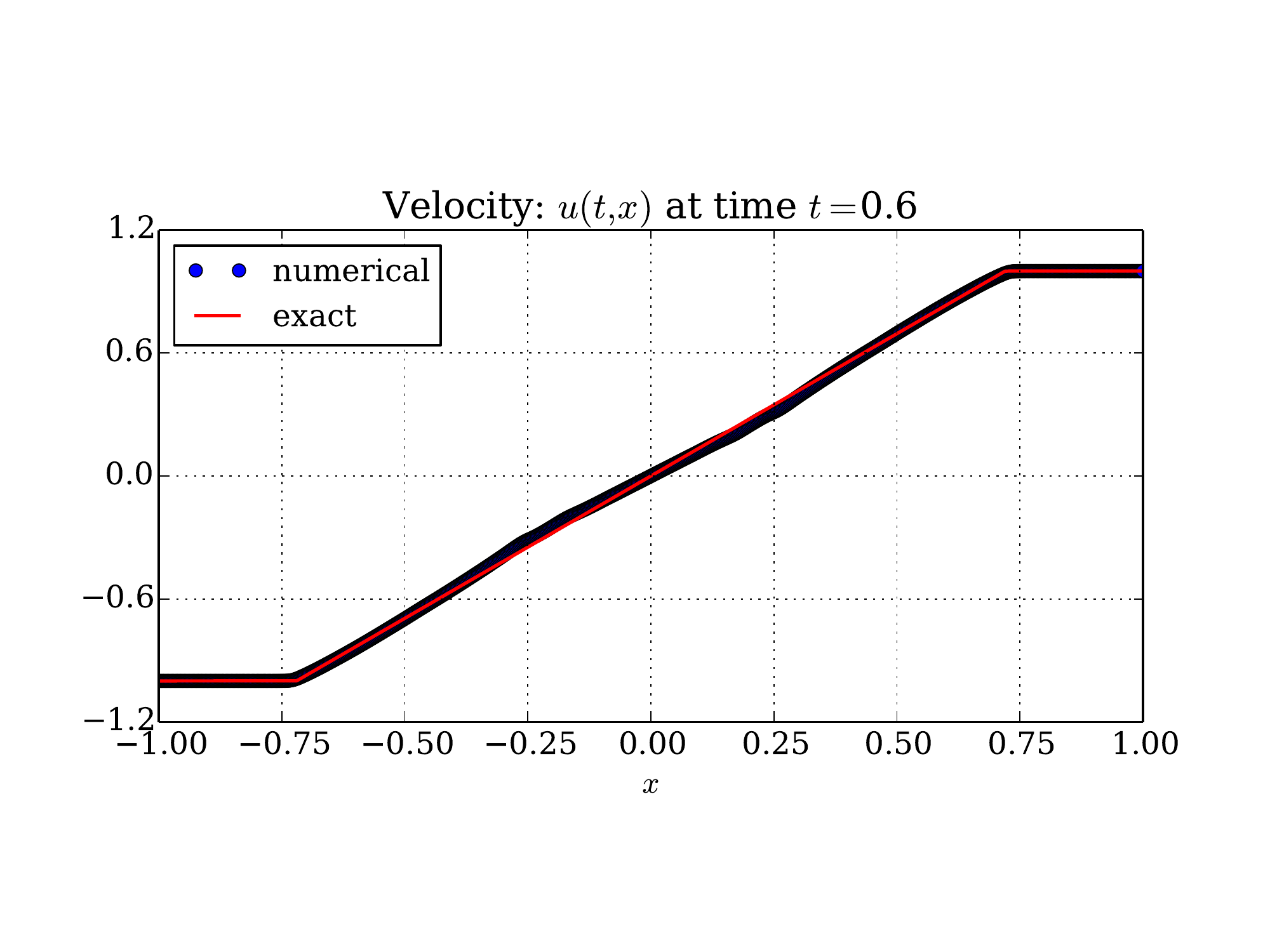} \\
\multicolumn{2}{c}{
(c)\includegraphics[width=0.43\textwidth,trim={1.0cm 1.8cm 1.3cm 1.8cm},clip]{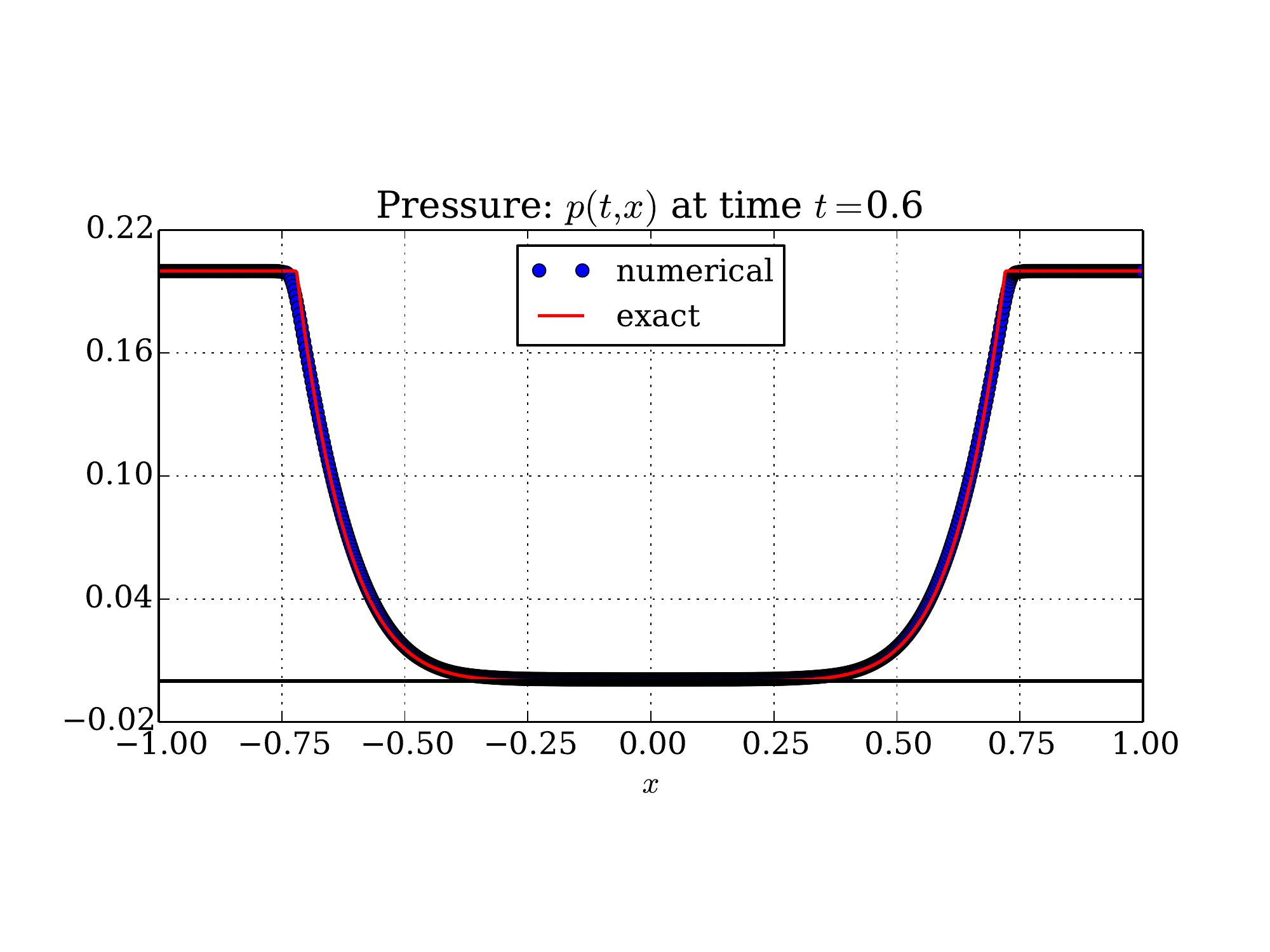}}
\end{tabular}
\caption{
Double rarefaction Riemann problem for the compressible Euler equations \cref{eqn:euler}. 
Shown are solutions computed with initial conditions
\cref{eqn:euler_double_raref}, outflow boundary conditions,
$\melems = 200$, and $\morder=4$ (i.e., fourth-order scheme with polynomial degree $\mdeg=3$). The solutions are shown at time $t=0.6$. 
In each panel we are plotting four points per element in order to clearly show
the subcell structure of the numerical solution. The individual panels show the numerical solution and the superimposed exact Riemann solution for the (a) density,
(b) fluid velocity, and (c) pressure. 
\label{fig:euler_doublerarefaction}}
\end{figure}

The results of solving the compressible Euler equations with initial conditions \cref{eqn:euler_double_raref}, on
the domain $[-1,1]$ with outflow boundary conditions, with the 
$\morder=4$ version of the scheme, and with $\melems=200$ elements is shown in \cref{fig:euler_doublerarefaction}.
The individual panels show the (a) density, (b) velocity, and (c) pressure at time $t=0.6$ with the exact
Riemann solution superimposed.
In each panel we are plotting four points per element in order to clearly show
the subcell structure of the numerical solution.
The proposed limiters are able to handle the near-vacuum solution in the center of the computational domain,
and the numerical solution remains stable and positivity-preserving in the sense of \cref{eqn:euler_pos_set_1}
and \cref{eqn:euler_pos_set_2}.

\subsubsection{Sedov blast problem}
Another standard test case for verifying the efficacy of positivity-preserving
limiters for the Euler equations is the 1D Sedov blast problem \cite{book:sedov1993similarity}.
For example, this test case is featured in the following papers:
\cite{article:ZhangShu10rectmesh,article:MoRoSe17,article:ZhangShu12}.
The initial conditions for the problem can be written on $[-1,1]$ as
\begin{equation}
\label{eqn:sedov_ic}
\rho(t=0,x) = 1, \quad
u(t=0,x) = 0, \quad
p(t=0,x) = \begin{cases}
(\gamma-1)\left(\frac{3.2 \times 10^6}{\Delta x} \right) & \bigl| x \bigr| \le \frac{\Delta x}{2}, \\
(\gamma-1) \left(10^{-12} \right) & \text{otherwise},
\end{cases}
\end{equation}
where $\Delta x$ is the mesh grid spacing. We assume here that there are an odd number of mesh elements
so that the middle element, $i=(\melems+1)/2$, is centered at the origin and defined by
\begin{equation}
{\mathcal T}_{(\melems+1)/2} = \left[-\frac{\Delta x}{2}, \frac{\Delta x}{2} \right].
\end{equation}
The initial conditions \cref{eqn:sedov_ic} represent a constant solution with almost zero
pressure everywhere, except in the middle element, where the pressure is many orders of magnitude larger than its surroundings. The resulting solution is a pressure blast that emenates
from the central element and propagates outward in both directions. 
These
initial conditions approximate a delta function of pressure.  Once the wave front propagates
away from the center of the domain, a post-shock region with near zero density is left behind.
This example can only be simulated with methods that either use extremely aggressive limiters or methods with guaranteed positivity-preservation.

\begin{figure}[!t]
\centering
\begin{tabular}{cc}
(a)\includegraphics[width=0.43\textwidth,trim={1.0cm 1.8cm 1.3cm 1.8cm},clip]{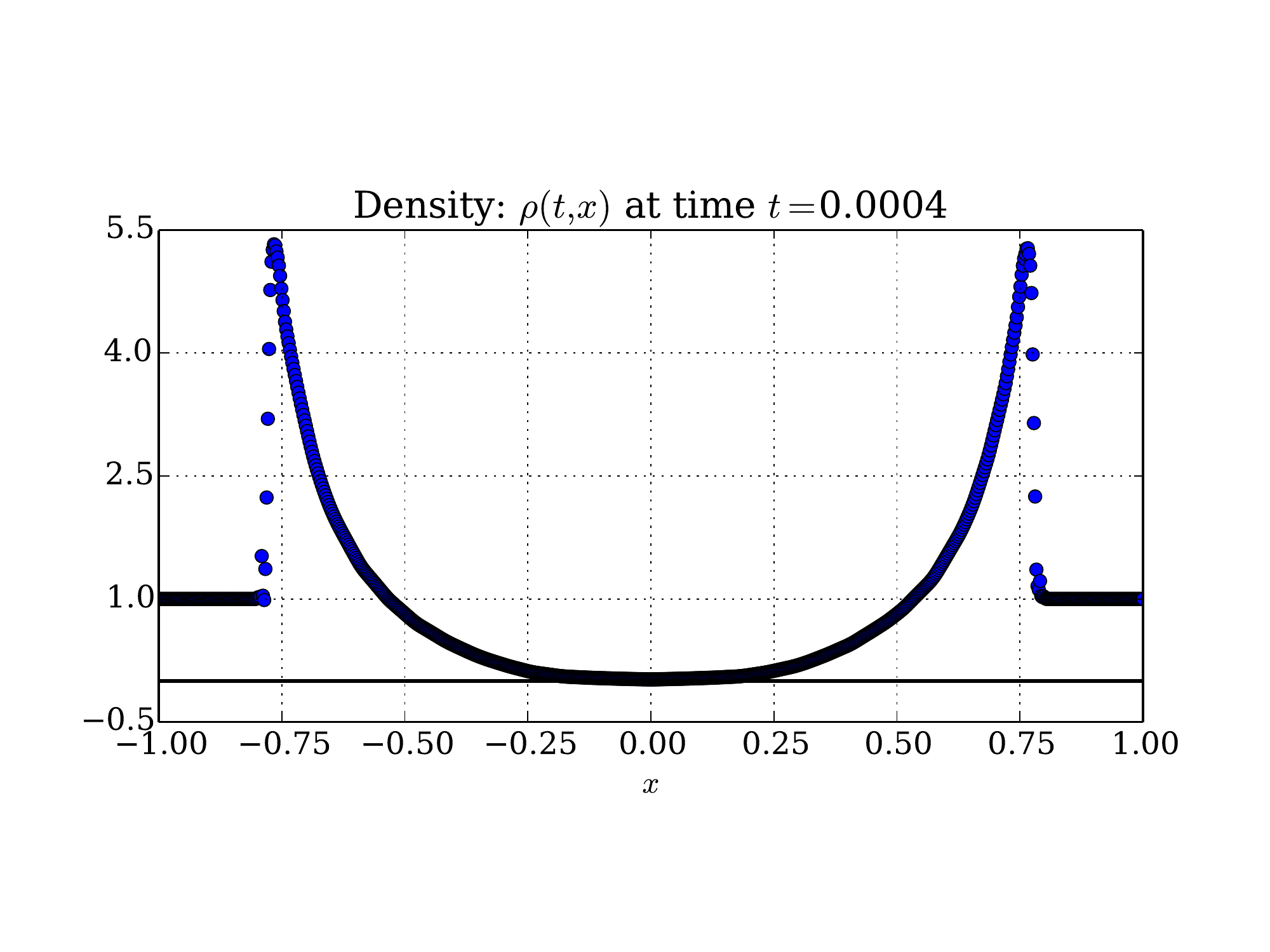} &
(b)\includegraphics[width=0.43\textwidth,trim={1.0cm 1.8cm 1.3cm 1.8cm},clip]{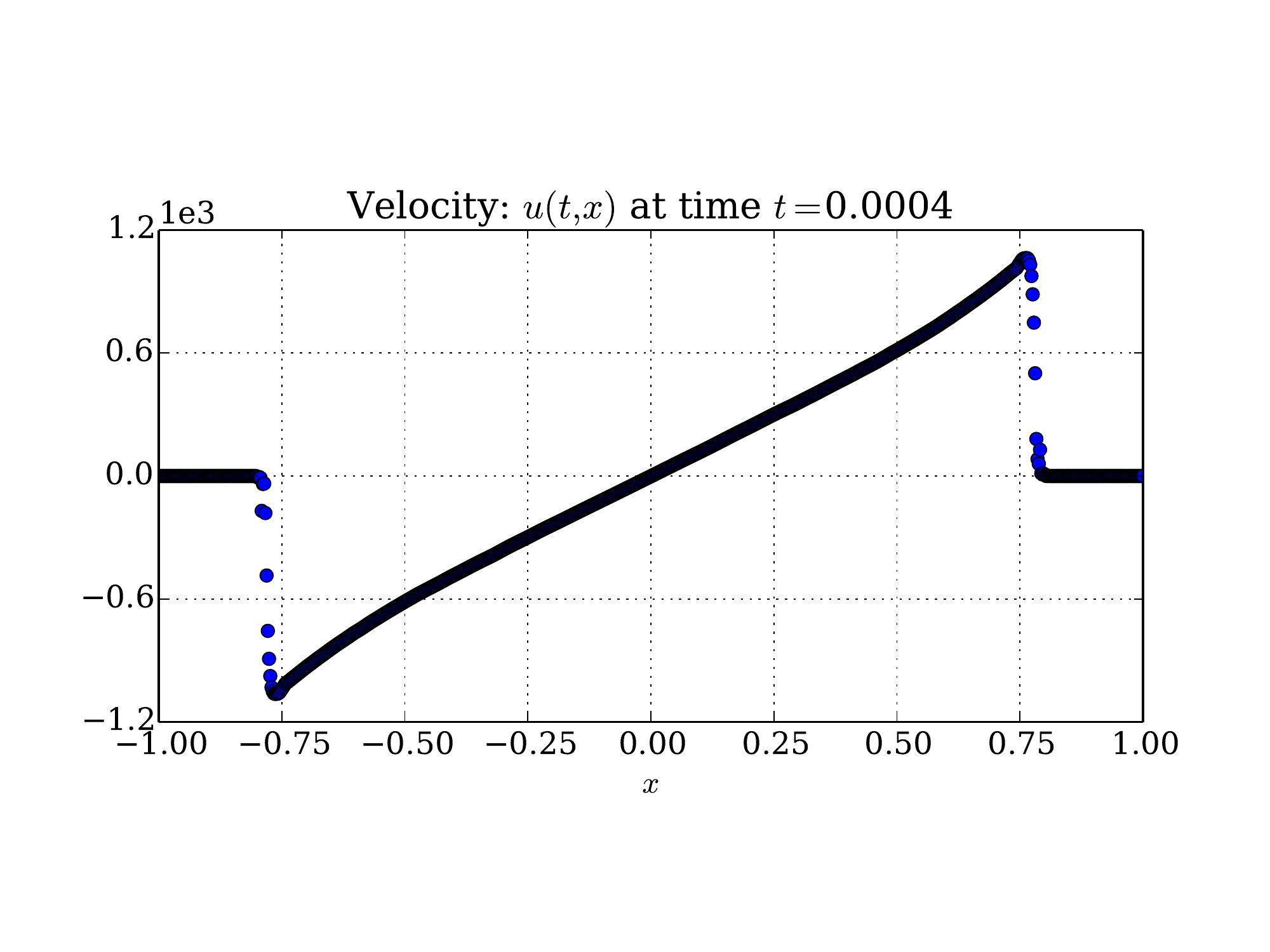} \\
\multicolumn{2}{c}{(c)\includegraphics[width=0.43\textwidth,trim={1.0cm 1.8cm 1.3cm 1.8cm},clip]{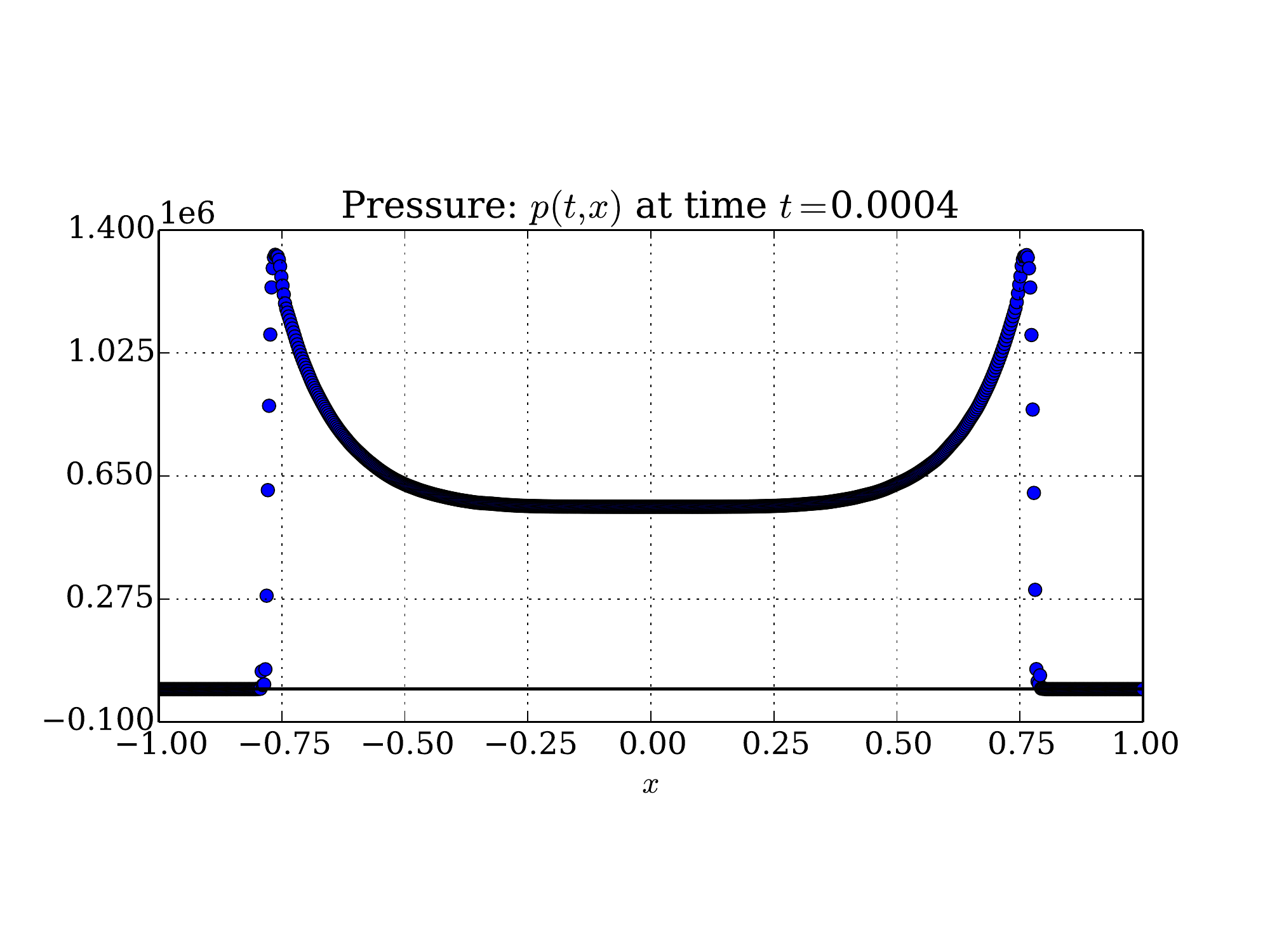}}
\end{tabular}
\caption{1D Sedov blast wave problem for the compressible Euler equations \cref{eqn:euler}. 
Shown are solutions computed with initial conditions
\cref{eqn:sedov_ic}, outflow boundary conditions,
$\melems = 201$, and $\morder=4$ (i.e., fourth-order scheme with polynomial degree $\mdeg=3$). The solutions are shown at time $t=0.0004$. 
In each panel we are plotting four points per element in order to clearly show
the subcell structure of the numerical solution. The individual panels show the (a) density,
(b) fluid velocity, and (c) pressure. 
\label{fig:euler_sedov}}
\end{figure}

The results of solving the compressible Euler equations with initial conditions \cref{eqn:sedov_ic}, on
the domain $[-1,1]$ with outflow boundary conditions, with the 
$\morder=4$ version of the scheme, and with $\melems=201$ elements is shown in \cref{fig:euler_sedov}.
The individual panels show the (a) density, (b) velocity, and (c) pressure at time $t=0.0004$.
In each panel we are plotting four points per element in order to clearly show
the subcell structure of the numerical solution.
The proposed limiters are able to handle the massive pressure jumps and the
near-vacuum solution in the center of the computational domain.
The numerical solution remains stable and positivity-preserving in the sense of \cref{eqn:euler_pos_set_1}
and \cref{eqn:euler_pos_set_2}.

\section{Conclusions}
\label{conclusions}
In this work we developed a new Lax-Wendroff discontinuous Galerkin (LxW-DG) method for solving hyperbolic conservation laws with a limiting strategy that keeps the solution non-oscillatory and positivity-preserving for relevant variables. For example, in the case of the shallow water equations, we guarantee positivity of the height, while in the case of the compressible Euler
equations, we guarantee positivity of the density and pressure.
The basic method was described in \cref{sec:LxW-DG}, while the various limiters were detailed
in \cref{sec:limiters}.

The scheme we developed is in the class of Lax-Wendroff DG schemes as introduced by
Qiu, Dumbser, and Shu \cite{article:Qiu05}, using the predictor-corrector 
interpretation developed by Gassner et al. \cite{article:GasDumHinMun2011}.
Each time-step of this new LxW-DG scheme is divided into two distinct phases:
\begin{description}
\item[{\bf Prediction step.}] In this phase, the equation and numerical solution are written in terms of primitive variables. A space-time DG approximation is applied on each element, but integration-by-parts is only performed on the time variable. This results in a system of nonlinear algebraic equations that are completely local on each element (i.e., no inter-element coupling). These nonlinear algebraic equations are approximately solved via a Jacobian-free Picard iteration, with the property that a sufficiently accurate solution is obtained after exactly $\morder$ iterations, where $\morder$ is the overall desired order of accuracy.
\item[{\bf Correction step.}] In this phase, the equation and numerical solution are written
in terms of conservative variables. A forward Euler-like step is applied to advance the solution from the old time, $t=t^n$, to the new time, $t=t^{n}+\Delta t$. This Euler-like step
is based on a DG scheme with proper integration-by-parts in the spatial variable, and requires the computation of temporal and spatiotemporal integrals of the predicted solution.
\end{description}
In order to guarantee positivity and to achieve numerical solutions without undue unphysical
oscillations, several limiters were introduced:
\begin{description}
\item[{\bf Prediction step positivity limiter (pointwise positivity).}]  Based on the celebrated Zhang and Shu \cite{article:ZhangShu11} limiter, we developed a completely local limiter that
minimally damps the high-order corrections to the primitive variables in order to get pointwise positivity of the predicted
solution at space-time quadrature points on each element. This limiter is applied once per Picard iteration for a total of $\morder$-times per time-step.
\item[{\bf Correction step positivity limiter I (positivity-in-the-mean).}] Following 
Moe et al. \cite{article:MoRoSe17}, we developed a limiter in which the high-order
numerical fluxes used to update the cell averages are minimally blended with a positivity-preserving low-order flux in such a way to obtain positivity of the high-order cell averages.
This limiter is applied once per time-step.
\item[{\bf Correction step positivity limiter II (pointwise positivity).}] Similar to what was done in the prediction step, a completely local limiter minimally damps the high-order corrections to the conserved variables in order to get pointwise positivity of the corrected
solution at spatial quadrature points on each element.
This limiter is applied once per time-step.
\item[{\bf Correction step unphysical oscillation limiter.}] Based on the  
 Krivodonova \cite{article:Kriv07} limited, we developed a hierarchical minmod limiter that
is applied to the characteristic variables in order to remove unphysical oscillations due to the
Gibbs phenomenon at shocks and rarefactions. 
This limiter is applied once per time-step.
\end{description}

The resulting Lax-Wendroff discontinuous Galerkin (LxW-DG) method was verified on a series
of standard test cases for the Burgers (\cref{subsec:burgers}), shallow water (\cref{subsec:shllw}), and compressible Euler (\cref{subsec:euler}) equations.
These test cases clearly showed that the overall scheme is successful at removing unphysical oscillations without overly diffusing the numerical solution, as well as keeping the solution
fully positivity-preserving (see \cref{sec:numerical_examples}). All of the presented methods and examples have been written in a freely available open-source Python code (see \cref{sec:pycode}).

\section*{Acknowledgments}
This research was carried out as part of the 2017 Summer REU (Research Experience for Undergraduates) Program at Iowa State University in Ames, Iowa. We would like to thank Iowa State University for their hospitality and the National Science Foundation for funding the REU program under the following grant: NSF Grant DMS--1457443.
Additionally, JAR was supported in part by NSF Grant DMS--1620128. We also thank the developers of Python
(\url{python.org}), Matplotlib (\url{matplotlib.org}), and PyPy (\url{pypy.org}) for developing the excellent open-source software tools that allow research projects such as this one to remain fully open-source.

\bibliographystyle{plain}

\end{document}